\documentclass[11pt,oneside]{amsart}
\usepackage{amsfonts}
\usepackage{hyperref}
\usepackage{etoolbox}
\usepackage{amsmath}
\usepackage{amssymb}
\usepackage{amsthm}
\usepackage{dsfont}
\usepackage{pifont}
\usepackage{mathrsfs}
\usepackage{graphicx}
\usepackage{enumerate}
\usepackage[all]{xy}
\usepackage{geometry}
\usepackage{amsfonts}
\usepackage{fancyhdr}
\usepackage{tikz}
\usepackage{tikz-cd}
\usetikzlibrary{arrows}
\usepackage[section]{algorithm}
\usepackage{algorithmic}

\usepackage[backend=biber,sorting=nty,maxnames=5]{biblatex}
\newcommand{\Aa}{\mathcal{A}}
\newcommand{\Bb}{\mathcal{B}}
\newcommand{\Cc}{\mathcal{C}}
\newcommand{\Dd}{\mathcal{D}}

\newcommand{\Ff}{\mathcal{F}}
\newcommand{\Gg}{\mathcal{G}}
\newcommand{\Hh}{\mathcal{H}}
\newcommand{\Ii}{\mathcal{I}}

\newcommand{\Mm}{\mathcal{M}}

\newcommand{\Oo}{\mathcal{O}}

\newcommand{\Rr}{\mathcal{R}}
\newcommand{\Ss}{\mathcal{S}}
\newcommand{\Tt}{\mathcal{T}}
\newcommand{\Uu}{\mathcal{U}}

\newcommand{\Xx}{\mathcal{X}}

\newcommand{\C}{\mathbb{C}}

\newcommand{\R}{\mathbb{R}}
\renewcommand{\S}{\mathfrak{S}}

\newcommand{\Z}{\mathbb{Z}}
\newcommand{\D}{\mathbb{D}}

\DeclareMathOperator\Ker{Ker}

\DeclareMathOperator\GL{GL}

\DeclareMathOperator\Id{Id}

\DeclareMathOperator\cyc{cyc}
\DeclareMathOperator\dec{dec}
\DeclareMathOperator\SSS{SSS}
\DeclareMathOperator\Ob{Ob}

\DeclareMathOperator\LL{LL}
\DeclareMathOperator\clbl{clbl}
\DeclareMathOperator\cc{cc}

\DeclareMathOperator\NCP{NCP}

\newcommand{\LLL}{\bbar{\LL}}
\renewcommand{\preceq}{\preccurlyeq}
\renewcommand{\succeq}{\succcurlyeq}

\newcommand{\xdownarrow}[1]{%
  {\left\downarrow\vbox to #1{}\right.\kern-\nulldelimiterspace}
}

\renewcommand{\matrix}[1]{\begin{pmatrix} #1\end{pmatrix}}

\newcommand{\infl}{\ar@{{(}->}}
\newcommand{\defl}{\ar@{->>}}

\newcommand{\intv}[1]{[\![#1]\!]}
\newcommand{\muls}[1]{ \left\{\kern-0.6em \left\{ #1\right\}\kern-0.6em \right\} }

\newcommand{\scal}[2]{\left\langle #1 ,#2 \right\rangle}
\newcommand{\epsi}{\varepsilon}

\newcommand{\pphi}{\varphi}

\newcommand{\spa}{\vspace*{1ex}}

\newcommand{\bbar}[1]{\overline{#1}}

\newcommand{\ttilde}[1]{\widetilde{#1}}
\newcommand{\nit}[1]{\textbf{\emph{#1}}}
\newcommand{\bbf}[1]{\mathbf{#1}}

\theoremstyle{plain}
\newtheorem{prop}{Proposition}[section]
\newtheorem{prop-def}[prop]{Proposition-Definition}
\newtheorem{lem}[prop]{Lemma}
\newtheorem{theo}[prop]{Theorem}
\newtheorem{cor}[prop]{Corollary}
\newtheorem{definition}[prop]{Definition}

\newtheorem*{theo*}{Theorem}
\newtheorem*{cor*}{Corollary}
\newtheorem*{prop*}{Proposition}

\newtheorem{theointro}{Theorem}
\newtheorem{propintro}[theointro]{Proposition}
\newtheorem{corintro}[theointro]{Corollary}

\theoremstyle{remark}
\newtheorem{rem}[prop]{Remark}
\newtheorem{exemple}[prop]{Example}
\newtheorem{notation}[prop]{Notation}

\title[Springer categories for regular centralizers in braid groups]{Springer categories for regular centralizers in well-generated complex braid groups}
\author{Owen Garnier}
\address{Departamento de \'Algebra e Instituto de Matem\'aticas de la Universidad de Sevilla, Spain.}
\email{owen.garnier@math.cnrs.fr}
\date{\today}


\begin{document}

\begin{abstract}In his proof of the $K(\pi,1)$ conjecture for complex reflection arrangements, Bessis defined  Garside categories suitable for studying braid groups of centralizers of Springer regular elements in well-generated complex reflection groups. We provide a detailed study of these categories, which we call Springer categories. 

We describe in particular the conjugacy of braided reflections of a regular centralizer in terms of the Garside structure of the associated Springer category. In so doing we obtain a pure Garside theoretic proof of a theorem of Digne, Marin and Michel on the center of finite index subgroups in complex braid groups in the case of a regular centralizer in a well-generated group. We also provide a ``Hurwitz-like'' presentation of Springer categories. To this aim we provide additional insights on noncrossing partitions in the infinite series. Lastly, we use this ``Hurwitz-like'' presentation, along with a generalized Reidemeister-Schreier method we introduce for groupoids, to deduce nice presentations of the complex braid group $B(G_{31})$.
\end{abstract}

\maketitle

\addtocontents{toc}{\protect\setcounter{tocdepth}{1}}
\tableofcontents

\section*{Introduction}
Let $W$ be an irreducible complex reflection group. In the case where $W$ is well-generated, that is, it can be generated by a number of reflections equal to its rank, Bessis introduced in \cite[Section 8]{beskpi1} a Garside monoid whose group of fractions is naturally isomorphic to the complex braid group $B(W)$ of $W$. It is the so-called dual braid monoid. 

Dual braid monoids have since been used to study group theoretic questions on complex braid groups. For instance, the center of finite index subgroups of complex braid groups is determined in \cite{dmm} through the use of various Garside monoids, including the dual braid monoid. More recently, Gonz\'alez-Meneses and Marin introduced in \cite{paratresses} a general concept of parabolic subgroups of a complex braid group. They also studied these subgroups using dual braid monoids in particular.

Unfortunately, the approach of using dual braid monoids is only available when $W$ is well-generated. This proves especially problematic for the exceptional group $G_{31}$, which is both badly-generated and of high rank. The lack of a well-studied Garside structure for the braid group $B(G_{31})$ was a serious obstacle in completing the study program of complex braid groups. For instance the proof of \cite[Theorem 1.4]{dmm} for the complex braid group $B(G_{31})$ is obtained by external representation-theoretic arguments using the linearity of spherical Artin groups.

However, Bessis also defined in \cite[Section 11]{beskpi1} a Garside category suitable for studying the braid group of the centralizer of a Springer regular element in a well-generated group. This approach applies in particular to $G_{31}$, which appears as the centralizer of a Springer regular element inside of the well-generated group $G_{37}$. The problem is then to understand these Garside categories. Indeed, the theory of Garside categories, which was in its infancy when Bessis wrote his article, is now much more developed, and allows for the generalization of many (if not all) properties of dual braid monoids to these categories. This article aims to be a first step in this program.

Let $W$ be a well-generated irreducible complex reflection group, and let $g\in W$ be a regular element in the sense of Springer. The centralizer $W_g:=C_W(g)$ is again a complex reflection group, and its braid group $B(W_g)$ is isomorphic to the centralizer in $B(W)$ of a so-called regular braid \cite[Theorem 12.4]{beskpi1}. The proof of this result by Bessis relies on heavy topological arguments, which boil down to the construction of a Garside category $\Cc$, whose enveloping groupoid $\Gg$ (obtained by formally inverting all morphisms) is equivalent to $B(W_g)$. Dual braid monoids can be seen as a particular case of this topological construction (associated with the identity element of $W$, which is $1$-regular). We refer to the category $\Cc$ (resp. to the groupoid $\Gg$) as the Springer category (resp. Springer groupoid) (associated with $W_g$ and $W$).

On the other hand, the construction of the Springer category can be formulated in a purely Garside theoretic way, as a category of periodic elements. We recall this construction in Section \ref{sec:garside}, before we restate the topological construction in Section \ref{sec:crg_bg_braid_cats}.

The interaction between these two constructions then allows us to describe the braided reflections of the group $B(W_g)$ in terms of the Garside structure of $\Cc$. In the classical case, the atoms of the dual braid monoid are known to be a set of braided reflections generating the braid group $B(W)$. The same result is not expected to hold in the general case, as an atom of the category $\Cc$ may have a different source and target, and thus it cannot be identified with an element of $B(W_g)$. In order to solve this issue, we define a family of endomorphisms in $\Cc$, which we call atomic loops: For every atom $a$ in $\Cc$, we define two morphisms $a^\flat$ and $a^\#$ which have the same target and same source as $a$, respectively. We write $a^{(2\#)}$ for $(a^{\#})^{\#}$ and so on. We show that there is a smallest integer $n$ such that $a^{(n\#)}=a^{\flat}$ and we define $\lambda(a)=aa^{\#}a^{(2\#)}\cdots a^{(n\#)}$. This element is by construction an endomorphism in $\Cc$, and we call it the atomic loop associated with $a$. Our first main result is that atomic loops do give the desired analogue of the braided reflections of $B(W_g)$.

\begin{theointro}[Theorem \ref{prop:braidreflectionconjugatetoatomicloop}] Let $u$ be an object of $\Cc$. Any atomic loop inside $\Cc(u,u)$ is a braided reflection of the group $B(W_g)\simeq\Gg(u,u)$. Conversely, any braided reflection $\sigma\in B(W_g)\simeq \Gg(u,u)$ is conjugate in the Springer groupoid to some atomic loop.
\end{theointro}

In the case where $\Cc$ is a monoid, the atomic loop associated with an atom is simply the atom itself, and we recover the classical situation. This explicit description of the braided reflections in terms of the Garside structure of $\Gg$ enables us to study the conjugacy of braided reflections and their powers using Garside theory. The following theorem ensures that the conjugacy of powers of braided reflections is ``the same as'' that of braided reflections.

\begin{theointro}[Theorem \ref{theo:conj_of_ato_loops}] Let $\lambda(s)\in \Cc(u,u)$ be an atomic loop of some object $u$, and let $f\in \Gg$. If there is some endomorphism $z\in \Cc$ such that $\lambda(s)^nf=fz$ for some $n\geqslant 1$, then $z=\lambda(s')^j$ for some atomic loop $\lambda(s')$ such that $\lambda(s)f=f\lambda(s')$.
\end{theointro}

This result can be seen as a generalization to Springer categories of \cite[Proposition 2.2]{dmm}, which was originally shown for a family of Garside monoids including dual braid monoids. As a corollary, we are able to complete the proof of a general result on centralizers of braided reflections in complex braid groups:

\begin{corintro}(\cite[Corollary 2.5]{dmm} and Corollary \ref{cor:centralizer_braided_reflection})\label{corintro:centralizer_braided_reflection}  Let $B$ be a complex braid group. The centralizer in $B$ of a braided reflection $\sigma$ is equal to the centralizer of any nontrivial power of $\sigma$.
\end{corintro}

Although this result is not stated explicitly in \cite{dmm}, it is actually proven during the proof of \cite[Corollary 2.5]{dmm} for braid groups of well-generated complex reflection groups (as well as for almost all groups of rank 2). We complete the proof for other badly-generated groups. Like in \cite{dmm}, we also obtain the analogue of \cite[Theorem 1.4]{dmm} in our context:

\begin{corintro}[Corollary \ref{cor:center_finite_index_subgroups}] Let $W$ be a well-generated irreducible complex reflection group, and let $g\in W$ be a regular element. If $U\subset B(W_g)$ is a finite index subgroup, then we have $Z(U)\subset Z(B(W_g))$.
\end{corintro}

Again, this last result was already proven in \cite{dmm} for all irreducible complex braid groups. However, the proof in the case of $B(G_{31})$ relied on the use of the generalized Krammer representation, while our proof is Garside theoretic. Moreover, our method actually proves the stronger statement of Corollary \ref{corintro:centralizer_braided_reflection}.

Our other point of focus in this article concerns presentations of Springer categories, and a way to deduce from such a presentation a presentation of the associated group $B(W_g)$. First, we have a presentation of the Springer category $\Cc$ associated with $B(W_g)$:

\begin{theointro}[Theorem \ref{prop:hurwitzpresentationforc31}] A Springer category is presented by its atoms, endowed with the Hurwitz relations, that is, all the relations induced by commutative squares of atoms.
\end{theointro}
Again, this is known to hold for dual braid monoids \cite[Lemma 8.8]{beskpi1}. The proof of this theorem relies on the particular case where the Springer category happens to be a monoid. In this case we show (Theorem \ref{theo:intervalles_identiques_p=1} and Corollary \ref{cor:dual_braid_relations_for_centralizer_of_c^q}) that the Springer category is naturally isomorphic to the dual braid monoid associated with $W_g$ (which is well-generated in this case). As a byproduct we obtain (Corollary \ref{cor:partitions_e1n}) an isomorphism between the lattice of simples of the dual braid monoid associated with $G(d,1,n)$ and the lattice of noncrossing partitions of type $(d,1,n)$ defined in \cite[Definition 1.11]{bescor}.

In the case of well-generated groups, following the indication of \cite[Remark 8.9]{beskpi1}, one can use the Hurwitz presentation of the dual braid monoid to obtain other presentations of the associated braid group, like the ones of \cite[Table 3]{bmr} or \cite[Section 2]{bessismichel}.

Following the idea of \cite[Remark 11.29]{beskpi1}, we propose in Section \ref{sec:presentation_of_b31} a generalization of this work to the case of categories, which we apply to the complex braid group $B(G_{31})$. The underlying method is an analogue of the Reidemeister-Schreier method adapted for groupoids. We define a Schreier transversal for a connected free groupoid $\Ff(\Ss)$ as a set of paths starting from a fixed object $u_0$ (the root of the transversal) and arriving at every object of $\Ff(\Ss)$. 

\begin{propintro}[Reidemeister-Schreier method for groupoids, Proposition \ref{prop:reidemeister_schreier}]~\newline 
Let $\Gg=\langle \Ss~|~ \Rr\rangle$ be a connected presented groupoid. A Schreier transversal $T$ for $\Ff(\Ss)$ rooted in the object $u_0$ induces an explicit presentation of the group $\Gg(u_0,u_0)$.
\end{propintro}

We consider the Springer groupoid $\Bb_{31}$ associated with the embedding of $G_{31}$ inside $G_{37}$ as the centralizer of a $i$-regular element. By choosing a particular Schreier transversal rooted in an object $u$, we are able to obtain a presentation of $\Bb_{31}(u,u)\simeq B(G_{31})$, where the generators are atomic loops. By considering a particular object of $\Bb_{31}$, we obtain the following theorem.
\begin{theointro}[Theorem \ref{theo:isomorphism_presentation_b31} and Section \ref{sec:4.3.1}]~\newline The complex braid group $B(G_{31})$ is given by the following presentation
\[B(G_{31})\simeq \left\langle s,t,u,v,w~\left| ~\begin{array}{l} ts=st,~vt=tv,~wv=vw,\\suw=uws=wsu,\\  svs=vsv,~vuv=uvu,~utu=tut,~twt=wtw. \end{array} \right. \right\rangle\]
where $s,t,u,v,w$ are braided reflections.
\end{theointro}

This is the presentation of $B(G_{31})$ that was conjectured in \cite[Table 3]{bmr} and \cite[Conjecture 2.4]{bessismichel}. Another approach to prove this presentation was proposed in \cite[Section 4]{beskpi1}, using the data of \cite[Figure 3]{bessismichel}. Unfortunately, this method cannot readily be carried out with these data as the formula for the discriminant of $G_{31}$ given by \cite[Figure 3]{bessismichel} does not appear to satisfy the conditions of \cite[Corollary 4.5]{beskpi1} (the discriminant considered in \cite[Corollary 4.5.(iii)]{beskpi1} is zero everywhere).

By considering other objects of $\Bb_{31}$, we obtain other new presentations of $B(G_{31})$. They are all positive homogeneous with braided reflections (atomic loops) as generators. They give explicit examples of presentations whose existence was proven in \cite[Theorem 0.1]{beszar} in the case of $G_{31}$. For each presentation, we give in particular an explicit isomorphism to the presentation given above in terms of images of the generators.

This work is part of my PhD thesis, with some results originating from my Master's thesis. Both are done under the supervision of Pr. Ivan Marin. I thank him for his precious advice during the preparation and redaction of this article.

\section{Preliminaries on Garside categories}\label{sec:garside}

Garside categories were originally introduced near the end of the 2000s as a natural generalization of Garside monoids (see for instance \cite{kgarcat} or \cite{besgar}). This early concept of Garside category was rapidly generalized into the more general concept of Garside family on a category. A comprehensive survey of the general theory of Garside families is given in \cite{ddgkm}. However, Garside categories as considered in \cite{besgar} are sufficient in the study of complex braid groups. We give in this section a quick summary of some results of \cite{ddgkm}, adapted to our context. Throughout this paper, all categories are assumed to be small categories.

\subsection{Definitions}\label{sec:1.1}In this section, we fix a category $\Cc$. The set of objects of $\Cc$ will be denoted by $\Ob(\Cc)$. For $u,v\in \Ob(\Cc)$, the set of morphisms from $u$ to $v$ in $\Cc$ will be denoted by $\Cc(u,v)$. We follow the usual convention for composition of arrows in Garside categories: the composition of the diagram 
\[\xymatrix{x\ar[r]^-f & y\ar[r]^-g & z}\]
will be denoted by $fg$. 

For $u\in \Ob(\Cc)$, we denote by $\Cc(u,-)$ (resp. $\Cc(-,u)$) the set of morphisms in $\Cc$ with source $u$ (resp. with target $u$). We define a relation $\preceq$ of \nit{left-divisibility} on $\Cc(u,-)$ by setting
\[\forall f,g\in \Cc(u,-), ~ f\preceq g\Leftrightarrow \exists h ~|~ fh=g.\]
In particular, the source of $h$ is the target of $f$, and the target of $h$ is the target of $g$. We say that $g$ is a right-multiple of $f$ and that $f$ left-divides $g$. Likewise, we define a relation $\succeq$ of \nit{right-divisibility} on $\Cc(-,u)$.

Recall that a \nit{length functor} on $\Cc$ is a functor $\ell:\Cc\to (\Z_{\geqslant 0},+)$ such that $\Cc$ is generated by morphisms with positive length. A category $\Cc$ endowed with a length functor $\ell$ is called \nit{homogeneous}. We will often abusively say that a category is homogeneous without referring to a specific length functor. A homogeneous category admits no nontrivial invertible morphism. Indeed if $f$ is invertible, then we have $\ell(f)+\ell(f^{-1})=\ell(1)=0$, so $\ell(f)=\ell(f^{-1})=0$. Since $\Cc$ is generated by elements of positive length this implies that $f$ is trivial.
A convenient way to prove that a category (or a monoid) is homogeneous is to define it by a homogeneous presentation (see Lemma \ref{lem:homogeneous_presentation}).

\begin{lem}
Let $\Cc=(\Cc,\ell)$ be a homogeneous category, and let $u\in \Ob(\Cc)$. The relations $\preceq$ and $\succeq$ are partial orders on $\Cc(u,-)$ and $\Cc(-,u)$, respectively.
\end{lem}
\begin{proof}
We prove the result for $\preceq$, the result for $\succeq$ is obtained by working in $\Cc^{\mathrm{op}}$. The relation $\preceq$ is obviously reflexive and transitive, it only remains to show that it is antisymmetric. Let $f,g\in \Cc(u,-)$ be such that $f\preceq g$ and $g\preceq f$. We have $\ell(f)\leqslant \ell(g)$ and $\ell(g)\leqslant \ell(f)$, so $\ell(g)=\ell(f)$. Let $h$ be such that $fh=g$. We have $\ell(h)=\ell(g)-\ell(f)=0$, thus $h$ is trivial and $f=g$.
\end{proof}

Let $\Cc$ be a category. A nontrivial element in $\Cc$ which admits no left-divisor (other than itself and the identity) is called an \nit{atom}. In a homogeneous category $(\Cc,\ell)$, a morphism $f$ can always be written as a composition of at most $\ell(f)$ atoms. 

\begin{definition}\cite[Definition 2.9]{ddgkm} 
Let $u$ be an object of a homogeneous category $\Cc$. 
\begin{itemize}
\item The \nit{left-gcd} of $f,g\in \Cc(u,-)$ is the meet of $f$ and $g$ in $(\Cc(u,-),\preceq)$ (should it exist), we denote it by $f\wedge g$. 
\item The \nit{right-lcm} of $f,g\in \Cc(-,u)$ is the join of $f$ and $g$ in $(\Cc(u,-),\succeq)$ (should it exist), we denote it by $f\vee g$.
\item Likewise, we define $f\wedge_R g$ and $f\vee_L g$ the \nit{right-gcd} and the \nit{left-lcm} of $f$ and $g$, respectively.
\end{itemize}
\end{definition}

Of course, gcds and lcms need not exist in $\Cc$. We need two more general definitions before we move on to the definition of a Garside category.

\begin{definition}\cite[Definition II.2.52]{ddgkm} Let $\Cc$ be a category. We say that $\Cc$ is \nit{cancellative} if every equality of the form $fgh=fg'h$ in $\Cc$ implies $g=g'$. This is equivalent to the statement that every morphism in $\Cc$ is both a monomorphism and an epimorphism.
\end{definition}

\begin{definition}\cite[Definition V.2.19]{ddgkm} Let $\Cc$ be a homogeneous cancellative category. A \nit{Garside map} in $\Cc$ is a map $\Delta:\Ob(\Cc)\to \Cc$ satisfying the following assumptions:
\begin{enumerate}[(1)]
\item For $u\in \Ob(\Cc)$, the source of $\Delta(u)$ is $u$. The target of $\Delta(u)$ is denoted by $\phi(u)$.
\item The families
\[\mathrm{Div}(\Delta):=\bigsqcup_{u\in \Ob(\Cc)} ([1_u,\Delta(u)],\preceq)~~\text{~~and~~}~~\mathrm{Div}_R(\Delta):=\bigsqcup_{u\in \Ob(\Cc)} ([1_{\phi(u)},\Delta(u)],\succeq)\]
are equal. We say that $\Delta$ is a \nit{balanced map}.
\item The family $\Ss:=\mathrm{Div}(\Delta)$ is finite and generates $\Cc$. We call its elements the \nit{simple morphisms}.
\item For every $g\in \Cc(u,-)$, the elements $g$ and $\Delta(u)$ admit a left-gcd.
\end{enumerate}
A homogeneous cancellative category $\Cc$ endowed with a Garside map $\Delta$ will be called a \nit{homogeneous Garside category}.
\end{definition}

\begin{rem}\label{1.6}
This definition is actually stronger than \cite[Definition V.2.19]{ddgkm}, we notably assume that $\Ss=\mathrm{Div}(\Delta)$ is finite. This implies in particular that $\Cc$ must have a finite number of objects. On the other hand, \cite[Definition V.2.19]{ddgkm} requires the map $u\mapsto \phi(u)$ to be injective. In our case, condition $(4)$ implies that the map $u\mapsto \phi(u)$ is surjective. Since $\Cc$ has finitely many objects, it is then automatically injective.
\end{rem}

When working in a Garside category $(\Cc,\Delta)$, we will often replace the assertion ``$f=\Delta(u)$ where $u$ is the source of $f$'' with ``$f\in \Delta$'' or even ``$f=\Delta$'' to avoid cumbersome expressions. This causes no confusion because $f=\Delta(u)$ implies that $u$ is the source of $f$, thus $\Delta(u)$ is the only morphism of the form $\Delta(x)$ to which $f$ may be equal.

\begin{notation}
Let $s$ be a simple morphism in a Garside category. If $u$ denotes the source of $s$, there is by definition a morphism $s'$ such that $ss'=\Delta(u)$. By cancellativity, the morphism $s'$ is unique and we denote it by $\bbar{s}$. Similarly, since $\Delta$ is a balanced map, there is a unique morphism $s^*$ such that $s^*s=\Delta(\phi^{-1}(v))$, where $v$ is the target of $s$. Note that we used that the map $u\mapsto \phi(u)$ is bijective in order to define $s^*$.
\end{notation}

\begin{prop}\cite[Proposition V.1.28 and Proposition V.2.32]{ddgkm}\label{prop:garside_automorphisme}\newline Let $(\Cc,\Delta)$ be a homogeneous Garside category. We define an automorphism $\phi$ of $\Cc$ by setting
\begin{itemize}
\item for $u\in \Ob(\Cc)$, $\phi(u)$ is the target of $\Delta(u)$.
\item for $f\in \Cc(u,v)$, $\phi(f)$ is the unique morphism in $\Cc$ such that $f\Delta(v)=\Delta(u)\phi(f)$.
\end{itemize}
Furthermore, $\phi$ has finite order and $\Delta$ is a natural transformation from the identity functor $1_\Cc$ of $\Cc$ to $\phi$.
\end{prop}
Let $s$ be a simple morphism. We have $\phi(s)=\bbar{\bbar{s}}$ because $s\Delta=(s\bar{s}) \bbar{\bbar{s}}=\Delta\bbar{\bbar{s}}$. In particular we see that $\phi$ maps $\Ss$ into itself. We also have that $\phi$ is uniquely determined by the property that, for all $s\in \Ss$, we have $\phi(s)=\bbar{\bbar{s}}$.

\begin{prop}\cite[Proposition V.2.35]{ddgkm} Let $(\Cc,\Delta)$ be a homogeneous Garside category, and let $u\in \Ob(\Cc)$. The posets $(\Cc(u,-),\preceq)$ and $(\Cc(-,u),\succeq)$ are lattices. The posets $(\Ss(u,-),\preceq)$ and $(\Ss(-,u),\preceq)$ are sublattices of these lattices.
\end{prop}

\begin{rem}\label{rem:garside_monoids}
If we consider the particular case of monoids, that is categories with only one object, we recover the classical definition of a homogeneous Garside monoid. The relations $\preceq$ and $\succeq$ are the classical left- and right-divisibility relations. The Garside map $\Delta$ corresponds to an element of the monoid, which is the Garside element.
\end{rem}

Let $(\Cc,\Delta)$ be a homogeneous Garside category with set of simples $\Ss$, and let $\phi$ be the automorphism of $\Cc$ introduced in Proposition \ref{prop:garside_automorphisme}. In the sequel, we will often be interested in subcategories of fixed points under some power of $\phi$. Let $q\geqslant 0$ be an integer, and let $\Cc^{\phi^q}$ be the subcategory of $\Cc$ consisting of $\phi^q$-invariant morphisms. We also consider the set $\Ss^{\phi^q}$ of simple morphisms which are $\phi^q$-invariant.

\begin{lem}\label{lem:deux_sur_trois} Let $a,b\in \Cc$. If two of $a,b$ and $ab$ lie in $\Cc^{\phi^q}$, then so does the third.\end{lem}
\begin{proof}
If both $a$ and $b$ are $\phi^q$-invariant, then $\phi^q(ab)=\phi^q(a)\phi^q(b)=ab$. If both $a$ and $ab$ are $\phi^q$-invariant, then $\phi^q(a)\phi^q(b)=\phi^q(ab)=ab=\phi^q(a)b$, and $\phi^q(b)=b$ by cancellativity. The third case is dual to the second.
\end{proof} 

\begin{lem}\label{lem:phi_q_invar_sature} Let $s\in \Ss$ and suppose that the source of $s$ is $\phi^q$-invariant. We denote by $\psi(s)$ the right-lcm in $\Ss$ of all the $\phi^{iq}(s)$ for $i\in \Z_{\geqslant 0}$. We have $\psi(s)\in \Ss^{\phi^q}$, and for every $y\in \Ss^{\phi^q}$, we have $s\preceq y$ in $\Ss$ if and only if $\psi(s)\preceq y$ in $\Ss^{\phi^q}$.
\end{lem}
\begin{proof}
Let $u$ be the source of $s$, and let $n$ denote the order of $\phi^q$. We have $\psi(s)=s\vee \phi^q(s)\vee\cdots \vee \phi^{(n-1)q}(s)$. As $\phi^q$ induces an automorphism of the lattice $\Ss(u,-)$, we have
\[\phi^q(\psi(s))=\phi^q(s)\vee \phi^{2q}(s)\vee\cdots \vee \phi^{nq}(s)=\psi(s)\]
and $\psi(s)\in \Ss^{\phi^q}$. Let now $y\in \Ss^{\phi^q}$, if $\psi(s)\preceq y$, then $s\preceq \psi(s)\preceq y$. Conversely, if $x\preceq y$, then for all $i\in \intv{1,n-1}$, we have $\phi^{iq}(s)\preceq \phi^{iq}(y)=y$, thus $\psi(s)\preceq y$ by definition of the right-lcm.
\end{proof}

By \cite[Proposition VII.4.2]{ddgkm}, the category $\Cc^{\phi^q}$, endowed with the restriction of the map $\Delta$, is again a homogeneous Garside category. Its simple morphisms are the elements of $\Ss^{\phi^q}$. In particular we see that $(\Ss^{\phi^q}(u,-),\preceq)$ is always a lattice.

\subsection{Normal forms, groupoid and conjugacy}

In this section, we fix a homogeneous Garside category $(\Cc,\Delta)$, with set of simples $\Ss$, and Garside automorphism $\phi$. Recall from \cite[Definition II.1.28]{ddgkm} that an $\Ss$-path is a sequence of composable elements of $\Ss$ in $\Cc$. By definition of a homogeneous Garside category, every morphism in $\Cc$ can be expressed by an $\Ss$-path. 

When working with Garside categories, almost every path we will consider will actually be an $\Ss$-path. Thus unless specified otherwise, we will simply say ``path'' when meaning ``$\Ss$-path''. When no confusion is possible, we will write a path as a sequence of (composable) elements of $\Ss$.

\begin{definition}\cite[Corollary V.1.54]{ddgkm}\label{def:greediness}
A path $st$ of length $2$ in $\Cc$ is called \nit{greedy} if $s$ is the left-gcd of $st$ with $\Delta$, or equivalently, if $\bbar{s}$ and $t$ are left-coprime. In general, a path $s_1\cdots s_r$ is called \nit{greedy} if each subpath $s_is_{i+1}$ is greedy for $i\in \intv{1,r-1}$, and if $s_{r}$ is nontrivial.
\end{definition}

\begin{rem}
Since $\phi$ is an automorphism of $\Cc$, it preserves left-gcds and thus, it preserves greediness. In other words, if $s_1\cdots s_r$ is a path in $\Cc$, then it is greedy if and only if the path $\phi(s_1)\cdots \phi(s_r)$ is greedy.
\end{rem}

\begin{prop}\cite[Proposition V.3.4]{ddgkm} Every morphism $f$ in $\Cc$ can be expressed by a unique greedy path. This path is called the \nit{greedy normal form} of $f$. 
\end{prop}

\begin{lem}\label{lem:greedylength2}
Let $s$ and $t$ be two composable simple morphisms in $\Cc$, and let $d:=\bbar{s}\wedge t$. We define $s':=sd\in \Ss$ and we set $t'$ to be the unique simple morphism such that $dt'=t$. We have $st=s't'$ in $\Cc$, and the path $s't'$ is greedy.
\end{lem}
\begin{proof}
First, we have $st=sdt'=s't'$. Then, let $p$ be such that $dp=\bbar{s}$. We have $sdp=\Delta$, thus $p=\bbar{s'}$ and
\[d(\bbar{s'}\wedge t')=d\bbar{s'}\wedge dt'=s\wedge t=d.\]
Thus $\bbar{s'}\wedge t'$ is trivial by cancellativity, and $s't'$ is a greedy path.
\end{proof}

This lemma provides an algorithmic way to compute the greedy normal form of a morphism in $\Cc$, provided that we know how to compute the left-gcd of two simple morphisms.

For every category $\Cc$, one can consider the \nit{enveloping groupoid} $\Gg(\Cc)$ of $\Cc$, defined by formally inverting all morphisms in $\Cc$ (see Lemma \ref{lem:enveloping_groupoid}). In the case where $(\Cc,\Delta)$ is a homogeneous Garside category, the natural functor $\Cc\to \Gg(\Cc)$ is an embedding, and the groupoid $\Gg(\Cc)$ can be described as a \nit{groupoid of fractions}. We call $\Gg(\Cc)$ a \nit{Garside groupoid}.

\begin{notation}
For a positive integer $m$, $\Delta^m(u)$ denotes the path $\Delta(u)\cdots \Delta(\phi^{m-1}(u))$. For a negative integer $m$, $\Delta^m(u)$ denotes the inverse in $\Gg(\Cc)$ of $\Delta^{-m}(\phi^{-m}(u))$. Occasionally, we will write $\Delta^m$ instead of $\Delta^m(u)$ when there is no need to specify the source $u$ explicitly, following \cite[Convention 3.7]{ddgkm}.
\end{notation}

\begin{prop-def}\cite[Definition V.3.17 and Proposition V.3.18]{ddgkm} Let $f$ be a morphism in $\Gg(\Cc)$. There is a unique way to express $f$ as a path of the form $f=\Delta^p s_1\cdots s_r$ such that $p \in \Z$, $s_1\neq \Delta$ and $s_1\cdots s_r$ is a greedy path. The path $\Delta^{p}s_1\cdots s_r$ in $\Gg(\Cc)$ is called the \nit{left-weighted factorization} of $f$.
\end{prop-def}

In particular we see that, for every morphism $f$ in $\Gg(\Cc)$, there is a positive integer $m$ such that $\Delta^m f$ lies in $\Cc$.

\begin{definition}\cite[Definition V.3.23]{ddgkm} Let $f$ be a morphism in $\Gg(\Cc)$, with left-weighted factorization $\Delta^ps_1\cdots s_r$. The \nit{infimum} and \nit{supremum} of $f$ are defined by $\inf(f):=p$ and $\sup(f):=p+r$, respectively.
\end{definition}

One of the original purposes of Garside groups was to provide an algorithmic solution to the conjugacy problem. The general idea is to define a characteristic and computable subset of the conjugacy class of a given element, such that two elements are conjugate if and only if their associated sets coincide. One classical such set is the so-called super-summit set. It translates well to the groupoid context.

\begin{definition}\cite[Definition VIII.1.1]{ddgkm} Let $x$ and $x'$ be two endomorphisms in $\Gg(\Cc)$. A morphism $f\in \Gg(\Cc)$ \nit{conjugates} $x$ to $x'$ if $xf=fx'$ in $\Gg(\Cc)$. The endomorphism $x'$ will be denoted by $x^f$. As usual, conjugacy in $\Gg(\Cc)$ induces an equivalence relation on endomorphisms, and we can consider the \nit{conjugacy class} in $\Gg(\Cc)$ of the endomorphism $x$.
\end{definition}

For the remainder of this section, we fix an endomorphism $x$ in the Garside groupoid $\Gg(\Cc)$. We also fix $\Delta^p s_1\cdots s_r$ to be the left-weighted factorization of $x$.

\begin{definition}\cite[Definition VIII.2.3 and Definition VIII.2.8]{ddgkm}\label{def:cycling} The \nit{initial factor} (resp. \nit{final factor}) of $x$ is defined as $\phi^{-p}(s_1)$ (resp. $s_r$).\newline The \nit{cycling}  of $x$ is defined as $\cyc(x):=x^{\phi^{-p}(s_1)}=\Delta^ps_2\cdots s_r \phi^{-p}(s_1).$\newline The \nit{decycling} of $x$ is defined as $\dec(x)=x^{s_r^{-1}}=s_r\Delta^ps_1\ldots s_{r-1}=\Delta^p\phi^{p}(s_r)s_1\ldots s_{r-1}.$
\end{definition}

Note that the expressions given for $\cyc(x)$ and $\dec(x)$ are not necessarily left-weighted factorizations.

\begin{definition}\cite[Definition 3.1]{rigidity}\label{def:rigid} The morphism $x$ is called \nit{rigid} if the path $s_r\phi^{-p}(s_1)$ is greedy, or if $r=0$.
\end{definition}

\begin{lem}\label{lem:cyclinganddecyclingofrigid}
If $x$ is rigid, then the left-weighted factorizations of $\cyc(x)$ and $\dec(x)$ are given by
\[\cyc(x)=\Delta^ps_2\cdots s_r \phi^{-p}(s_1) \text{~and~} \dec(x)=\Delta^p\phi^{p}(s_r)s_1\ldots s_{r-1}.\]
Furthermore, the morphisms $\cyc(x)$ and $\dec(x)$ are also rigid.
\end{lem}

\begin{proof}
By definition, the path $s_2\cdots s_r \phi^{-p}(s_1)$ is greedy, and $s_2\neq \Delta$ (otherwise the path $s_1s_2$ would not be greedy). Thus the path given for $\cyc(x)$ is its left-weighted factorization. The same reasoning holds for $\dec(x)$. Lastly, the fact that $s_1\cdots s_r$ is a greedy path gives directly that $\cyc(x)$ and $\dec(x)$ are both rigid. 
\end{proof}

\begin{prop-def}\cite[Definition VIII.2.12]{ddgkm}\label{propdef:sss}
The conjugacy class of $x$ in $\Gg(\Cc)$ admits a well-defined subset $\SSS(x)$  on which each one of $\inf$ and $\sup$ takes a constant value. Furthermore, for every conjugate $x'$ of $x$ in $\Gg(\Cc)$, we have
\[\inf(x')\leqslant \inf(\SSS(x))\text{~and~} \sup(x')\geqslant \sup(\SSS(x)).\]
The set $\SSS(x)$ is called the \nit{super-summit set} of $x$.
\end{prop-def}

Note that $\SSS(x)$ must be finite, because there are finitely many morphisms in $\Gg(\Cc)$ with given $\inf$ and $\sup$. Also, if $x\in \Cc$, then $\inf(x)\geqslant 0$, and $\SSS(x)$ is included in $\Cc$.

\begin{prop}\cite[Proposition VII.2.16]{ddgkm}\label{prop:sssreachedbycycling} One can go from $x$ to an element of $\SSS(x)$ by a finite sequence of cycling, followed by a finite sequence of decycling.
\end{prop}

This proposition, combined with Lemma \ref{lem:cyclinganddecyclingofrigid}, shows that a rigid element always lies inside its own super-summit set. 

\begin{prop}\cite[Lemma VIII.2.19 and Proposition VIII.2.20]{ddgkm}\label{prop:sssconnectedbysimples} Let $x',x''$ be in $\SSS(x)$. Let also $f$ be a morphism in $\Cc$ with $x'^f=x''$. If $f=s_1\ldots s_r$ is the greedy normal form of $f$, then for all $i\in \intv{1,r}$ the morphism $x'^{s_1\cdots s_i}$ lies in $\SSS(x)$.
\end{prop}

Notice that, as $\phi$ preserves left-weighted factorizations, it stabilizes super-summit sets. In particular, the last proposition also applies to the case where the conjugating element lies in $\Gg(\Cc)$.

\subsection{Periodic elements and divided categories}
In this section, we fix a homogeneous Garside category $(\Cc,\Delta)$, with set of simples $\Ss$ and Garside automorphism $\phi$. We denote by $\Gg$ the enveloping groupoid of $\Cc$. We also fix two coprime integers $p,q>0$, along with positive integers $\eta,\mu$ such that $p\eta-q\mu=1$.

The study of super-summit sets provides a solution to the conjugacy problem in $\Gg$ for any pair of endomorphisms. If we restrict our attention to the so-called periodic elements of $\Gg$, then we obtain a more convenient solution to the conjugacy problem through the construction of a particular Garside category.

\begin{definition}\cite[Definition V.3.2]{ddgkm}\label{def:periodic} \\An endomorphism $\rho$ of $\Gg(u,u)$ is called $(p,q)$-\nit{periodic} if $\rho^p=\Delta^q(u)$.
\end{definition} 

In particular, if $\Gg(u,u)$ contains $(p,q)$-periodic elements, then $\Delta^q$ is an endomorphism and $\phi^q(u)=u$.

\begin{rem}\label{rem:pqcoprime}
The above definition makes sense for arbitrary integers $p,q$. However for a nonzero integer $n$, any $(np,nq)$-periodic element in $\Gg$ is conjugate in $\Gg$ to a $(p,q)$-periodic element lying in $\Cc$ by \cite[Corollary VIII.3.31 and Proposition VIII.3.34]{ddgkm}. Thus, understanding the conjugacy of periodic elements with coprime parameters is enough to understand the conjugacy of arbitrary periodic elements.
\end{rem}

The study of periodic elements in a general Garside context was introduced in \cite{besgar}, in which the concept of divided category is introduced. The $p$-divided category of a (homogeneous) Garside category $(\Cc,\Delta)$ is a (homogeneous) Garside category $(\Cc_p,\Delta_p)$, with Garside automorphism $\phi_p$. The conjugacy of $(p,q)$-periodic elements in $\Gg$ can be understood by considering the category $(\Cc_p)^{\phi_p^q}$ of $(\phi_p)^q$-invariant morphisms in the divided category $\Cc_p$ (see also \cite[Section XIV.1.1]{ddgkm}). In this section, we will provide a more direct version of this construction. 

First, let us follow the exposition given in \cite[Section XIV.1.1]{ddgkm} about divided categories. For $m\in \Z_{\geqslant 1}$, define
\[D_m(\Delta):=\{u:=(u_0,\ldots,u_{m-1})\in \Ss^m~|~ u_0\cdots u_{m-1}=\Delta\}.\]
In particular, we require the $u_i$ to be composable: the target of $u_i$ must be the source of $u_{i+1}$ for $i\in \intv{0,m-2}$. For every positive integer $m$, we define an action of the automorphism $\phi$ on $m$-tuples by setting
\[\phi.(u_0,\ldots,u_{m-1}):=(u_1,\ldots,u_{m-1},\phi(u_0))\in D_m(\Delta).\]
This is not formally an action of the group $\langle \phi\rangle$, as $\phi^k=1_\Cc$ does not necessarily mean that $\phi^k.u=u$ for all $p$-tuple $u$. This is rather an action of the free group $\Z$, which we denote by $\phi$ for convenience. Note that this action restricts to an action on $D_m(\Delta)$. Indeed, for $u=(u_0,\ldots,u_{m-1})\in D_m(\Delta)$, we have $u_1\cdots u_{m-1}=\bbar{u_0}$ by definition and thus $\phi.u\in D_m(\Delta)$. Now, for two positive integers $m,n$, we define
\[D_m^n(\Delta):=\left\{u\in D_m(\Delta)~|~ \phi^n.u=u \right\}.\]

In order to study $(p,q)$-periodic elements, we are going to define a categorical presentation using the sets $D_{kp}^{kq}(\Delta)$ with $k\in \{1,2,3\}$ (see Appendix \ref{sec:presentations_categories} for reminders on categorical presentations). 

 The subcategory $\Cc^{\phi^q}$ introduced at the end of Section \ref{sec:1.1} is useful for giving a more efficient description of the sets $D_{kp}^{kq}(\Delta)$. Those sets are a priori described by $kp$ parameters lying in $\Cc$, but since $p$ and $q$ are coprime, we show that they depend only on $k$ parameters lying in $\Cc^{\phi^q}$.

\begin{lem}\label{lem:descriptionpointfixesphiq}
For all $f:=(f_0,\ldots,f_{p-1})\in \Cc^p$, we have
\[\phi^q.f=f\Leftrightarrow \phi^q(f_0)=f_0\text{~and~}\forall i\in \intv{1,p-1},~f_i=\phi^{i\eta}(f_0).\]
In particular, $f$ depends only on $f_0$.
\end{lem}
\begin{proof}
Let $q=pm+r$ be the euclidean division of $q$ by $p$. For $f\in \Cc^p$, we have
\[\phi^q.f=(\phi(f_0),\ldots,\phi(f_{p-1})),\]
and thus,
\[\phi^q.f=(\phi^m(f_r),\ldots,\phi^{m}(f_{p-1}),\phi^{m+1}(f_0),\ldots,\phi^{m+1}(f_{r-1})).\]
$(\Leftarrow)$ Assume that $\phi^q(f_0)=f_0$, and that $\forall i\in \intv{1,p-1}$, we have $f_i=\sigma^{i\eta}(f_0)$. We then have
\begin{align*}
f&=\left(f_0,\sigma^\eta(f_0),\ldots,\sigma^{(p-1)\eta}(f_0)\right),\\
\phi^q.f&=\left(\sigma^{r\eta+m}(f_0),\ldots,\sigma^{(p-1)\eta+m}(f_0),\sigma^{m+1}(f_0),\ldots,\sigma^{m+(r-1)\eta+1}(f_0)\right).
\end{align*}
Since we have $pm+r=q$, we have $pm+r\equiv 0[q]$ and $\eta pm+r\eta\equiv 0[q]$. However, since $p\eta-q\mu=1$ implies $p\eta\equiv 1[q]$, we deduce that $m+r\eta\equiv 0[q]$. Since $f_0$ is invariant under $\phi^q$, we have
\[\begin{cases} \forall j\in \intv{0,p-r-1},~\phi^{(r+j)\eta+m}(f_0)=\phi^{j\eta}(f_0),\\
\forall j\in \intv{p-r,p-1},~ \phi^{m+(j-p+r)\eta+1}(f_0)=\phi^{j\eta-p\eta+1}(f_0)=\phi^{j\eta}(f_0).\end{cases}\]
and thus $f=\phi^q.f$.\\ $(\Rightarrow)$ Conversely, assume that $\phi^q.f=f$. We have by definition
\[\forall i\in \intv{0,p-1}, f_i=\begin{cases} \phi^m(f_{i+r})&\text{if }i<p- r,\\\phi^{m+1}(f_{i+r-p})&\text{if }i\geqslant p-r . \end{cases}\]
By an immediate induction, we obtain that 
\[f_0=\phi^{m}(f_r)=\ldots=\phi^{pm+r}(f_{pr})=\phi^q(f_0),\]
and $f_0$ is $\phi^q$-invariant. Now, we show that, for all $i\in \intv{0,p-1}$, we have $f_i=\phi^{i\eta}(f_0)$. The result is obvious for $i=0$. Now, assume that $f_i=\phi^{i\eta}(f_0)$ for some $i\in \intv{0,p-1}$. If $i>r$, then we have
\[f_{i-r}=\phi^m(f_i)=\phi^{-r\eta}(f_i)=\phi^{(i-r)\eta}(f_0)\]
since $m+r\eta\equiv 0[q]$, and since the value of $\phi^n(f_0)$ depends only on the value of $n$ modulo $q$. Likewise, if $i<r$, then we have
\[f_{i-r}=f_{i-r+p}=\phi^{m+1}(f_i)=\phi^{-r\eta+1+i\eta}(f_0)=\phi^{-r\eta+p\eta+i\eta}(f_0)=\phi^{(p-r+i)\eta}(f_0).\]
Since $r$ and $p$ are coprime, every $i$ in $\intv{0,p-1}$ can be reached from $0$ by successively subtracting $r$ (mod $p$). This terminates the proof.
\end{proof}

Notice that if $f_0\in \Cc$ is $\phi^q$-invariant, then $\phi^{\eta}(f)$ does not depend on the choice of $\eta$. Indeed, if $\eta',\mu'$ are such that $p\eta'-q\mu'=1$, then $\eta$ and $\eta'$ are equal modulo $q$. Using the above lemma, we give alternative descriptions of the sets $D_{kp}^{kq}(\Delta)$ for $k\in \{1,2,3\}$.

\begin{prop}\label{prop:descriptiondivset}
\begin{enumerate}[(a)]
\item The map $(u_0,\ldots,u_{p-1})\mapsto u_0$ induces a bijection between $D_p^q(\Delta)$ and the set 
\[\Dd_1:=\left\{u\in \Ss^{\phi^q}~|~ u\phi^{\eta}(u)\ldots\phi^{(p-1)\eta}(u)=\Delta\right\}.\]
\item The map $(a_0,b_0,\ldots,a_{p-1},b_{p-1})\mapsto (a_0,b_0)$ induces a bijection between $D_{2p}^{2q}(\Delta)$ and the set 
\[\Dd_2:=\{(a,b)\in (\Ss^{\phi^q})^2~|~ ab\in \Dd_1\}.\]
\item The map $(x_0,y_0,z_0,\ldots,x_{p-1},y_{p-1},z_{p-1})\mapsto (x_0,y_0,z_0)$ induces a bijection between $D_{3p}^{3q}(\Delta)$ and the set 
\[\Dd_3:=\{(x,y,z)\in (\Ss^{\phi^q})^3~|~ xyz\in \Dd_1\}.\]
\end{enumerate}
\end{prop}
\begin{proof}The first claim is a direct consequence of Lemma \ref{lem:descriptionpointfixesphiq} and the fact that $u_0\cdots u_{p-1}=\Delta$. For the second claim, let $(a_0,\ldots,b_{p-1})$ be a $2p$-tuple. This $2p$-tuple is $\phi^{2q}$-invariant if and only if the two $p$-tuples $(a_0,\ldots,a_{p-1})$ and $(b_0,\ldots,b_{p-1})$ are both $\phi^q$-invariant. By Lemma \ref{lem:descriptionpointfixesphiq}, this is equivalent to $\phi^q(a_0)=a_0,\phi^q(b_0)=b_0$ and 
\[(a_0,b_0,\ldots,a_{p-1},b_{p-1})=(a_0,b_0,\phi^{\eta}(a_0),\ldots,\phi^{(p-1)\eta}(b_0)),\]
which shows that the map $(a_0,b_0,\ldots,a_{p-1},b_{p-1})\mapsto (a_0,b_0)$ is injective and maps $D_{2p}^{2q}(\Delta)$ into $\Dd_2$. \newline Conversely, for $(a,b)$ in $\Dd_2$, the $2p$-tuple $(a_0,b_0,\phi^{\eta}(a_0),\ldots,\phi^{(p-1)\eta}(b_0))$ lies inside $D_{2p}^{2q}(\Delta)$ and gives the preimage of $(a,b)$. The same reasoning applies to the third claim.
\end{proof}

From now on, we will freely identify the sets $D_{kp}^{kq}(\Delta)$ and $\Dd_k$ for $k\in \{1,2,3\}$ using Proposition \ref{prop:descriptiondivset}. Under these identifications, the action of $\phi$ on $\Dd_1,\Dd_2,\Dd_3$ is given by
\begin{align*}
&\forall u \in \Dd_1,~\phi.u=\phi^\eta(u),\\
&\forall (a,b)\in \Dd_2,~\phi.(a,b)=(b,\phi^\eta(a)),\\
&\forall (x,y,z)\in \Dd_3,~\phi.(x,y,z)=(y,z,\phi^\eta(x)).
\end{align*}

We can now translate the construction of Bessis using the sets $\Dd_1,\Dd_2,\Dd_3$. We start by constructing an oriented graph which will serve as the generators for our categorical presentation.

\begin{definition}\label{def:graph_of_periodic_simples}
We denote by $\Ss_p^q$ the following oriented graph
\begin{itemize}
\item The objects are the elements of $\Dd_1$,
\item The arrows are the elements of $\Dd_2$,
\item The source (resp. target) of $(a,b)\in \Dd_2$ is $ab$ (resp. $b\phi^{\eta}(a)$).
\end{itemize}
We call $\Ss_p^q$ the \nit{graph of simples} (attached to the integers $p,q$).
\end{definition}

Note that, for $(a,b)\in \Ss_p^q$, the source of $(a,b)$ is in $\Dd_1$ by definition of $\Dd_2$, and its target is in $\Dd_1$ since $\phi.(a,b)=(b,\phi^{\eta}(a))$ is in $\Dd_2$. We then endow the graph $\Ss_p^q$ with a set of relations $R_p^q$, in bijection with $\Dd_3$. An element $(x,y,z)$ of $\Dd_3$ induces the relation $(x,yz)(y,z\phi^{\eta}(x))=(xy,z)$ between elements of $\Ss_p^q$.

\begin{definition}\label{def:periodic_elements_category} The presented category $\Cc_p^q:=\langle \Ss_p^q~|~R_p^q \rangle^+$ is called the \nit{category of $(p,q)$-periodic elements} of $\Cc$. The enveloping groupoid of $\Cc_p^q$ will be denoted by $\Gg_p^q$ and called the groupoid of $(p,q)$-periodic elements of $\Gg$.
\end{definition}

We now proceed to show that $\Cc_p^q$ is a homogeneous Garside category. The Garside part boils down to the fact that we translated the construction proposed by Bessis, which is already known to induce a Garside category. It then remains to prove homogeneity.

\begin{lem}\label{lem:lengthforperiodiccat}
The function $\ell$ defined on arrows of $\Ss_p^q$ by $\ell(a,b):=\ell(a)$ extends to a length functor on $\Cc_p^q$, making $\Cc_p^q$ into a homogeneous category.\end{lem}
\begin{proof}
First, the function $\ell$ extends to $\Cc_p^q$ as it is compatible with $R_p^q$. Indeed, for $(x,y,z)\in D_{3p}^{3q}(\Delta)$, we see that $\ell(xy,z)=\ell(xy)=\ell(x)+\ell(y)=\ell(x,yz)+\ell(y,z\phi^{\eta}(x))$. We then have to show that $\Cc_p^q$ is generated by elements of positive length. As $\Cc_p^q$ is generated (by definition) by $\Ss_p^q$, we have to show that the elements of length $0$ in $\Ss_p^q$ are trivial in $\Cc_p^q$. By definition of $\ell$, the only elements of $\Ss_p^q$ of length $0$ are elements of the form $(1_x,u)$ for $u\in D_p^q(\Delta)$ with source $x$ in $\Cc$. Since we are working with tuples of composable morphisms in $\Cc$, we will omit the subscript when writing identity morphisms in $\Cc$, writing $(1,u)$ instead of $(1_x,u)$. We claim that the element $(1,u)\in \Ss_p^q$ is the identity morphism of the object $u$ in $\Cc_p^q$. The source and target of $(1,u)$ are both $u$. Furthermore, for $(a,b)\in \Ss_p^q$ with source $u$, the element $(1,a,b)\in D_{3p}^{3q}(\Delta)$ expresses the relation $(1,u)(a,b)=(a,b)$. Similarly, for $(d,e)\in \Ss_p^q$ with target $u$, the element $(d,1,e)$ in $D_{3p}^{3q}(\Delta)$ expresses the relation $(d,e)(1,u)=(d,e)$. Thus $(1,u)$ represents the morphism $1_u$ in $\Cc_p^q$.
\end{proof}

\begin{theo}\label{theo:periodic_category_is_garside}\cite[Proposition XIV.1.5]{ddgkm} and \cite[Proposition VII.4.2]{ddgkm}\newline The category $\Cc_p^q$ is cancellative. Furthermore, defining
\[\forall u \in \Ob(\Cc_p^q),~~\Delta_p(u):=(u,1)\in \Cc_p^q(u,\phi_p(u))\]
induces a Garside map $\Delta_p$ on $\Cc_p^q$. The simple morphisms associated with $\Delta_p$ are exactly the images in $\Cc_p^q$ of the graph $\Ss_p^q$.
\end{theo}
\begin{proof}
Following the discussion before \cite[Definition 9.3]{besgar}, the $p$-divided category is defined via a categorical presentation involving the sets $D_{kp}(\Delta)$ for $k\in \{1,2,3\}$. The set $D_p(\Delta)$ forms the objects of an oriented graph whose arrows are the elements of $D_{2p}(\Delta)$. The set $D_{3p}(\Delta)$ then gives a set of relations on the resulting graph, providing the desired presentation. The result of \cite[Proposition XIV.1.5]{ddgkm} is that the category $\Cc_p$ defined by this presentation is a Garside category, whose simple morphisms are exactly the (images of) the elements of $D_{2p}(\Delta)$. The Garside map $\Delta_p(u)$ of an object $u=(u_0,\ldots,u_{p-1})\in D_p(\Delta)$ is $(u_0,1,\ldots,u_{p-1},1)$. Moreover, the action of the Garside automorphism $\phi_p$ of $(\Cc_p,\Delta_p)$ on simple morphisms is the action of $\phi^2$.

Now, by \cite[Proposition VII.4.2]{ddgkm}, the subcategory $(\Cc_p)^{\phi_p^q}$ of $(\phi_p)^q$-invariant elements of $\Cc_p$ is generated by the subgraph of $(\phi_p)^q$ invariant elements of $D_{2p}(\Delta)$, that is $D_{2p}^{2q}(\Delta)$. Moreover, a presentation is obtained by only considering relations induced by elements of $D_{3p}^{3q}(\Delta)$. Lastly, the category $(\Cc_p)^{\phi^q}$ is a Garside category whose Garside map is the restriction of the Garside map of $\Cc_p$ and whose simple morphisms are the (images of) elements of $D_{2p}^{2q}(\Delta)$. The translation of the presentation of $(\Cc_p)^{\phi^q}$ to our context using Proposition \ref{prop:descriptiondivset} is precisely the presentation we gave of $\Cc_p^q$, which gives the result.
\end{proof}

As mentioned in the proof, the action of the Garside automorphism on simple morphisms is the action of $\phi^2$. We then obtain the value of the Garside automorphism on simple morphisms of $\Cc_p^q$:

\begin{cor}
The Garside automorphism $\phi_p$ of the Garside category $(\Cc_p^q,\Delta_p)$ is given on simple morphisms by $\phi_p(a,b)=(\phi^{\eta}(a),\phi^{\eta}(b))$.
\end{cor}

We finish this section by describing properties of categories of $(p,q)$-periodic elements. We do so by relating the category $\Cc_p^q$ with the category $\Cc$ using the collapse functor introduced by Bessis.

\begin{definition}\cite[Definition B.23]{beskpi1} The map defined on $\Ss_p^q$ by $(a,b)\mapsto a$ induces a functor $\pi_p:\Cc_p^q\to \Cc$, which is called the \nit{collapse functor}.
\end{definition}

\begin{lem}\label{lem:preceqindivided}
Let $u$ be an object of $\Cc_p^q$, and let $x:=\pi_p(u)$ be the source of $u$ in $\Cc$. The collapse functor restricts to an isomorphism of posets 
\[(\Ss_p^q(u,-),\preceq)\simeq (\{s\in \Ss^{\phi^q}~|~ s\preceq u\},\preceq)\subset \Ss^{\phi^q}(x,-)\subset \Ss(x,-).\]
In particular, if $s:=(a,b)$ and $s':=(a',b')$ are two simples in $\Cc_p^q$ with source $u$, then $s\preceq s'$ in $\Cc_p^q$ if and only if $a\preceq a'$ in $\Cc$.
\end{lem}
\begin{proof}
Let $\pi$ denote the restriction of the collapse functor $\pi_p$ to the set $\Ss_p^q(u,-)$. For $s:=(a,b)\in \Ss_p^q(u,-)$, we have $ab=u$ by definition, so $a\preceq u$ and $\phi^q(a)=a$. Conversely, if $a\preceq u$ is such that $\phi^q(a)=a$, then there is a unique $b\in \Ss$ such that $ab=u$. We have $\phi^q(b)=b$ by Lemma \ref{lem:deux_sur_trois}, and the simple morphism $(a,b)$ is then the unique preimage by $\pi$ of $a$ in $\Ss_p^q(u,-)$.

Now, let $s:=(a,b)$ and $s':=(a',b')$ be two elements of $\Ss_p^q(u,-)$. If $s\preceq s'$, then there is some $(x,y,z)\in D_{3p}^{3q}(\Delta)$ such that $(x,yz)=(a,b)$ and $(xy,z)=(a',b')$. We have in particular $ay=xy=a'$ and $a\preceq a'$. Conversely, if $a\preceq a'$, there is some $y\in \Ss$ such that $ay=a'$. Again, we have $\phi^q(y)=y$ by Lemma \ref{lem:deux_sur_trois}, and the relation induced by $(a,y,b)$ then gives $s\preceq s'$.
\end{proof}

\begin{cor}\label{cor:greedy_normal_form_length_2_in_periodic}
Let $s:=(a,b)$ and $t:=(\alpha,\beta)$ be two composable simple morphisms in $\Cc_p^q$. The greedy normal form of the path $st$ in $\Cc_p^q$ is given by
\[st=(ad,d^{-1}b)(d^{-1}\alpha,\beta \phi^\eta(d)),\]
where $d=\alpha\wedge b$. In particular, the path $st$ is greedy if and only if $b$ and $\alpha$ are coprime.
\end{cor}
\begin{proof}
The composability of $s$ and $t$ is equivalent to the fact that the equality $\alpha\beta=b\phi^\eta(a)$ holds. By Lemma \ref{lem:preceqindivided}, we have $t\wedge \bbar{s}=(\alpha,\beta)\wedge (b,\phi^\eta(a))=(d,e)$, where $e$ is defined by $de=\alpha\beta$. We have 
\begin{align*}
st&=(a,b)(\alpha,\beta)=(a,b)(d,e)(d^{-1}\alpha,\beta)\\
&=(ad,d^{-1}b)(d^{-1}\alpha,\beta\phi^\eta(d)).
\end{align*}
This last term is the greedy normal form of $st$ because of Lemma \ref{lem:greedylength2}.
\end{proof}

During the proof of Theorem \ref{theo:periodic_category_is_garside}, we also proved that the category of $(p,q)$-periodic elements is isomorphic to the construction proposed by Bessis (that is, category of $(\phi_p)^q$-invariant elements of the $p$-divided category $\Cc_p$). Thus results on the latter also apply to the former. In particular we have the following theorem:

\begin{theo}\cite[Proposition XIV.1.8]{ddgkm}\label{theo:categoriediviseeetelementsreguliers}
\begin{enumerate}[(a)]
\item Let $u$ be an object of $\Gg_p^q$. The element $\pi_p(\Delta_p^q(u))$ is a $(p,q)$-periodic element of $\Gg$.
\item Every $(p,q)$-regular element of $\Gg$ is conjugate to an element of the form $\pi_p(\Delta_p^q(u))$.
\item The collapse functor induces a functor $\Gg_p^q\to \Gg$, which in turn induces a group isomorphism between $\Gg_p^q(u,u)$ and the centralizer of $\pi_p(\Delta_p^q(u))$ in $\Gg(\pi_p(u),\pi_p(u))$.
\end{enumerate}
\end{theo}

\subsubsection{The case $p=1$}\label{sec:1.3.1} If $p=1$, then the graph $\Ss_p^q$ is simply the subgraph $\Ss^{\phi^q}$ of $\Ss$. Indeed, in this case, we have
\begin{align*}
\Dd_1&=\{u\in \Ss^{\phi^q}~|~u=\Delta\}=\{\Delta(u)\in \Ss~|~\phi^q(\Delta(u))=\Delta(\phi^q(u))=\Delta(u)\}\\
&\simeq \{u\in \Ob(\Cc)~|~ \phi^q(u)=u\}.\\
\Dd_2&=\{(a,b)\in (\Ss^{\phi^q})^2~|~ab=\Delta\}\\
&=\{(a,\bbar{a})\in (\Ss^{\phi^q})^2\}\simeq \Ss^{\phi^q}.
\end{align*}
Let $(x,y,z)$ be in $\Dd_3$. We have $z=\bbar{xy}$, thus $(x,y,z)$ expresses the relation $x.y=(xy)$ between elements of $\Ss^{\phi^q}$. 

By \cite[Proposition VI.1.11]{ddgkm}, the category $\Cc_1^q$ is simply the subcategory $\Cc^{\phi^q}$. The Garside map of Theorem \ref{theo:periodic_category_is_garside} is then the restriction to $\Cc^{\phi^q}$ of the Garside map of $\Cc$. The same goes for the automorphism $\phi_1$.

The collapse functor $\pi_1$ sends $a\in \Ss^{\phi^q}$ to $a\in \Ss$: it is the inclusion $\Cc^{\phi^q}\hookrightarrow \Cc$. Lemma \ref{lem:preceqindivided} expresses that, for an object $u$ of $\Cc^{\phi^q}$, the inclusion $\Ss^{\phi^q}(u,-)\subset \Ss(u,-)$ is a morphism of lattices, which is also a consequence of Lemma \ref{lem:deux_sur_trois}. Corollary \ref{cor:greedy_normal_form_length_2_in_periodic} ensures that the greedy normal form of a path in $\Cc^{\phi^q}$ is the same as its greedy normal form in $\Cc$.

Lastly, Theorem \ref{theo:categoriediviseeetelementsreguliers} expresses that every $(1,q)$-periodic element of $\Cc$ is conjugate to some $\Delta^q(u)$, and that $\Gg(\Cc)^{\phi^q}(u,u)$ is the centralizer of $\Delta^q(u)$ in $\Gg(\Cc)(u,u)$.

\section{Complex reflection groups, braid groups and braid groupoids}\label{sec:crg_bg_braid_cats}

\subsection{Reminders on complex reflection groups and their braid groups}
\subsubsection{Definitions, basic invariants}
In this section we fix a positive integer $n$ and a complex vector space $V$ of dimension $n$. We follow the exposition of \cite{lehrertaylor}. A nontrivial element $s\in \GL(V)$ is called a \nit{reflection} if it pointwise fixes some hyperplane of $V$ and has finite order. The hyperplane $\Ker(s-\Id)$ is then called the \nit{reflecting hyperplane} of $s$. A subgroup $W< \GL(V)$ is called a \nit{complex reflection group} if it is finite and generated by reflections of $V$. The codimension in $V$ of the subspace $V^W$ of fixed points under the action of $W$ is the \nit{rank} of $W$.

A complex reflection group $W< \GL(V)$ is \nit{irreducible} if there are no nontrivial $W$-stable subspaces in $V$. Every complex reflection group decomposes as a direct product of irreducible complex reflection groups. This allows us to restrict our attention to irreducible groups, which were classified by Shephard and Todd in \cite{shetod}. The classification separates on the one hand an infinite series of groups related to groups of monomial matrices, and on the other hand a list of $34$ exceptional cases, labeled $G_{4}$ to $G_{37}$.

An irreducible complex reflection group of rank $n$ is called \nit{well-generated} if it can be generated by a set of $n$ reflections. If this is not the case then the group is said to be \nit{badly-generated}. A badly-generated irreducible complex reflection group of rank $n$ can always be generated by $n+1$ reflections.

We fix an irreducible complex reflection group $W<\GL(V)$ for the remainder of this section. Since $W$ is irreducible, the space $V^W$ is trivial and thus $n$ is the rank of $W$. 

The action of $W$ on $V$ extends to an action of $W$ on the algebra $S$ of polynomial functions on $V$. By the Chevalley-Shephard-Todd Theorem \cite[Theorem 3.20]{lehrertaylor}, the algebra $S^W$ of $W$-invariant elements of $S$ is a polynomial algebra in $n$ variables. An $n$-tuple $f=(f_1,\ldots,f_n)$ of homogeneous algebraically independent generators of $S^W$ is called a \nit{system of basic invariants}. The degrees $d_i$ of the $f_i$ do not depend on the choice of $f$. They are the \nit{degrees} of the reflection group $W$.

Likewise, one can define the \nit{codegrees} $d_1^*,\ldots,d_n^*$ of $W$ by considering the module of invariant derivations on the algebra $S$ \cite[Definition 1.2]{beskpi1}.
\subsubsection{Braid groups, braided reflections}
We consider the subset $X\subset V$ consisting of points not belonging to any of the reflecting hyperplanes associated with reflections of $W$. By Steinberg's Theorem \cite[Theorem 9.44]{lehrertaylor}, the action of $W$ on $X$ is free and induces a covering map from $X$ to $X/W$. We fix a basepoint $x_0\in X$ and we set
\[P(W):=\pi_1(X,x_0) \text{~and~} B(W):=\pi_1(X/W,W.x_0)\]
the \nit{pure braid group} of $W$ and the \nit{braid group} of $W$, respectively. The covering map $X\twoheadrightarrow X/W$ induces a short exact sequence
\[\xymatrix{1 \ar[r] & P(W)\ar[r] & B(W)\ar[r] & W\ar[r] & 1}.\]
Note that, as $X$ is path connected, a change of basepoint $x_0\to x_1$ yields a (non canonical) isomorphism of short exact sequences
\[\xymatrix{1\ar[r] & P(W)\ar[r] \ar[d]^-{\simeq} & B(W)\ar[d]^-{\simeq}\ar[r] & W\ar@{=}[d] \ar[r] & 1\\
1\ar[r] & \pi_1(X,x_1)\ar[r] & \pi_1(X/W,W.x_1)\ar[r] & W\ar[r] & 1.}\]

Let $f:=(f_1,\ldots,f_n)$ be a system of basic invariants. The isomorphism between $S^W$ and $\C[f_1,\ldots,f_n]$ induces in turn an algebraic isomorphism $V/W\simeq \C^n$, which sends the orbit $W.v$ of $v\in V$ to $(f_1(v),\cdots,f_n(v))\in \C^n$. Let $\Hh$ denote the image in $V/W$ of the union of the reflecting hyperplanes of $W$. It is an algebraic hypersurface of $V/W\simeq \C^n$. The \nit{discriminant} of $W$ relative to $f$ is the polynomial $\mathrm{Disc}(W,f)$ defining the image of $\Hh$ under the isomorphism $V/W\simeq \C^n$ induced by $f$. 

Since the braid group $B(W)$ is defined as the fundamental group of the complement of $\Hh$ in $V/W$, it is generated by particular elements called \nit{braided reflections} (around the irreducible divisors of $\Hh$).

We quickly recall the definition of braided reflections from \cite[Appendix 1]{bmr} , as it will be useful in Section \ref{sec:braid_reflections_and_atomic_loops}. Let $H$ be a reflecting hyperplane of $W$, and let $L\subset \Hh$ be the image of $H$ in $V/W$. Let also $\ttilde{L}$ be the image in $V/W$ of points of $V$ which belong to $H$ and to no other reflecting hyperplane of $W$.

\begin{definition}
A \nit{path from $x_0$ to $L$ in $X/W$} is a path $\gamma:[0,1]\to V/W$ such that $\gamma(0)=x_0$, $\gamma(1)\in \ttilde{L}$ and $\gamma(t)\in X/W$ for $t<1$. Two paths $\gamma,\gamma'$ from $x_0$ to $L$ in $X/W$ are $L$-\nit{homotopic} if there exists a homotopy $T:[0,1]\times[0,1]\to V/W$ from $\gamma$ to $\gamma'$ such that
\begin{itemize}
\item For all $t\in [0,1[$ and $u\in[0,1]$ $T(t,u)\in X/W$.
\item For all $u\in [0,1]$, $T(0,u)=x_0$ and $T(1,u)\in \ttilde{L}$.
\end{itemize}
The $L$-homotopy class of $\gamma$ is denoted by $[\gamma]$.
\end{definition}

Let $\gamma$ be a path from $x_0$ to $L$ in $X/W$. Let also $U$ be a connected open neighborhood of $\gamma(1)$ in $X/W\cup \ttilde{L}$ such that $U\cap X/W$ has a fundamental group free of rank $1$. Let $u\in [0,1[$ be such that $\gamma(t)\in U$ for $t\geqslant u$. The orientation of $U\cap X$ allows us to choose a ``positive'' generator $\lambda$ of $\pi_1(U\cap X/W,\gamma(u))$. 

We set $\gamma_u(t):=\gamma(ut)$ for $t\in [0,1]$ and $\rho_{\gamma,\lambda}:=\gamma_u*\lambda*\gamma_u^{-1}$. The homotopy class of $\rho_{\gamma,\lambda}$ depends only on $[\gamma]$ and is denoted by $\rho_{[\gamma]}$. It is by definition a braided reflection of $B(W)$ (around $L$). The image of $\rho_{[\gamma]}$ inside $W$ is a reflection with hyperplane $H$.

\begin{prop}\cite[Appendix 1]{bmr} All braided reflections around $L$ form a conjugacy class of $B(W)$. In particular the set of all braided reflections is stable under conjugacy in $B(W)$.
\end{prop}

\begin{rem}\label{rem:basepoints_and_braid_reflections}
Let $x_0\to x_1$ be a change of basepoint. The induced isomorphism between $B(W)$ and $\pi_1(X/W,W.x_1)$ maps braided reflections to braided reflections.
\end{rem}

Lastly, we give the definition of a particular element in $B(W)$, which plays an important role when considering centers and centralizers.
\begin{definition}\cite[Notation 2.3]{bmr} The \nit{full-twist} is the homotopy class in $B(W)$ of the loop
\[\begin{array}{cccl}\gamma_\pi:&[0,1]&\longrightarrow &X/W\\
 &t&\longmapsto & \exp(2i\pi t)x_0.\end{array}\]
It lies inside $P(W)\cap Z(B(W))$.
\end{definition}

\subsubsection{Regular elements, regular braids}
Let $d$ be a positive integer. We denote by $\mu_d$ the group of $d$-th roots of unity in $\C$. We also consider $\mu_d^*$ the subgroup of $\mu_d$ consisting of primitive $d$-th roots of unity, and $\zeta_d:=\exp(\frac{2i\pi}{d})\in \mu_d^*$.

\begin{definition}\cite[Definition 11.21]{lehrertaylor} Let $\zeta$ be a root of unity in $\C$. An element $g\in W$ is called $\zeta$\nit{-regular} if it admits a $\zeta$-eigenvector lying in $X$. An integer $d$ is \nit{regular} for $W$ if there is a $\zeta_d$-regular element $g$ in $W$.
\end{definition}

If $g\in W$ is a $\zeta$-regular element, then for every integer $k$, $g^k$ is $\zeta^k$-regular. Thus, if $g$ is $\zeta$-regular for some $\zeta\in \mu_d^*$, then $d$ is a regular integer for $W$. It is known that $\zeta$-regular elements (should they exist) form a conjugacy class in $W$.

\begin{theo}\label{theo:degrees_of_regular_centralizer}\cite[Theorem 11.24]{lehrertaylor} Let $g\in W$ be a $\zeta$-regular element for some $\zeta\in \mu_d^*$. The centralizer $C_W(g)$ of $g$ in $W$ acts on the eigenspace $\Ker(g-\zeta\Id)$ as a complex reflection group. Its degrees (resp. codegrees) coincide with those degrees (resp. codegrees) of $W$ which are divisible by $d$.
\end{theo}

Let $g\in W$ be a $\zeta_d$-regular element, and let $W_g:=C_W(g)$. It is possible to study the braid group $B(W_g)$ by embedding it inside $B(W)$. Let $V_g$ denote the eigenspace $\Ker(g-\zeta_d \Id)$ on which $W_g$ acts as a reflection group. Let also $X_g$ denote the space of regular vectors inside $V_g$. By \cite[Theorem 11.33]{lehrertaylor}, we have $X_g=X\cap V_g$.

The scalar action of $\C^*$ on $V$ induces in particular an action of $\mu_d$ on $V/W$. By \cite[Theorem 1.9]{beskpi1}, the embedding $V_g\to V$ induces two homeomorphisms $V_g/W_g\simeq (V/W)^{\mu_d}$ and $X_g/W_g\simeq (X/W)^{\mu_d}$. In particular $B(W_g)$ identifies with the fundamental group of $(X/W)^{\mu_d}$. It was shown in \cite[Theorem 12.4]{beskpi1} and \cite[Theorem 1.2]{regularbraids} that the embedding $(X/W)^{\mu_d}\to X/W$ identifies the fundamental group $(X/W)^{\mu_d}$ with the centralizer in $B(W)$ of some $d$-th root of the full-twist. Such $d$-th roots of the full-twist will be called $d$-\nit{regular braids}. Note that \cite[Theorem 12.4]{beskpi1} and \cite[Theorem 1.2]{regularbraids} also prove that $d$-regular braids in $B(W)$ (should they exist) form a conjugacy class in $B(W)$.

\begin{exemple}\label{ex:=g31_dans_g37}
A guiding example is given by the groups $G_{31}$ and $G_{37}$. The degrees (resp. codegrees) of $G_{37}$ are
\[\begin{array}{cccccccc}2&8&12&14&18&20&24&30\\
0&6&10&12&16&18&22&28\end{array}\]
The integer $4$ is regular for $G_{37}$. The centralizer of a $i=\zeta_4$-regular element in $G_{37}$ is isomorphic to $G_{31}$. Its degrees (resp. codegrees) are $8,12,20,24$ (resp. $0,12,16,28$). Note that $4$ is the gcd of the degrees of $G_{37}$ which are divisible by $4$.
\end{exemple}

\begin{lem}\label{lem:reflectionsetconjcommute}
Let $W< \GL_n(\R)\subset \GL_n(\C)$ be a real reflection group which contains $-\Id$ and for which $4$ is regular. If $g\in W$ is a $i$-regular element and $r$ is a reflection of $W$, then $r$ and $r^g$ commute.
\end{lem}
\begin{proof}
First, the element $-\Id$ is a $-1$-regular element in $W$. Since $-\Id$ is central in $W$, and since all $-1$-regular elements in $W$ are conjugate to $-\Id$, we get that $-\Id$ is in fact the only $-1$-regular element in $W$. Now, since $g$ is $i$-regular, $g^2$ is $i^2=-1$-regular, hence equal to $-\Id$.

We assume that $r\neq r^g$ (otherwise our claim is immediate). We fix $\scal{-}{-}$ a $W$-invariant scalar product on $\R^n$. Let $\alpha$ be a root for $r$ (that is, a generator of $\Ker(r-\Id)^\perp$). We have
\[g^2rg^{-2}=r\Rightarrow grg^{-1}=g^{-1}rg.\]
Thus $g^{-1}(\alpha)$ and $g(\alpha)$ are two roots of the same reflection $r^g$. Let then $\lambda\in \R$ be such that $g^{-1}(\alpha)=\lambda g(\alpha)$. We have $\alpha=\lambda g^2(\alpha)=-\lambda \alpha$ and $\lambda=-1$, thus
\[\scal{g(\alpha)}{\alpha}=\scal{\alpha}{g^{-1}(\alpha)}=\scal{\alpha}{-g(\alpha)}=-\scal{g(\alpha)}{\alpha}.\]
The two roots $\alpha$ and $g(\alpha)$ are then orthogonal for $\scal{-}{-}$ and $r$ and $r^g$ commute.
\end{proof}

This lemma applies to the case of $G_{31}$ seen inside $G_{37}$ as the centralizer of a $i$-regular element. Indeed, $G_{37}$ is the complexified version of the exceptional Coxeter group $E_8$, which contains $-\Id$.

\subsection{Dual braid monoid for well-generated groups} \label{dbm&dc} In this section we fix a finite dimensional complex vector space $V$, along with a well-generated complex reflection group $W<\GL(V)$. In \cite{besdbr} and \cite{beskpi1}, Bessis defines a particular Garside monoid, the \nit{dual braid monoid}, which admits the braid group $B(W)$ as its group of fractions. The construction of this monoid is the first step required in order to define a Garside category suitable for studying centralizers of regular braids in $B(W)$.

Here we consider the combinatorial definition of the dual braid monoid as an interval monoid. We will be discussing the topological definition of the dual braid monoid in the next section.  We follow \cite[Section 8]{beskpi1}.

Let $R$ be the set of reflections of $W$. We have $R=R^{-1}$, and $R$ generates $W$ by definition. For $w\in W$, we denote by $\ell_R(w)$ the minimal length of a decomposition of $w$ as a product of reflections. Because the set $R$ is globally invariant under conjugacy, the function $\ell_R$ is constant on conjugacy classes in $W$. The function $\ell_R$ induces relations $\preceq$ and $\succeq$ on $W$ defined by
\[\forall w,v\in W,~v\preceq w~\Leftrightarrow~ \ell_R(v)+\ell_R(v^{-1}w)=\ell_R(w),\]
\[\forall w,v\in W,~w\succeq v~\Leftrightarrow~ \ell_R(wv^{-1})+\ell_R(v)=\ell_R(w).\]
Let $v,w$ be in $W$. Since $\ell_R$ takes constant values on conjugacy classes in $W$, we get that $v\preceq w\Leftrightarrow w\succeq v$. Thus we will always work with $\preceq$ and only consider $\succeq$ for readability purposes.

As $W$ is well-generated, the highest degree $h$ of $W$ is regular for $W$ (see the proof of \cite[Theorem 2.4]{beskpi1}). A \nit{Coxeter element} of $W$ is then an element $c\in W$ which is $\zeta_h$-regular. Let $c$ be a Coxeter element in $W$. We define 
\[R_c:=\{r\in R~|~r\preceq c\} \text{~and~} I_c:=[1,c]_{\preceq}=\{w\in W~|~ w\preceq c\}.\]
We also consider copies $\Rr\subset \Ii$ of $R_c\subset I_c$, whose elements will be denoted in bold font. 
The \nit{dual braid monoid} associated with $R$ and $c$ is then defined by the following monoid presentation:
\[M(c):=\left\langle \Ii~|~ \mathbf{st}=\mathbf{u}\Leftrightarrow (st=u\text{ and }\ell_R(s)+\ell_R(t)=\ell_R(u))\right\rangle^+.\]
The enveloping group of $M(c)$ is denoted by $G(c)$. The above monoid presentation, when seen as a group presentation, provides a presentation of $G(c)$.

The poset $(I_c,\preceq)$ is a lattice by \cite[Lemma 8.6]{beskpi1}. We will denote by $s\wedge t$ (resp. $s\vee t$) the gcd (resp. lcm) of $s$ and $t$ in $I_c$. We do not have to distinguish between left- and right-gcds and lcms since all simple elements are balanced (their left- and right-divisors coincide).

Let $c$ and $c'$ be two Coxeter elements of $W$. Since they are regular elements of $W$ for the same eigenvalue $\zeta_h$, there is some $w\in W$ such that $wcw^{-1}=c'$. Since $\ell_R$ is invariant under conjugacy, one readily sees that conjugation by $w$ induces an isomorphism between the two intervals $I_c$ and $I_{c'}$. In particular we have $M(c)\simeq M(c')$.

\begin{theo}\cite[Theorem 8.2]{beskpi1} The monoid $M(c)$ is a homogeneous Garside monoid, with Garside element $\mathbf{c}$ and set of simples $\Ii$. The atoms of $M(c)$ are the elements of $\Rr$. 
\end{theo}

The \nit{Hurwitz relations} with respect to $W$ and $c$ are the formal relations of the form
\[\mathbf{rr'}=\mathbf{r'r''},\]
where $r,r',r''\in R_c$ are such that $r\neq r'$, $rr'\in I_c$ and $rr'=r'r''$ holds in $W$. We know by \cite[Lemma 8.8]{beskpi1} that the monoid $M(c)$ admits a particular presentation with $\Rr$ as a set of generators, endowed with the Hurwitz relations.

\begin{rem}\label{rem:egalitesdesimpledelongueur2} We can be a little less specific and say that $M(c)$ (and thus $G(c)$) is presented by relations of the form
\[\mathbf{ss'}=\mathbf{rr'},\]
with $s,s',r,r'\in R$, and $ss'\in I_c$. That is $M(c)$ is presented by its atoms and the equality between decompositions of simple elements of length $2$. It is this precise rephrasing which we will show holds in Springer categories (see Theorem \ref{prop:hurwitzpresentationforc31}).
\end{rem}

\begin{lem}\label{2.4}
For $\bbf{a}\in \Ii$, we have $\bbf{a}^2\in \Ii$ if and only if $\bbf{a}=1$.
\end{lem}
\begin{proof}
The presentation of $M(c)$ is obtained by a germ structure on the interval $I_c$. This germ structure is a so-called bounded Garside germ by \cite[Proposition IX.2.4]{ddgkm}. This shows that $(M(c),\bbf{c})$ is a Garside monoid with set of simples $\Ii$, but also that for any $\bbf{s},\bbf{t}\in \Ii$, if $m:=\bbf{s}\bbf{u}\in \Ii$, then $m=\bbf{s}\bbf{u}$ is a defining relation of $M(c)$.

Assume now that $\bbf{a}^2\in \Ii$, and write $\bbf{b}:=\bbf{a}^2$. By the above argument, $\bbf{a}^2=\bbf{b}$ is a defining relation of $M(c)$. We then have $a^2=b$ and $2\ell_R(a)=\ell_R(b)=\ell_R(a^2)$. If $a\neq 1$, then let $s\in R_c$ be a divisor of $a$ with $sa'=a$ and $a''s=a$. We have $a^2=a''ssa=a''s^2a$, so $\ell_R(a^2)\leqslant \ell_R(a'')+\ell_R(a')+1=2\ell_R(a)-1$, which contradicts $2\ell_R(a)=\ell_R(a^2)$.
\end{proof}

\begin{prop}\label{prop:produit_de_simple_simple_dans_dual}
Let $\bbf{a,a'}$ be in $\Ii$. If the product $\bbf{aa'}$ is also in $\Ii$, then $\bbf{aa'}=\bbf{a}\vee\bbf{a'}$ and $\bbf{a}\wedge \bbf{a'}$ is trivial.
\end{prop}
\begin{proof}
First, note that in $I_c$, we have $aa'=a'a^{a'}={}^aa'a$. Now, because $\ell_R$ is invariant under conjugacy, we have
\[\ell_R(aa')=\ell_R(a)+\ell_R(a')=\ell_R(a')+\ell_R\left(a^{a'}\right)=\ell_R\left({}^aa'\right)+\ell_R(a).\]
This means that $aa'\succeq a^{a'}$ and ${}^a a'\preceq aa'$: we have $a^{a'},{}^aa'\in I_c$, and $\bbf{aa'}$ is a common left and right multiple of $\bbf{a}$ and $\bbf{a'}$. Now let $\bbf{ab}=\bbf{a'b'}$ be the right-lcm of $\bbf{a}$ and $\bbf{a'}$. By definition, there is some $\bbf{x}\in \Ii$ such that $\bbf{bx}=\bbf{a'}$ and $\bbf{b'x}=\bbf{a^{a'}}$. In particular, we have $\bbf{a'}\succeq \bbf{x}$, $\bbf{a^{a'}}\succeq \bbf{x}$ and $\bbf{x}\preceq \bbf{a^{a'}}$ since every element is balanced. So $\bbf{aa'}$ admits $\bbf{x}^2$ as a factor and $\bbf{x}^2$ is a simple. By Lemma \ref{2.4}, $\bbf{x}=1$, $\bbf{a'}=\bbf{b}$, $\bbf{b'}=\bbf{a^{a'}}$ and $\bbf{aa'}=\bbf{a}\vee \bbf{a'}$. 

If $\bbf{x}$ is a common divisor of $\bbf{a}$ and $\bbf{a'}$, then we have $\bbf{a}\succeq \bbf{x}$ and $\bbf{x}\preceq \bbf{a'}$. We get that $\bbf{x}^2$ is a simple element in $M(c)$, it must then be trivial.
\end{proof}

\begin{rem}\label{rem:divided_set_interval_monoid}
By definition of an interval monoid, there is a natural bijection between the simples of $M(c)$ and the interval $I_c=([1,c],\preceq_R)$ in $W$. This bijection induces in turn a bijection between the sets $D_m(\bbf{c})$ and the sets $D_m(c)$, defined by
\[D_m(c):=\left\{(u_0,\ldots,u_{p-1})\in W^p~\left|~ \begin{cases} u_0\cdots u_{p-1}=c \\ \ell_R(u_0)+\cdots+\ell_R(u_{p-1})=\ell_R(c) \end{cases} \right. \right\}.\]
The automorphism $\phi$ of $M(c)$ acts on $D_m(c)$ by $\phi.(u_0,\ldots,u_{p-1}):=(u_1,\ldots,u_{p-1},u_0^c)$. The bijection between $D_m(\bbf{c})$ and $D_m(c)$ restricts to a bijection between $D_m^n(\bbf{c})$ and 
\[D_m^n(c):=\{(u_0,\ldots,u_{p-1})\in D_m(c)~|~ \phi^q.(u_0,\ldots,u_{p-1})=(u_0,\ldots,u_{p-1})\}.\]
The results of Proposition \ref{prop:descriptiondivset} also allow for an alternative description of the sets $D_n^m(c)$.
\end{rem}

\subsection{Lyashko-Looijenga map, cyclic labels}\label{sec:2.3} So far we have only given a combinatorial definition of dual braid monoids as interval monoids. In order to relate the dual braid monoid attached to a well-generated complex reflection group with its braid group, we also need a topological definition of dual braid monoids. We explain this topological definition in this section. We also explain the topological construction of the Springer category associated with a regular number for a well-generated complex reflection group. Most of this section comes from \cite[Section 5,6]{beskpi1}.

In this section, we fix a finite dimensional complex vector space $V$ with $n:=\dim(V)$, along with an irreducible well-generated complex reflection group $W<\GL(V)$. We denote by $h$ the highest degree of $W$. We otherwise keep the notation from the previous sections. We also fix a system $f$ of basic invariants for $W$ such that the attached discriminant has the following form:
\[\mathrm{Disc}(W,f)=X_n^n+\alpha_2X_n^{n-2}+\cdots+\alpha_n,\]
where $\alpha_i\in \C[X_1,\cdots,X_{n-1}]$ for $i\in \intv{2,n}$ (such a system always exists for well-generated groups by \cite[Theorem 2.4 and Section 2.5]{beskpi1}).

Recall that the system of basic invariants $f$ provides an isomorphism $V/W\simeq \C^n$, which sends the orbit $W.v$ of $v\in V$ to $(f_1(v),\cdots,f_n(v))\in \C^n$. We denote $v_i:=f_i(v)$.

\begin{definition}\cite[Definition 7.23]{beskpi1}
Let $x=W.v\in V/W$. The multiset $\LLL(x)$ is defined as the solutions in $T$ of the univariate polynomial
\[\mathrm{Disc}(W,f)(v_1,\cdots,v_{n-1},T+v_n)\in \C[T].\]
It is by definition an element of $E_n:=\C^n/\S_n$. The map $\LLL:V/W\to E_n$ is called the (extended) \nit{Lyashko-Looijenga morphism}.
\end{definition}

\begin{rem}
The original approach in \cite{beskpi1} was to study a slightly different application, denoted by $\LL$ \cite[Definition 5.1]{beskpi1}. Bessis later suggested that the application $\bbar{\LL}$ was in fact a far better choice.
\end{rem}

For $x\in V/W$, we have $x\in X/W$ if and only if $0\notin \LLL(x)$. By definition this is equivalent to the statement that $\bbar{\LL}(x)$ belongs to $E_n^\circ:=(\C^*)^n/\S_n$. We endow $V/W$ with the quotient of the scalar action of $\C$ on $V$. We also endow the space of configurations $E_n$ with the scalar action of $\C$.

\begin{lem}\label{lem:lll_homog}\cite[Lemma 11.1]{beskpi1} Let $x\in V/W$ and let $\lambda\in \C^*$. We have $\bbar{\LL}(\lambda x)=\lambda^h \bbar{\LL}(x)$.
\end{lem}

The \nit{fat basepoint} for $X/W$ is then defined as the following set
\[\Uu:=\{x \in V/W~|~ \bbar{\LL}(x)\cap i\R_{\geqslant 0}=\varnothing\}.\]
That is, $x\in \Uu$ if no point of $\bbar{\LL}(x)$ is directly above $0$. This subset of $X/W$ is open and contractible (\cite[Lemma 6.3]{beskpi1}). It can be used as a fat basepoint for defining the braid group of $W$ (\cite[Definition 6.4]{beskpi1}). From now on we consider $B(W):=\pi_1(X/W,\Uu)$.

Let $x \in X/W$, the points of $\bbar{\LL}(x)$ are ordered clockwise starting from ``right after noon''. In other words, we consider the argument of a point as an element of $[\pi/2,5\pi/2[$ and we order points by decreasing argument. Points with the same argument are ordered by increasing modulus. Following \cite[Definition 11.8]{beskpi1}, we define the \nit{cyclic support} of some $x\in V/W$ as the sequence $(x_1,\ldots,x_k)$ of points of $\LLL(x)$, ordered as explained above.

\begin{rem}\label{desingcyc}
The set of points $x \in \Uu$ whose elements of the cyclic support all have distinct arguments is dense in  $X/W$. The ordering we chose allows us to define univalent desingularization, which preserves the ordering of the points of the cyclic support.
\end{rem}

\begin{definition}\cite[Definition 11.24]{beskpi1} A \nit{circular semitunnel} is a couple $T=(x,L)\in \Uu\times [0,\frac{2\pi}{h}]$. 
The path $\gamma_T$ associated with a semitunnel $T$ is the path
\[\begin{array}{cccc} \gamma_T:&[0,1]&\longrightarrow& X/W \\ & t&\longmapsto &e^{itL}x.\end{array}\]
We say that $T$ is a \nit{circular tunnel} if it satisfies the additional condition $\gamma_T(1)\in \Uu$.
\end{definition}

Let $T:=(x,L)$ be a circular (semi)tunnel. Since $\LLL$ is homogeneous of degree $h$, the path $\LLL\circ \gamma_T$ corresponds to a continuous rotation of angle $hL$. Since the path $\gamma_T$ starts and ends in $\Uu$, it induces a well defined element in $B(W)=\pi_1(X/W,\Uu)$. From now on we amalgamate a circular tunnel $T$ with the path $\gamma_T$ it induces in $X/W$.

\begin{theo}\cite[Proposition 8.5 and Lemma 11.10]{beskpi1}\label{monomono}
\begin{enumerate}[(i)]
\item The homotopy class $\delta$ of the circular tunnel $T=(x,\frac{2\pi}{h})$ does not depend on $x\in \Uu$.
\item The image $c$ of $\delta$ in $W$ is a Coxeter element of $W$.
\item Consider $S$ the set of homotopy classes of circular tunnels in $B(W)$. It is endowed with the relation $\preceq$ defined by $s\preceq s'$ if $s^{-1}s'$ is homotopic to a circular tunnel. The projection map $B(W)\twoheadrightarrow W$ induces an isomorphism of posets $(S,\preceq)\simeq (I_c,\preceq)$. This isomorphism induces in turn an isomorphism $B(W)\simeq G(c)$ sending $\delta$ to $\bbf{c}$.
\end{enumerate}
\end{theo}

From now on we fix the Coxeter element $c$ so that we have the isomorphism $G(c)\simeq B(W)$ induced by the above theorem. By \cite[Lemma 6.13]{beskpi1}, the element $\delta^h$ in $B(W)$ represents the full-twist.

\begin{rem} In \cite{beskpi1}, the notion of circular tunnel is preceded by the notion of tunnel (cf \cite[Definition 6.6]{beskpi1}). However, \cite[Lemma 11.10]{beskpi1} and \cite[Corollary 6.18]{beskpi1} show that the two notions are actually synonymous, in the sense that an element of $B(W)$ is represented by a circular tunnel if and only if it is represented by a tunnel.
\end{rem}

Circular tunnels allow us to define the notion of cyclic label, which will prove to be crucial in the study of $(X/W)^{\mu_d}$. Let $x\in \Uu$, and let $(x_1,\ldots,x_k)$ be the cyclic support of $x$. We assume that all the $x_i$ have distinct arguments. There is a minimal $\theta\in \R_{>0}$ such that $e^{i\theta}x\notin \Uu$. The \nit{head} of $x$ is defined as the element $c_1$ of $B(W)$ represented by the circular tunnel $T=(x,\theta+\epsi)$, which does not depend on $\epsi>0$ small enough. We can then consider the head of $\gamma_T(1)=e^{i(\theta+\epsi)t}x$ and so on, until all the points of the cyclic support of $x$ have been labeled.

\begin{definition}\cite[Definition 11.9]{beskpi1} Let $x\in \Uu$ be such that all the points in $\bbar{\LL}(x)$ have distinct arguments. The sequence $(s_1,\ldots,s_k)$ defined above is the \nit{cyclic label} of $x$, denoted by $\clbl(x)$.\newline If different points in the cyclic support of $x$ have the same argument, we define $\clbl(x)$ as the cyclic label of some desingularization of $x$, as in Remark \ref{desingcyc}.
\end{definition}

Since the elements of $\clbl(x)$ are defined as circular tunnels, they are simple elements in $M(c)$. Thus we can see them as elements of $I_c\subset W$. Theorem \ref{monomono} proves that, for any $x\in X/W$, the product of all the terms of $\clbl(x)$ is equal to $\delta$ in $B(W)$, as it is represented by the circular tunnel $(x,\frac{2\pi}{h})$. Seeing circular tunnels as elements of $W$, if the cyclic support of $x$ contains $k$ points, then $\clbl(x)$ lies inside $D_k(c)$.

\begin{theo}\cite[Proposition 11.13]{beskpi1}\label{theo:imagestandard} A pair $(x,(a_1,\ldots,a_k))\in E_n^\circ \times D_k(c)$ is \nit{compatible} if the cyclic support of $x$ contains $k$ points, and their respective multiplicities coincide with $(\ell_R(a_1),\ldots,\ell_R(a_k))$.\newline
The map $(\bbar{LL},\clbl)$ induces a bijection between $X/W$ and the set $E_n^\circ \boxdot D(c)$ of compatible pairs. This bijection induces a topology on the latter.
\end{theo}

By \cite[Remark 7.21]{beskpi1}, a path $\gamma$ in $E_n^\circ$ admits a unique lift $\ttilde{\gamma}$ in $X/W$ with fixed starting point provided that points are only merged and not unmerged in $\gamma(t)$ when $t$ increases. We want to understand how the cyclic label of $\ttilde{\gamma}(t)$ changes depending on $\gamma$.

Let $x\in X/W$ denote our chosen starting point for $\ttilde{\gamma}$. We begin by studying two cases where the homotopy class of $\ttilde{\gamma}(t)$ does not depend on the choice of $x$, at least under certain conditions:

Let $\tau:[0,1]\to \C$ be a path. If $x\in X/W$ is such that for all $t\in [0,1]$, we have $\tau(t)\notin \LLL(x)$, then we can define a path $\Tt_{\gamma,x}$ as follows: recall that we can identify $V/W$ with $\C^n$, and let us denote by $(v_1,\ldots,v_{n})$ the coordinates of $x$ under this identification. We set
\[\forall t\in [0,1],~ \Tt_{\gamma,x}(t)=(v_1,\ldots,v_{n-1},v_n+\tau(t)).\]
The condition that $\tau(t)\notin \LLL(x)$ for $t\in [0,1]$ ensures that $\Tt_{\tau,x}$ is indeed a path in $X/W$. By construction, we have $\Tt_{\gamma,x}=\Tt_{\gamma-\gamma(0),\Tt_{\gamma,x}(0)}$. In particular if $\tau(0)=0$, we have $\Tt_{\tau,x}(0)=x$. Moreover, by definition of $\LLL$, the multiset $\LLL(\Tt_{\tau,x}(t))$ is equal to $\LLL(x)-\tau(t)$, that is the translation of $\LLL(x)$ by $-\tau(t)$.

The image under $\LLL$ of a path of the form $\Tt_{\tau,x}$ is a translation by the continuous parameter $\tau(t)$. In particular such a path never merges points together and it admits a unique lift in $X/W$ starting from any point $x'$ such that $\LLL(x)=\LLL(x')$. The following lemma ensures that if $\gamma$ is a path of the form $\Tt_{\tau,x}$, then for $x'$ and $x$ close enough, the paths $\Tt_{\tau,x'}$ and $\Tt_{\tau,x}$ are homotopic.

\begin{lem}[\textbf{First Hurwitz rule}]
Let $\tau:[0,1]\to \C$ be a path and let $\eta:[0,1]\to X/W$ be a path such that
\begin{itemize}
\item For all $s,t\in [0,1]$, $\tau(t)\notin \LLL(\eta(s))$.
\item For all $s\in [0,1]$, no point of $\LLL(\eta(s))$ is directly above $\tau(0)$ or $\tau(1)$.
\end{itemize}
Then the paths $\Tt_{\tau,\eta(0)}$ and $\Tt_{\tau,\eta(1)}$ represent the same element of $B(W)$. 
\end{lem}
\begin{proof}
The second condition ensures that $\Tt_{\tau,\eta(0)}$ and $\Tt_{\tau,\eta(1)}$ are paths from $\Uu$ to $\Uu$. In particular they both represent elements of $B(W)$. Let now $s,t\in [0,1]$. Defining $H(s,t)=\Tt_{\tau,\eta(s)}(t)$ yields a homotopy between the paths $\Tt_{\tau,\eta(0)}$ and $\Tt_{\tau,\eta(1)}$. The first condition ensures that $H(s,t)$ remains in $X/W$ for all values of $s,t$, while the second condition ensures that the endpoints of the homotopy remain in $\Uu$ for all values of $s$.
\end{proof}

This Lemma is a generalization of \cite[Lemma 6.15]{beskpi1}, which corresponds to the case where $\tau$ is a path of the form $t\mapsto t\alpha$ for $\alpha>0$ and $t\in [0,1]$.

Another similar result states that the homotopy class of a circular tunnel $(x,L)$ does not depend on $x'$ close enough to $x$.

\begin{lem}[\textbf{Second Hurwitz Rule}] Let $x\in X/W$, and let $T:=(x,L)$ be a circular tunnel. Let also $\eta:[0,1]\to X/W$ be a path starting at $x$ and such that, for all $t\in [0,1]$, $(\eta(t),L)$ is a circular tunnel. For all $t\in [0,1]$, the elements in $B(W)$ represented by $(\eta(t),L)$ and $T$ are equal.
\end{lem}
\begin{proof}
Let $t\in [0,1]$. Defining $H(r,s):=e^{isL}\eta(rt)$ yields a homotopy between the paths associated with $T$ and to the circular tunnel $(\eta(t),L)$.
\end{proof}

This result, while easy to prove, is quite useful. For instance, it implies directly that altering the modulus of one (or several) point(s) of the cyclic support of some $x\in X/W$ does not affect the cyclic label provided that no two points are merged together.

Let now $x$ be in $X/W$ and let $(x_1,\cdots,x_k)$ be its cyclic support. Suppose that we swap two consecutive points $x_i$ and $x_{i+1}$ of $\LLL(x)$. Up to homotopy in $E_n^\circ$, there are two natural ways to do so, one point going ``farther'' than the other:

\spa
\begin{center}\begin{tikzpicture}
\draw [dashed] (0,0)--(0,2);

\node[draw,circle,inner sep=1pt,fill] (P0) at (00:1) {};
\node[draw,circle,inner sep=1pt,fill] (P1) at (60:1) {};
\node[draw,circle,inner sep=1pt,fill] (P2) at (120:1) {};
\node[draw,circle,inner sep=1pt,fill] (P3) at (180:1) {};
\node[draw,circle,inner sep=1pt,fill] (P4) at (240:1) {};
\node[draw,circle,inner sep=1pt,fill] (P5) at (300:1) {};

 \draw[->,>=latex] (P0) to[bend right] (P1);
 \draw[->,>=latex] (P1) to[bend right] (P0);

\begin{scope}[shift={(5,0)}]
\draw [dashed] (0,0)--(0,2);

\node[draw,circle,inner sep=1pt,fill] (P0) at (00:1) {};
\node[draw,circle,inner sep=1pt,fill] (P1) at (60:1) {};
\node[draw,circle,inner sep=1pt,fill] (P2) at (120:1) {};
\node[draw,circle,inner sep=1pt,fill] (P3) at (180:1) {};
\node[draw,circle,inner sep=1pt,fill] (P4) at (240:1) {};
\node[draw,circle,inner sep=1pt,fill] (P5) at (300:1) {};

 \draw[->,>=latex] (P0) to[bend left] (P1);
 \draw[->,>=latex] (P1) to[bend left] (P0);
\end{scope}
\end{tikzpicture}\end{center}

By \cite[Remark 7.21]{beskpi1}, these path both admits a unique lift in $X/W$. We denote these respective lifts by $\gamma_i^{+}$ and $\gamma_i^-$.

\begin{prop}\label{prop:cyclic_hurwitz_moves}
Let $x\in X/W$, with cyclic support $(x_1,\ldots,x_k)$ and cyclic label $(s_1,\ldots,s_k)$. The cyclic labels of $\gamma_i^+(1)$ and $\gamma_i^{-}(1)$ are given by
\[\clbl(\gamma_i^+(1))=(s_1,\ldots,s_{i-1},s_{i+1},s_{i}^{s_{i+1}},s_{i+2},\ldots,s_k),\]
\[\clbl(\gamma_i^-(1))=(s_1,\ldots,s_{i-1},{}^{s_i}s_{i+1},s_{i},s_{i+2},\ldots,s_k).\]
\end{prop}
\begin{proof}
The concatenation of $\gamma_i^-$ and $\gamma_i^+$ gives a homotopically trivial path from $x$ to itself. In particular the assertion about $\clbl(\gamma_i^-(1))$ follows from the assertion about $\clbl(\gamma_i^{+}(1))$. By construction of the cyclic label, we can assume that all the arguments of the points in $\LLL(x)$ are distinct. Let then $\theta$ be some angle strictly between the arguments of $x_{i-1}$ and $x_i$. By \cite[Lemma 11.1]{beskpi1}, one can replace $x$ by $e^{i\theta/h}x$ and assume that $i=1$. We will also denote $\gamma$ instead of $\gamma_1^+$ for readability.

Now if $\theta'$ denotes an angle strictly included between the arguments of $x_{2}$ and $x_3$, then the circular tunnels $(x,\theta')$ and $(\gamma(1),\theta')$ represent the same element in $B(W)$ by the second Hurwitz rule. That is the products of the first two terms of $\clbl(x)$ and of $\clbl(\gamma(1))$ are equal. The second Hurwitz rule also shows that the terms of $\clbl(x)$ and $\clbl(\gamma(1))$ are equal for $i>2$.

It remains to show that the first two terms of $\clbl(\gamma(1))$ are $(s_2,s_1^{s_2})$. Since we know that the product of the first two terms is equal to $s_1s_2$. It is sufficient to show that the first term of $\clbl(\gamma(1))$ is $s_2$. We will do so by using the Hurwitz rules but we first give some reductions.

As we mentioned above, changing the modulus of points of the cyclic support does not affect the cyclic label. Since we assumed that all points of the cyclic support have different arguments, we can then assume that all points of the cyclic support have the same modulus $R>0$. 

Let us write $\theta_k$ for the argument of $x_k$. Recall that we see $\theta_k$ as an element of $[\pi/2,5\pi/2[$. For a positive integer $n$ and $\theta\in [\pi/2,5\pi/2[$, we can consider the following paths 
\[A^-_{n,\theta}:t\mapsto \frac{3\pi}{2}+\frac{\theta-\frac{3\pi}{2}}{nt+1},~A^+_{n,\theta}:t\mapsto \frac{5\pi}{2}+\frac{\theta-\frac{5\pi}{2}}{nt+1},\]
both defined for $t\in [0,1]$. Now, consider the path $\rho_n$ in $E_n^\circ$ starting from $\LLL(x)$ and such that the $k$-th point of the cyclic support of $\rho_n(t)$ is $R\exp(iA^+_{n,\theta_k})$ if $k\in \{1,2\}$ and $R\exp(iA^-_{n,\theta_k})$ otherwise. In other words, the path $\rho_n$ consists in bringing the first two points of the cyclic support closer to the positive vertical half-line, while bringing all the other points of the cyclic support closer to the negative vertical half-line, as in the following example:
\spa
\begin{center}\begin{tikzpicture}[scale=0.8]
\draw [dashed] (0,0)--(0,2);

\node[draw,circle,inner sep=1pt,fill] (P0) at (00:1) {};
\node[draw,circle,inner sep=1pt,fill] (P1) at (60:1) {};
\node[draw,circle,inner sep=1pt,fill] (P2) at (120:1) {};
\node[draw,circle,inner sep=1pt,fill] (P3) at (180:1) {};
\node[draw,circle,inner sep=1pt,fill] (P4) at (240:1) {};
\node[draw,circle,inner sep=1pt,fill] (P5) at (300:1) {};

\draw[->,>=latex] (2,0) to (3,0);

\begin{scope}[shift={(5,0)}]
\draw [dashed] (0,0)--(0,2);

\node[draw,circle,inner sep=1pt,fill] (P0) at (435:1) {};
\node[draw,circle,inner sep=1pt,fill] (P1) at (445:1) {};
\node[draw,circle,inner sep=1pt,fill] (P2) at (245:1) {};
\node[draw,circle,inner sep=1pt,fill] (P3) at (255:1) {};
\node[draw,circle,inner sep=1pt,fill] (P4) at (265:1) {};
\node[draw,circle,inner sep=1pt,fill] (P5) at (275:1) {};
\end{scope}
\end{tikzpicture}\end{center}
The path $\rho_n$ preserves the ordering of the points of the cyclic support. Moreover, considering the unique lift $\ttilde{\rho_n}$ of $\rho_n$ in $X/W$ starting from $x$, we see by the second Hurwitz rule that the cyclic labels of $x$ and of $\ttilde{\rho_n}(1)$ are equal and that it is enough to prove the result for $\ttilde{\rho_n}(1)$. For all $\epsi>0$, choosing $n$ sufficiently large, we can assume that $\theta_1$ and $\theta_2$ lie in $]5\pi/2-\epsi,5\pi/2[$, while the other $\theta_k$ lie in $]3\pi/2-\epsi,3\pi/2+\epsi[$.

Now, let $x_1=a_1+ib_1$ and let $x_2=a_2+ib_2$. Since $x_1$ and $x_2$ have the same modulus, and for $\epsi>0$ small enough, we have $a_1<a_2$ and $b_2<b_1$. We set $\alpha:=(a_1+a_2)/2$ and $\beta:=(b_1+b_2)/2$.

Up to homotopy, the path $\gamma$ can be decomposed as the concatenation of two paths $\eta_1,\eta_2$, where the points $x_1,x_2$ move first vertically, then horizontally:
\begin{center}\begin{tikzpicture}[scale=0.8]
\node[draw,circle,inner sep=1pt,fill] (P0) at (0,1) {};
\node[draw,circle,inner sep=1pt,fill] (P1) at (0,0) {};
\node[draw,circle,inner sep=1pt,fill] (P2) at (3,0) {};
\node[draw,circle,inner sep=1pt,fill] (P3) at (3,1) {};

\node[left] at (P0) {$x_1$};
\node[right] at (P2) {$x_2$};

\draw[->,>=latex] (P0) to (P1);
\draw[->,>=latex] (P2) to (P3);
\draw[->,>=latex] (P1) to (P2);
\draw[->,>=latex] (P3) to (P0);

\node[above] at (1.5,1) {$\eta_2$};
\node[below] at (1.5,0) {$\eta_2$};

\node[left] at (0,0.5) {$\eta_1$};
\node[right] at (3,0.5) {$\eta_1$};
\end{tikzpicture}\end{center}
Let $\ttilde{\eta_1}$ be the lift of $\eta_1$ starting from $x$ and let $x'$ be its endpoint. Let also $\ttilde{\eta_2}$ be the lift of $\eta_2$ starting from $x'$. Since $\theta_1,\theta_2$ can be assumed to lie in $]5\pi/2-\pi/4,5\pi/2[$, the path $\eta_1$ preserves the ordering of the points of the cyclic support. Moreover by the second Hurwitz rule, the cyclic labels of $x$ and of $x'$ are equal and we just have to understand how $\ttilde{\eta_2}$ affects the cyclic label.

If $z,z'$ are two complex numbers, we will denote by $[z;z']$ the path in $\C$ defined by $t\mapsto z(1-t)+z't$ for $t\in [0,1]$. The first term of the cyclic label of $x'$ is represented by some circular tunnel homotopic to the path $\Tt_{[0,\alpha],x}$, while the product of the first two terms of the cyclic label of $x'$ is represented by some circular tunnel homotopic to the path $\Tt_{[0;2\alpha],x'}$. 

We have $\Tt_{[0;2\alpha],x'}=\Tt_{[0;\alpha],x'}\Tt_{[\alpha;2\alpha],x'}$. Since the path $\Tt_{[0;\alpha],x'}$ represents $s_1$ in $B(W)$, the path $\Tt_{[\alpha;2\alpha],x'}$ represents $s_2$ in $B(W)$. We will conclude by showing that the path $\Tt_{[\alpha;2\alpha],x'}$ is homotopic to the circular tunnel defining the first term of the cyclic label of $\gamma(1)$. We will do so by repetitive applications of the first Hurwitz rule, as in the following picture (where we represented $\tau$ instead of the motion $\Tt_{\tau,x'}$):

\begin{center}\begin{tikzpicture}[scale=0.4]
\node[draw,circle,inner sep=1pt,fill] (P1) at (1,1) {};
\node[draw,circle,inner sep=1pt,fill] (P3) at (4,2) {};

\draw[->,>=latex] (2.5,0) to (5,0);
\node at (5.5,1) {$\sim$};

\begin{scope}[shift={(6.5,0)}]
\node[draw,circle,inner sep=1pt,fill] (P1) at (1,1) {};
\node[draw,circle,inner sep=1pt,fill] (P3) at (4,2) {};

\draw[->,>=latex] (2.5,1.5) to (5,1.5);
\node at (5.5,1) {$\sim$};
\end{scope}

\begin{scope}[shift={(13,0)}]
\node[draw,circle,inner sep=1pt,fill] (P1) at (1,1) {};
\node[draw,circle,inner sep=1pt,fill] (P3) at (4,2) {};

\draw[-,>=latex] (0,0) to (0,1.5);
\draw[->,>=latex] (0,1.5) to (5,1.5);
\node at (5.5,1) {$\sim$};
\end{scope}

\begin{scope}[shift={(19.5,0)}]
\node[draw,circle,inner sep=1pt,fill] (P1) at (4,1) {};
\node[draw,circle,inner sep=1pt,fill] (P3) at (1,2) {};

\draw[-,>=latex] (0,0) to (0,1.5);
\draw[->,>=latex] (0,1.5) to (5,1.5);
\node at (5.5,1) {$\sim$};
\end{scope}

\begin{scope}[shift={(26,0)}]
\node[draw,circle,inner sep=1pt,fill] (P1) at (4,1) {};
\node[draw,circle,inner sep=1pt,fill] (P3) at (1,2) {};

\draw[-,>=latex] (0,0) to (0,1.5);
\draw[->,>=latex] (0,1.5) to (2.5,1.5);
\node at (5.5,1) {$\sim$};
\end{scope}

\begin{scope}[shift={(32.5,0)}]
\node[draw,circle,inner sep=1pt,fill] (P1) at (4,1) {};
\node[draw,circle,inner sep=1pt,fill] (P3) at (1,2) {};

\draw[->,>=latex] (0,0) to (2.5,0);
\end{scope}
\end{tikzpicture}\end{center}

By the first Hurwitz rule, $\Tt_{[\alpha;2\alpha],x'}$ is homotopic to $\Tt_{[\alpha+i\beta;2\alpha+i\beta],x'}$. Then, consider the path $\tau:=[0:i\beta]*[i\beta;\alpha+i\beta]$. The path $\Tt_{\tau,x'}$ is a path in $\Uu$ since no point of the cyclic support of $x'$ is directly above $\tau(t)$ for any $t$. Thus the path $\Tt_{\tau,x'}$ is homotopically trivial and the path $\Tt_{[\alpha+i\beta;2\alpha+i\beta],x'}$ is homotopic to
\[\Tt_{\tau,x'}*\Tt_{[\alpha+i\beta;2\alpha+i\beta],x'}=\Tt_{\tau*[\alpha+i\beta;2\alpha+i\beta],x'}=\Tt_{[0;i\beta]*[i\beta;2\alpha+i\beta],x'}.\]
Now, applying the first Hurwitz rule to $\tau':=[0;i\beta]*[i\beta;2\alpha+i\beta]$ and $\eta=\eta_2$, we obtain that the above path is homotopic to
\[\Tt_{[0;i\beta]*[i\beta;2\alpha+i\beta],\gamma(1)}\sim \Tt_{[0;i\beta]*[i\beta;\alpha+i\beta],\gamma(1)}.\]
Since no point of the cyclic support of $\gamma(1)$ lies below the segment $[i\beta;\alpha+i\beta]$, this path is homotopic to $\Tt_{[0;\alpha],\gamma(1)}$. As for $x'$, the path $\Tt_{[0;\alpha],\gamma(1)}$ is homotopic to a circular tunnel representing the first term of the cyclic label of $\gamma(1)$, as we wanted to show.
\end{proof}

By induction, we see that in general, moving a point of the cyclic support of some $x\in X/W$ may only affect the terms of the cyclic label corresponding to points of the support with lower modulus.

\subsection{Regular numbers and topological groupoids}\label{sec:2.4} We are now ready to give the topological definition of the Springer category associated with a regular number for a well-generated irreducible complex reflection group.

We keep the definitions and notation of the last section. Moreover, we fix a regular number $d$ for the group $W$, along with a $\zeta_d$-regular element $g$. We denote by $W_g:=C_{W}(g)$ the centralizer in $W$ of $g$ and we set
\[p:=\frac{d}{d\wedge h},~~q:=\frac{h}{d\wedge h},\]
where $d\wedge h$ is the gcd of $d$ and $h$. 

The categories we are going to consider are defined as fundamental groupoids with fat basepoints having several contractible connected components as in \cite[Appendix A]{beskpi1}. Define
\begin{align*}
D&:=\bigcup_{\zeta\in \mu_{p}} \zeta i\R_{\geqslant 0},\\
\Uu_{p}&:=\{x\in X/W~|~ \LLL(x)\cap  D=\varnothing\},\\
\Uu^{\mu_d}&:=\{x\in (X/W)^{\mu_d}~|~ \LLL(x)\cap D=\varnothing\}.
\end{align*}
The set $D$ is composed of $p$ half-lines starting at $0$. It cuts the plane in $p$ sectors $P_1,\ldots P_{p}$, labeled clockwise starting from the vertical halfline $i\R_+\subset D$. \newline In the following example, we have $p=3$:

\begin{center}\begin{tikzpicture}[scale=1]

\draw [dashed] (0,0)--(0,2);
\draw [dashed] (0,0)--(2*0.866,-1);
\draw [dashed] (0,0)--(-2*0.866,-1);

\node[right] at (0,2) {$D$};

\node at (0.866,0.5) {$P_1$};
\node at (0,-1) {$P_2$};
\node at (-0.866,0.5) {$P_3$};
\end{tikzpicture}\end{center}

\begin{lem}\cite[Lemma 11.22]{beskpi1} 
Let $x\in \Uu_p$ with cyclic label $(c_1,\ldots,c_k)$. For $i\in \intv{1,p}$, define $u_i$ as the product of the $c_j$ corresponding to points inside the sector $P_i$. The \nit{cyclic content} of $x$ is defined as $\cc_{p}(x):=(u_1,\ldots,u_{p})$.
\newline The map $\cc_{p}$ induces a bijection between the connected components of $\Uu_{p}$ (resp. $\Uu^{\mu_d}$) and $D_{p}(c)$ (resp. $D_{p}^{q}(c)$). Furthermore the connected components of $\Uu_{p}$ (resp. $\Uu^{\mu_d}$) are contractible.
\end{lem}

Since the connected components of $\Uu_{p}$ and $\Uu^{\mu_d}$ are contractible, we can use them as fat basepoints for groupoids in the sense of \cite[Definition A.4]{beskpi1}.

\begin{definition}\cite[Definition 11.23]{beskpi1} The \nit{Springer groupoid} associated with $W$ and $d$ is the groupoid
\[B_{p}^{q}(W)=\pi_1((X/W)^{\mu_d},\Uu^{\mu_d}).\]
The functoriality of $\pi_1$ gives a natural functor $B_{p}^{q}(W)\to B(W)$.
\end{definition}
Note that $B_p^q(W)$ is equivalent to the braid group of $W_g$ by construction. In particular, it is a connected groupoid. 

We now show that the Springer groupoid $B_p^q(W)$ coincides with the groupoid $\Gg$ of $(p,q)$-periodic elements of the dual braid monoid $M(c)$ in the sense of Definition \ref{def:periodic_elements_category}. Let $\Cc:=(M(c))_p^q$ denote the category of $(p,q)$-periodic elements of $M(c)$, so that $\Gg$ is the enveloping groupoid of $\Cc$.

Combining Remark \ref{rem:divided_set_interval_monoid} with Proposition \ref{prop:descriptiondivset}, we obtain that the simple morphisms of $\Cc$ are particular couples of elements belonging to the interval $I_c$ in $W$. Let $s:=(a,b)$ be such a simple morphism. By Theorem \ref{theo:imagestandard}, there is a unique element $x_s\in \Uu^{\mu_d}$ such that $\clbl(x)=s$ and $\LLL(x_s)$ consists of the points $\exp\left(i\pi(\frac{1}{2}-\frac{2j+1}{2p})\right)$ such that the $j$-th term of $s$ is nontrivial. The circular tunnel $(x_s,\frac{\pi}{ph})$ then defines an element of $B_{p}^q(W)$ which we denote by $b_s$.

\begin{theo}\label{theo:groupoidtressequivgroupoiddiv}\cite[Theorem 11.28]{beskpi1} The map $s\mapsto b_s$ extends to a groupoid isomorphism $\Gg\to B_{p}^q(W)$.
\end{theo}

Using this result, we will now also call \nit{Springer groupoid} the groupoid of $(p,q)$-periodic elements of the dual braid monoid $M(c)$. Moreover, the corresponding category of $(p,q)$-periodic elements will be called the \nit{Springer category}.

We have the following diagram of functors
\[\xymatrix{\Cc \ar@{^(->}[d] \ar[r]^-{\pi_p} & M(c) \ar@{^(->}[d]\\ 
\Gg \ar[d]^-{\simeq} \ar[r] & G(c) \ar[d]^-{\simeq}\\ 
B_{p}^{q}(W) \ar[r] & B(W) }\]

As $\bbf{c}^h\in M(c)$ represents the full-twist in $B(W)$, a $d$-th root of the full-twist is a $(d,h)$-periodic element of $M(c)$. Such an element is conjugate to a $(p,q)$-periodic element of $M(c)$ (see Remark \ref{rem:pqcoprime}). Let $u$ be an object of $B_{p}^{q}(W)$, that is an object of $M(c)_{p}^{q}$ or an element of $D_{p}^{q}(c)$. By Theorem \ref{theo:categoriediviseeetelementsreguliers}, the morphism $B_{p}^{q}(W)\to B(W)$ sends the automorphism group $B_{p}^{q}(W)(u,u)$ to the centralizer in $B(W)$ of some $(p,q)$-periodic element of $M(c)$, in particular a $d$-th root of the full twist in $B(W)$. 

\begin{lem}\label{futw}
Let $u$ be an object of $B_{p}^{q}(W)$. The group isomorphisms
\[B(W_g)\simeq B_{p}^{q}(W)(u,u)\simeq \Gg(u,u)\]
map the full-twist in $B(W_g)$ to $\Delta_p^{ph}(u)$, where $\Delta_p$ denotes the Garside map of the Springer category $\Cc$.
\end{lem}
\begin{proof}
The morphism $\Delta_p(u)$ in $\Cc=(M(c))_{p}^{q}$ corresponds to a rotation of angle $\frac{2\pi}{ph}$ in $X/W$, while the full-twist corresponds to a rotation of angle $2\pi$.
\end{proof}

\section{Springer categories which are monoids}\label{sec:categories_which_are_monoids} 
In this section, we fix a finite dimensional complex vector space $V$ with $n:=\dim(V)$, along with a well-generated irreducible complex reflection group $W<\GL(V)$. We denote by $h$ the highest degree of $W$. We also fix a regular number $d$ for $W$, which we will assume to divide $h$ after Lemma \ref{lem:caractérisation_springer_groupoid_one_object}. We will then set $q:=h/d$. Lastly, we fix some Coxeter element $c$ in $W$. Exceptionally, we denote by $I_c(W)$ the interval associated with $c$ in $W$ with respect to the reflection length instead of $I_c$.

Even though Springer groupoids have more than one object in general, there are several cases in which they only have one object. For instance, if $d=1$, then we have $(X/W)^{\mu_d}=X/W$, and the Springer category attached to $W$ and $d$ is simply the dual braid monoid $M(c)$. We would then like to know when do Springer categories have only one object, and what are the monoids that arise in these cases. 

The answer to the first question is given by the following lemma:

\begin{lem}\label{lem:caractérisation_springer_groupoid_one_object}
The Springer category associated with $W$ and $d$ is a monoid if and only if $d$ divides $h$.
\end{lem}
\begin{proof}
We write $p:=d/d\wedge h$ and $q:=h/d\wedge h$. Let $\eta,\mu$ be positive integers such that $p\eta-q\mu=1$. Let also $\Cc:=(M(c))^q_p$ denote the considered Springer category. If $p=1$, then $\Cc=M(c)^{\phi^q}$ is by definition a submonoid of $M(c)$ (where $\phi$ denotes the Garside automorphism of $M(c)$, that is, conjugation by $\bbf{c}$). In particular it has only one object. Conversely, assume that $\Cc$ has one object. Following Definition \ref{def:periodic_elements_category} 
, the object set of $\Cc$ is in bijection with the set 
\[\{\bbf{u}\in M(c)^{\phi^q}~|~ \bbf{u}\bbf{u}^{\bbf{c}^\eta}\cdots \bbf{u}^{\bbf{c}^{(p-1)\eta}}=\bbf{c}\}.\]
In particular, if $\bbf{u}$ is an object, then so is $\bbf{u}^{\bbf{c}^{\eta}}$. If $\Cc$ has one object, we then have that $\bbf{u}$ commutes with $\bbf{c}^\eta$. Since $\bbf{u}$ commutes with $\bbf{c}^q$ by definition, we obtain that $\bbf{u}$ commutes with $\bbf{c}$, as $\eta$ and $q$ are coprime. By \cite[Lemma 12.2 and Theorem 12.3]{beskpi1}, the only nontrivial simple element of $M(c)$ which commutes with $\bbf{c}$ is $\bbf{c}$ itself. We then have $\bbf{u}=\bbf{c}$ and thus $p=1$.
\end{proof}

We can now assume that $d$ divides $h$, and we set $q:=h/d$. In this case, the Springer category is the category of $(1,q)$-periodic elements of $M(c)$, that is the submonoid $M(c)^{\phi^q}$ of $M(c)$ consisting of elements which commute with $\bbf{c}^q$. Its set of simple elements corresponds to the set $I_c(W)^{c^q}$ of simple elements which commute with $c^q$ in $W$. We have by Theorem \ref{theo:groupoidtressequivgroupoiddiv} that $M(c)^{q}$ is a Garside monoid for the braid group of the centralizer $W_{c^q}:=C_W(c^q)$ (note that $c^q$ is a $\zeta_d$-regular element of $W$). We will show that $M(c)^{\phi^q}$ is actually isomorphic to the dual braid monoid for the group $W_{c^q}$. First, we show that $W_{c^q}$ is well-generated.

\begin{lem}
The group $W_{c^q}=C_W(c^q)$ (acting on the eigenspace $\Ker(c^q-\zeta_h^q)$) is well-generated.
\end{lem}
\begin{proof}
By \cite[Theorem 2.4]{beskpi1}, we know that an irreducible complex reflection group is well-generated if and only if the sum of its $i$-th degree (in increasing order) with its $i$-th codegree (in decreasing order) is constant and equal to the highest degree. We denote the degrees (resp. codegrees) of $W$ by $d_1,\ldots,d_n$ in increasing order (resp. $d_1^*,\ldots,d_n^*$ in decreasing order). We know by Theorem \ref{theo:degrees_of_regular_centralizer} that the degrees (resp. codegrees) of $W_{c^q}$ are precisely the $d_i$ (resp. the $d_i^*$) which are divisible by $d$. Since $d$ divides $h=d_i+d_i^*$ for all $i\in \intv{1,n}$, we have that $d$ divides $d_i$ if and only if it divides $d_i^*$. Thus, if the degrees of $W_{c^q}$ are $d_{i_1},\ldots,d_{i_k}$, then the codegrees of $W_{c^q}$ are $d_{i_1}^*,\ldots,d_{i_k}^*$. In particular for $j\in \intv{1,k}$, we have $d_{i_j}+d_{i_j}^*=h=d_{i_k}$ and $W_{c^q}$ is well-generated.
\end{proof}

Since $c$ commutes with $c^q$, we have $c\in W_{c^q}$ and, as $\Ker(c-\zeta_h)\subset \Ker(c^q-\zeta_h^q)$, it is a Coxeter element of $W_{c^q}$. By the above lemma, we can consider the interval $I_c(W_{c^q})$ associated with $c$ in $W_{c^q}$ for the reflection length in $W_{c^q}$. This section is devoted to the proof of the following theorem.

\begin{theo}\label{theo:intervalles_identiques_p=1}
The posets $I_c(W)^{c^q}=I_c(W_{c^q})$ are equal. That is, we have $I_c(W)^{c^q}=I_c(W_{c^q})$ as subsets of $W_{c^q}$ and, for $s,t\in I_c(W_{c^q})$, we have $s\preceq t$ in $W$ if and only if $s\preceq t$ in $W_{c^q}$, for the respective reflection lengths of $W$ and $W_{c^q}$.
\end{theo}

\begin{cor}\label{cor:dual_braid_relations_for_centralizer_of_c^q}
The monoid $M(c)^{\phi^q}$ is isomorphic to the dual braid monoid $M'(c)$ associated with $W_{c^q}$.
\end{cor}
The corollary comes from the fact that the defining presentation of an interval monoid only depends on the poset structure of the associated interval $(I,\preceq)$: the defining relations are all those of the form $s(s^{-1}t)=t$ for $s,t\in I$ with $s\preceq t$.

The proof of Theorem \ref{theo:intervalles_identiques_p=1} will ultimately rely on a case by case analysis, but we can do some easy reductions:

\begin{itemize}
\item If $d$ divides all the degrees of $W$, then $c^q$ is central in $W$ and we have $W_{c^q}=W$, Theorem \ref{theo:intervalles_identiques_p=1} is obvious in this case. 
\item If $h$ is the only degree of $W$ divisible by $d$, then $W_{c^q}$ is a complex reflection group of rank $1$: it is cyclic and equal to $\langle c\rangle$. We have $I_c(W_{c^q})=\{\Id, c\}$. On the other hand, any $s\in I_c^q$ lies in $W_{c^q}=\langle c\rangle$. We then have $I_c^q=\{\Id, c\}$ by \cite[Lemma 12.2]{beskpi1}. The poset structure is induced by $\Id\preceq c$ in both cases and Theorem \ref{theo:intervalles_identiques_p=1} holds.
\end{itemize}

Note that these two extreme cases are sufficient to prove Theorem \ref{theo:intervalles_identiques_p=1} when $W$ has rank $2$. We can now assume that $n\geqslant 3$. In this case, using the classification of irreducible complex reflection groups \cite[Theorem 8.29]{lehrertaylor}, along with \cite[Theorem 2.4]{beskpi1}, we are in exactly one of the following cases:

\begin{itemize}
\item $W$ is one of the exceptional groups $G_i$ with $i\in \intv{23,37}\setminus \{31\}$. These cases are handled by computer.
\item $W\simeq G(m,1,n)$ for $m\geqslant 2$. This is the group of $n\times n$ monomial matrices whose nonzero entries lie in $\mu_m$.
\item $W\simeq G(e,e,n)$ for $e\geqslant 2$. This group is the kernel of the morphism $G(m,1,n)\to \mu_m$ which sends $w$ to the product of its nonzero entries (i.e. to $(-1)^n\det(w)$).
\item $n\geqslant 4$ and $W$ is the group $\S_{n+1}$ acting by permutation on the hyperplane 
\[H:=\{(x_1,\ldots,x_{n+1})\in \C^{n+1}~|~ x_1+\cdots +x_n=0\}.\]
The group $G(1,1,n+1)=\S_{n+1}$ acting on the whole space $\C^{n+1}$ is not irreducible (it admits an invariant line). However, the sets of reflections of $W$ and of $G(1,1,n+1)$ are equal (they are the transpositions in $\S_{n+1}$), and their sets of Coxeter elements are the same (they are the $n+1$-cycles in $\S_{n+1}$). Thus the interval monoids given by $W$ and $G(1,1,n+1)$ are identical. It is then sufficient to prove the statement of Theorem \ref{theo:intervalles_identiques_p=1} for the latter group.
\end{itemize}

\subsection{The case $W=G(1,1,n+1)$ for $n\geqslant 5$}\label{sec:G(1,1,n)}
Let $\ell_W$ denote the reflection length in $W$.  Our approach is largely modeled on that of \cite[Section 3 and Section 4]{bradywatt2}, which covers the case where $n$ is even and $q=\frac{n}{2}$.
\begin{lem}\label{3.13}Let $w=c_1\cdots c_k$ be a product of disjoint cycles in $W$.
We have
\[\ell_W(w)=\sum_{i=1}^k\ell_W(c_i)=\sum_{i=1}^k (\ell(c_i)-1),\]
where $\ell(c)$ denotes the size of the support of $c$.
\end{lem}
\begin{proof}
The group $W$ is a complexified real reflection group by \cite[Lemma IX.2.19]{ddgkm}. The reflection length $\ell_W(w)$ is then given by the codimension of the fixed space of $w$ acting on $\C^n$. Since the $c_i$ all have disjoint supports, the first equality is immediate. For the second equality, let $c=(i_1~\cdots~i_k)$ be a $k$ cycle in $W$. The fixed space of $c$ acting on $\C^n$ is given by $\{(x_1,\ldots,x_n) \in \C^n~|~x_{i_1}=x_{i_2}=\cdots=x_{i_k}\}$ and has codimension $k-1$.
\end{proof}
The highest degree of $W$ is $n$. We label the canonical basis of $\C^n$ in the following way
\[\{e(0,1),e(0,2),\ldots,e(0,q),e(1,1),\ldots,e(d-1,q)\}.\]
A guiding idea is to think of $e(j,i)$ as $\zeta_d^je_i$ in the vector space $\C^q$ (with canonical basis $e_1,\ldots,e_q)$. We consider the following element of $W$:
\[c(1,1,n):=(e(0,1)~e(0,2)~\cdots~e(0,q)~e(1,1)~\cdots~e(d-1,q)).\]
It is a Coxeter element of $W$ as the vector $\sum_{i=1}^q\sum_{j=0}^{d-1} \zeta_n^{n-i-qj}e(j,i)$ is a $\zeta_n$-eigenvector for $c(1,1,n)$ which is regular. We have
\[c(1,1,n)^q=(e(0,1)~e(1,1)~\cdots~e(d-1,1))(e(0,2)~\cdots~e(d-1,2))\cdots(e(0,q)~\cdots~e(d-1,q)).\]

\begin{lem}\label{lem:cest_gd1q}The eigenspace $V:=\Ker(c(1,1,n)^q-\zeta_d\Id)$ has dimension $q$ and admits a basis $v:=\{v_1,\ldots,v_q\}$, where
\[v_i:=\sum_{j=0}^{d-1}\zeta_d^{d-j}e(j,i).\]
The isomorphism of vector spaces $V\simeq \C^q$ given by the basis $v$ induces an isomorphism between $W_{c(1,1,n)^q}$ and $G(d,1,q)$.
\end{lem}
\begin{proof}
The degrees of $W_{c(1,1,n)^q}$ acting on $V$ are the degrees of $G(1,1,n)$ which are divided by $d$, that is $d,2d,\ldots,dq$. In particular, $W_{c(1,1,n)^q}$ has rank $q$ and $V$ has dimension $q$. As the $v_i$ are clearly linearly independent, $v$ is a basis for $V$. 

The groups $W_{c(1,1,n)^q}$ and $G(d,1,q)$ share the same degrees (cf. \cite[Table 2]{bmr}), and the same cardinality (since the order of a complex reflection group is the product of its degrees). We then only have to show that the image of $W_{c(1,1,n)^q}$ in $\GL_q(\C)$ contains a generating set of $G(d,1,q)$, like the one given in \cite[Section 2.7]{lehrertaylor}. This is a direct check:
\begin{itemize}
\item The cycle $(e(0,1)~e(1,1)\cdots e(d-1,1))$ acts trivially on $v_j$ for $j\neq 1$ and sends $v_1$ to $\zeta_d v_1$.
\item For $i\in \intv{1,q-1}$, the permutation $(e(0,i)~e(0,i+1))\cdots(e(d-1,i)~e(d-1,i+1))$ swaps $v_i$ and $v_{i+1}$.
\end{itemize}
\end{proof}
From now on, we set $W'=W_{c(1,1,n)^q}$ and $\ell_{W'}$ the reflection length for elements of $W'$ regarding the reflections of $W'$ acting on $V$. The action of $c(1,1,n)$ on $V$ is given (in the basis $v$) by the matrix
\[c(d,1,q):=\matrix{0&  & & \zeta_d \\ 1 &  \ddots & & \\ &  \ddots &0 & \\ & & 1 & 0},\]
which is indeed a Coxeter element for the group $G(d,1,q)$ (see Section \ref{sec:G(m,1,n)}). The group $G(d,1,q)$ is endowed with a character $\chi$, sending a monomial matrix to the product of its nonzero entries. The isomorphism $W'\simeq G(d,1,q)$ of Lemma \ref{lem:cest_gd1q} allows us to define $\chi$ on $W'$.

Let $c:=(e(j_1,i_1)~\cdots~e(j_k,i_k))\in W$ be a cycle. We set 
\[c^{(1)}:=c^{c(1,1,n)^q}=(e(j_1-1,i_1)~\cdots~e(j_k-1,i_k)).\]

\begin{prop}\label{prop:elements_de_gd1q_dans_sn}
An element $\sigma\in W$ lies in $W'$ if and only if it can be written as decomposition as a product of disjoint cycles of the form
\[\sigma=c_1c_1^{(1)}\cdots c_1^{(d-1)}c_2c_2^{(1)}\cdots c_2^{(d-1)}\cdots c_ac_a^{(1)}\cdots c_a^{(d-1)} \gamma_1\cdots \gamma_b,\]
where $\gamma_i^{(1)}=\gamma_i$ for $i\in \intv{1,b}$.
\end{prop}
\begin{proof}
Our proof is an adaptation of the proof of \cite[Proposition 3.1]{bradywatt2}, which deals with the case $d=2$.
First, it is clear that elements of the given form lie in $W'$. Conversely, let $\sigma=c_1\cdots c_r$ be a product of disjoint cycles in $\S_n$. We have that $c(1,1,n)^q$ centralizes $\sigma$ if and only if $c_1\cdots c_r=c_1^{(1)}\cdots c_r^{(1)}$. By uniqueness (up to reordering) of cycle decompositions in $\S_n$, for each $i$ either $c_i=c_j^{(1)}$ for some $j\neq i$ or else $c_i=c_i^{(1)}$. An immediate induction then gives that $\sigma$ has the required decomposition.
\end{proof}

We note in particular that the cycles $c_i\cdots c_i^{(d-1)}$ in the decomposition of Proposition \ref{prop:elements_de_gd1q_dans_sn} are disjoint. 

\begin{definition}
Let $c\in W$ be a cycle such that $c,c^{(1)},\ldots,c^{(d-1)}$ are disjoint. The product $cc^{(1)}\cdots c^{(d-1)}$ will be denoted by $\ttilde{c}$ and called a \nit{saturated cycle}. If $\gamma\in W$ is a cycle lying in $W'$, we say that $\gamma$ is a \nit{balanced cycle}.
\end{definition}

Proposition \ref{prop:elements_de_gd1q_dans_sn} states that elements of $W'$ are the products of disjoint saturated cycles and disjoint balanced cycles.

\begin{lem}\label{lem:saturated_cycles_and_balanced_cycles}
Let $\ttilde{c}$ be a saturated cycle in $W'$. We have $\chi(\ttilde{c})=1$ and $\ell_{W'}(c)=\ell(c)-1$.\newline Let $\gamma$ be a balanced cycle in $W'$. We have $\chi(\gamma)=\zeta_d$ and $\ell_{W'}(\gamma)=\ell(\gamma)/d$.
\end{lem}
\begin{proof}Let $c=(e(j_1,i_1)~\cdots~e(j_k,i_k))$ be a cycle such that $c,\ldots,c^{(d-1)}$ are all disjoint. By assumption, $i_1,\ldots,i_k$ are all distinct. One readily sees that $\ttilde{c}$ acts on $V$ by 
\[\begin{cases} \ttilde{c}.v_{i_m}=\zeta_d^{j_{m+1}-j_m}v_{i_{m+1}}&\forall m\in \intv{1,k-1},\\ \ttilde{c}.v_{i_k}=\zeta_d^{j_1-j_k}v_{i_1}.\end{cases}\]
In particular, we have $\chi(\ttilde{c})=\zeta_d^{j_2-j_1}\cdots\zeta_d^{j_k-j_{k-1}}\zeta_d^{j_1-j_k}=1$. The fixed space of $\ttilde{c}$ acting on $V$ is generated by all the $v_i$ with $i\notin\{i_1,\ldots,i_k\}$ and $v_{i_1}+\cdots+v_{i_k}$. Thus $\ell_{W'}(\ttilde{c})=k-1=\ell(c)-1$ by Lemma \ref{lem:cest_gd1q} and \cite[Theorem 2.1]{reflection_length_in_monomial}.

Let $\gamma$ be a balanced cycle. It can be written as 
\[\gamma=(e(0,i_1)~e(j_2,i_2)~\cdots~e(j_k,i_k)~e(1,i_1)~\cdots~e(j_k+d-1,i_k)).\]
The action of $\gamma$ on $V$ is then given by
\[\begin{cases} \gamma.v_{i_1}=\zeta_d^{j_2}v_{i_2},\\ \gamma.v_{i_m}=\zeta_d^{j_{m+1}-j_m}v_{i_{m+1}}&\forall m\in \intv{2,k-1},\\ \gamma.v_{i_k}=\zeta_d^{1-j_k}v_1.\end{cases}\]
In particular we have $\chi(\gamma)=\zeta_d^{j_2}\cdots \zeta_d^{j_k-j_{k-1}}\zeta_d^{1-k}=\zeta_d$. The fixed space of $\gamma$ acting on $V$ is generated by all the $v_i$ with $i\notin\{i_1,\ldots,i_k\}$. Thus $\ell_{W'}(\gamma)=k=\ell(\gamma)/d$ again by Lemma \ref{lem:cest_gd1q} and \cite[Theorem 2.1]{reflection_length_in_monomial}.
\end{proof}

By combining Lemma \ref{lem:saturated_cycles_and_balanced_cycles} and Lemma \ref{3.13}, we obtain the following proposition.
\begin{prop}\label{lem:longueurs_reflections}
Let $\sigma=\ttilde{c_1}\cdots\ttilde{c_a}\gamma_1\cdots\gamma_b\in W'$. The reflection length of $\sigma$ in $W$ and $W'$ are given by
\[\ell_W(\sigma)=\sum_{i=1}^a d(\ell(c_i)-1)+\sum_{j=1}^b (\ell(\gamma_i)-1) \text{~and~}\ell_{W'}(\sigma)=\sum_{i=1}^a (\ell(c_i)-1)+\sum_{j=1}^b \frac{\ell(\gamma_i)}{d}.\]
In particular, we have $\ell_W(\sigma)+b=d\ell_{W'}(\sigma)$.
\end{prop}
This relation between reflection lengths in $W$ and in $W'$ will be the key element in our proof of Theorem \ref{theo:intervalles_identiques_p=1} for the case $W=G(1,1,n)$.

\begin{lem}\label{lem:division_matrice_diagonale}
Let $D$ be a nontrivial diagonal matrix in $G(d,1,q)$. We have $D\preceq c(d,1,q)$ in $G(d,1,q)$ if and only if $D$ is a diagonal reflection $s$ with $\chi(s)=\zeta_d$.
\end{lem}
\begin{proof}
By \cite[Theorem 2.1]{reflection_length_in_monomial}, the reflection length $\ell_{G(d,1,q)}(D)$ of $D$ in $G(d,1,q)$ is the number of nontrivial diagonal entries of $D$. The underlying permutation of $D^{-1}c(d,1,q)$ is an $n$-cycle. Thus, again by \cite[Theorem 2.1]{reflection_length_in_monomial}, we have
\[\ell_{G(d,1,q)}(D^{-1}c(d,1,q))=\begin{cases} n-1&\text{if }\chi(D)=\chi(c(d,1,q))=\zeta_d,\\ n&\text{otherwise.}\end{cases}\]
If $D$ is nontrivial (i.e. if its reflection length is nonzero), we get $D\preceq c(d,1,q)$ if and only if $\chi(D)=\zeta_d$ and $D$ has exactly one nontrivial diagonal entry.
\end{proof}

\begin{prop}\label{prop:keylemma}
Let $\sigma=\ttilde{c_1}\cdots \ttilde{c_a}\gamma_1\cdots\gamma_b \in W'$. If $\sigma\in I_{c(1,1,n)}(W')$ or if $\sigma \in I_{c(1,1,n)}(W)$, then $b\in \{0,1\}$ and $\chi(\sigma)=\zeta_d^{b}$.
\end{prop}
\begin{proof}Let $\gamma=(e(0,i_1)~e(j_2,i_2)~\cdots~e(j_k,i_k)~e(1,i_1)~\cdots~e(j_k+d-1,i_k))$ be a balanced cycle. The diagonal reflection of $W'$ sending $v_i$ to $\zeta_d v_i$ is given by $s_i:=(e(0,i)\cdots e(d-1,i))$. We have that $s_{i_1}^{-1}\gamma=(e(0,i_1)\cdots e(j_k,i_k))\cdots (e(d-1,i_1)\cdots e(j_k+d-1,i_k))$ is a saturated cycle. Proposition \ref{lem:longueurs_reflections} then gives that $\ell_{W'}(s_{i_1})+\ell_{W'}(s_{i_1}^{-1}\gamma)=1+k-1=k=\ell_{W'}(\gamma)$ and $s_{i_1}\preceq \gamma$ in $W'$.

Assume that $\sigma\in I_{c(1,1,n)}(W')$ and $b>1$. Since $\gamma_1$ and $\gamma_2$ have disjoint support, we deduce that there are two diagonal reflections $s,s'$ of $W'$ such that $s\preceq \gamma_1,s'\preceq \gamma_2$, $ss'\neq \Id$ and $\chi(ss')=\zeta_d^2$. We then have $ss'\preceq c(1,1,n)$ in $W'$ which contradicts Lemma \ref{lem:division_matrice_diagonale}.

Assume now that $\sigma\in I_{c(1,1,n)}(W)$. By \cite[Proposition IX.2.7]{ddgkm}, the set theoretical partition of the set $\mu_n$ induced by $\sigma$ is noncrossing (in the sense of \cite[Section 1.2]{bescor}). If $b>1$, then the orbits of $\gamma_1$ and $\gamma_2$ induce two parts of $\mu_n$ whose convex hull contains $0$. Thus the partition of $\mu_n$ induced by $\sigma$ is crossing and $\sigma\notin I_{c(1,1,n)}(W)$.

The assumption on $\chi(\sigma)$ is a direct consequence of Lemma \ref{lem:saturated_cycles_and_balanced_cycles}.
\end{proof}

\begin{prop}
The two posets $I_{c(1,1,n)}(W)^{c(1,1,n)^q}$ and $I_{c(1,1,n)}(W')$ are equal.
\end{prop}
\begin{proof}
First, we show that the two sets $I_{c(1,1,n)}(W)^{c(1,1,n)^q}$ and $I_{c(1,1,n)}(W')$ are equal. Let $\sigma\in I_{c(1,1,n)}(W')$. By Proposition \ref{prop:keylemma}, we have $\chi(\sigma)\in \{1,\zeta_d\}$ and $\chi(\sigma^{-1}c(1,1,n))=\chi(\sigma)^{-1}\zeta_d$. Proposition \ref{lem:longueurs_reflections} then gives
\begin{align*}
\ell_W(\sigma)+\ell_W(\sigma^{-1}c(1,1,n))&=d\ell_{W_{c(1,1,n)^q}}(\sigma)+d\ell_{W_{c(1,1,n)^q}}(\sigma^{-1}c(1,1,n))-1\\
&=d(\ell_{W_{c(1,1,n)^q}}(\sigma)+\ell_{W_{c(1,1,n)^q}}(\sigma^{-1}c(1,1,n)))-1\\
&=d\ell_{W_{c(1,1,n)^q}}(c(1,1,n)^q)-1\\
&=\ell_W(c(1,1,n))
\end{align*}
and $\sigma\in I_{c(1,1,n)}(W)$. The same reasoning gives $I_{c(1,1,n)}(W)^{c(1,1,n)^q}\subset I_{c(1,1,n)}(W')$.

Let now $\sigma,\tau\in I_{c(1,1,n)}(W')$, we have to show that $\sigma\preceq \tau$ in $W'$ if and only if $\sigma\preceq \tau$ in $W$. We consider four cases
\begin{itemize}
\item $\chi(\sigma)=1$ and $\chi(\tau)=1$. We have $\ell_W(\sigma)=d\ell_{W'}(\sigma)$ and $\ell_W(\tau)=d\ell_{W'}(\tau)$. If $\sigma\preceq \tau$ in $W$ (resp. in $W'$), then $\sigma^{-1}\tau\in I_{c(1,1,n)}(W')$ is such that $\chi(\sigma^{-1}\tau)=1$. We then have $\ell_W(\sigma^{-1}\tau)=d\ell_{W'}(\sigma^{-1}\tau)$ and $\sigma\preceq \tau$ in $W'$ (resp. in $W$).
\item $\chi(\sigma)=\zeta_d$ and $\chi(\tau)=1$. We cannot have $\sigma\preceq \tau$ in either $W$ or $W'$ since this would imply that $\sigma^{-1}\tau$ is an element of $I_{c(1,1,n)}(W')$ with $\chi(\sigma^{-1}\tau)=\zeta_d^{-1}$ (we can assume $\zeta_d\neq -1$ since the case $d=2$ is known by \cite[Lemma 4.8]{bradywatt2}).
\item $\chi(\sigma)=1$ and $\chi(\tau)=\zeta_d$. We have $\ell_W(\sigma)=d\ell_{W'}(\sigma)$ and $\ell_W(\tau)+1=d\ell_{W'}(\tau)$. If $\sigma\preceq \tau$ in $W$ (resp. in $W'$), then $\sigma^{-1}\tau\in I_{c(1,1,n)}(W')$ is such that $\chi(\sigma^{-1}\tau)=\zeta_d$. We then have $\ell_W(\sigma^{-1}\tau)+1=d\ell_{W'}(\sigma^{-1}\tau)$ and $\sigma\preceq \tau$ in $W'$ (resp. in $W$).
\item $\chi(\sigma)=\zeta_d$ and $\chi(\tau)=\zeta_d$. We have $\ell_W(\sigma)+1=d\ell_{W'}(\sigma)$ and $\ell_W(\tau)+1=d\ell_{W'}(\tau)$. If $\sigma\preceq \tau$ in $W$ (resp. in $W'$), then $\sigma^{-1}\tau\in I_{c(1,1,n)}(W')$ is such that $\chi(\sigma^{-1}\tau)=1$. We then have $\ell_W(\sigma^{-1}\tau)=d\ell_{W'}(\sigma^{-1}\tau)$ and $\sigma\preceq \tau$ in $W'$ (resp. in $W$).
\end{itemize}
\end{proof}

This closes the proof of Theorem \ref{theo:intervalles_identiques_p=1} in this case. We have the following corollary, which substantiates \cite[Remark at the end of Section 7]{bescor}, for which we could not find a proof in the literature.

\begin{cor}\label{cor:partitions_e1n}
The map from $G(d,1,q)$ to the set of partitions of the regular $n$-gon that sends every element of $G(d,1,q)$ to the partition given by the orbits of its image in $\S_n$ induces a poset isomorphism from $I_{c(d,1,q)}(G(d,1,q))$ to $\NCP(d,1,q):=\NCP(1,1,n)^{\mu_d}$, where $\NCP(1,1,n)$ denotes the set of noncrossing partitions of the regular $n$-gon (in the sense of \cite[Section 1.2]{bescor}).
\end{cor}
\begin{proof}
The embedding $G(d,1,q)\to G(1,1,n)$ as the centralizer of $c(1,1,n)^q$ identifies the posets $I_{c(d,1,q)}(G(d,1,q))$ and $I_{c(1,1,n)}(W')=I_{c(1,1,n)}(W)^{c(1,1,q)}$. On the other hand, the map from $G(1,1,n)$ to the set of partitions of the regular $n$-gon that sends every element of $G(1,1,n)$ to the partition given by its orbits induces a poset isomorphism from $I_{c(1,1,n)}(G(1,1,n))$ to $\NCP(1,1,n)$ \cite[Proposition IX.2.7]{ddgkm}. Under this isomorphism, the action of $c(1,1,n)$ by conjugation corresponds to a rotation of angle $\frac{2\pi}{n}$. Thus the isomorphism $I_{c(1,1,n)}(W)\simeq \NCP(1,1,n)$ induces an isomorphism between $I_{c(1,1,n)}(W)^{c(1,1,n)^q}$ and $\NCP(1,1,n)^{\mu_d}$ as claimed.
\end{proof}

\subsection{The case $W=G(m,1,n)$ for $m\geqslant 2$}\label{sec:G(m,1,n)}
The highest degree of $W$ is $mn$. We endow $W$ with the Coxeter element 
\[c(m,1,n)=\matrix{0& & & & \zeta_m \\ 1 & 0 &  & & \\ & 1 &\ddots & & \\ & & \ddots &0 & \\ & & & 1 & 0},\]
A $\zeta_{mn}$-regular eigenvector for $c(m,1,n)$ which is given by $\zeta_{mn}^{n-1} e_1+\zeta_{mn}^{n-2} e_2+\cdots+e_n$.

We have seen in Section \ref{sec:G(1,1,n)} that $G(m,1,n)$ can be realized as the centralizer in $W':=G(1,1,mn)$ of $c(1,1,mn)^{m}$, and this identification maps $c(m,1,n)$ to $c(1,1,mn)$. We have
\[W_{c(m,1,n)^q}\simeq (W_{c(1,1,mn)^m})_{c(1,1,mn)^{q}}=W'_{c(1,1,mn)^{q\vee m}}.\]
This identification induces the following isomorphisms of posets
\begin{align*}
I_{c(m,1,n)}(W)^{c(m,1,n)^q}&\simeq (I_{c(1,1,mn)}(W'_{c(1,1,mn)^m}))^{c(1,1,mn)^q},\\
I_{c(m,1,n)}(W_{c(m,1,n)^q})\simeq I_{c(1,1,mn)}(W'_{c(1,1,mn)^{q\vee m}}).
\end{align*}
By Section \ref{sec:G(1,1,n)}, we have equalities of posets
\begin{align*}
(I_{c(1,1,mn)}(W'_{c(1,1,mn)^m}))^{c(1,1,mn)^q}&=(I_{c(1,1,mn)}(W')^{c(1,1,mn)^m})^{c(1,1,mn)^q}\\
&=I_{c(1,1,mn)}(W')^{c(1,1,mn)^{m\vee q}}\\
&=I_{c(1,1,mn)}(W'_{c(1,1,mn)^{q\vee m}}),
\end{align*}
which give the desired results.

\subsubsection{The case $W=G(e,e,n)$ for $e\geqslant 2$} The highest degree of $W$ is $e(n-1)$. We start by considering the group $W':=G(e,1,n-1)$. The character $\chi:W'\to \C^*$ giving the product of the nonzero entries allows us to define an embedding
\[\begin{array}{rccc}i:&W'&\longrightarrow & W\\ & M&\longmapsto & \matrix{\chi(M)^{-1}& 0 \\ 0 & M}.\end{array}\]

The element $c(e,e,n)=i(c(e,1,n-1))$ is a Coxeter element for $W$. For instance $\zeta_{e(n-1)}^{n-2}e_2+\cdots+\zeta_{e(n-1)}e_{n-1}+e_n$ is a $\zeta_{e(n-1)}$-eigenvector for $c(e,e,n)$ which is regular. We have
\[c(e,e,n)^q=\matrix{\zeta_e^{-q} & 0 \\0 &c(e,1,n-1)^q}.\]
We can compute directly the centralizer of this element in $W$. 

\begin{prop}
\begin{enumerate}[(a)]
\item If $\frac{e(n-1)}{e\wedge n}$ divides $q$, then $W=W_{c(e,e,n)^q}$. In particular Theorem \ref{theo:intervalles_identiques_p=1} holds.
\item If $\frac{e(n-1)}{e \wedge n}$ does not divide $q$, then $i$ induces an isomorphism between $W'_{c(e,1,n-1)^q}$ and $W_{c(e,e,n)^q}$.
\end{enumerate}
\end{prop}
\begin{proof}
$(a)$ Assume that $\frac{e(n-1)}{e\wedge n}j=q$ for some integer $j$. Since $dq=e(n-1)$, we obtain that $d=\frac{e\wedge n}{j}$ divides $e\wedge n$, which is the gcd of the degrees of $W$ (see \cite[Table 2]{bmr}). The degrees of $W_{c(e,e,n)}$ are then the same as that of $W$.\newline $(b)$ Let $w\in W$. We write $w$ as a block matrix.
\[w=\matrix{X & Y \\ Z & T},\]
with $X\in \Mm_1(\C),Y\in \Mm_{1,n-1}(\C), Z\in \Mm_{n-1,1}(\C)$ and $T\in \Mm_{n-1,n-1}(\C)$. We find
\[w\in W_{c(e,e,n)^q}\Leftrightarrow \begin{cases}c(e,1,n-1)^qZ=\zeta_e^{-q}Z,\\ Yc(e,1,n-1)^q=\zeta^{-q}_e Y,\\ Tc(e,1,n-1)^q=c(e,1,n-1)^q T.\end{cases}\]
Since $w$ is a monomial matrix, $Z$ and $Y$ have at most one nonzero coefficient. Suppose that $Z\neq 0$, then $c(e,1,n-1)^q$ has a diagonal coefficient equal to $\zeta_d^{-q}$. The underlying permutation of $c(e,1,n-1)$ is a $(n-1)$-cycle. Thus $c(e,1,n-1)^q$ has a nonzero diagonal coefficient if and only if $(n-1)$ divides $q$.

Let $j$ be an integer such that $j(n-1)=q$. Since $dq=e(n-1)$, we have $dj=e$. We also have $c(e,1,n-1)^q=(c(e,1,n-1)^{n-1})^j=\zeta_e^{j}\Id$. Thus $c(e,1,n-1)^q Z=\zeta_e^{-q}Z$ if and only if $\zeta_e^j=\zeta_e^{-q}$, that is, $j=-q$ mod $e$. Since $j(n-1)=q$, we obtain that there is some integer $k$ with $jn=kjd$. Thus $kd=n$ and $d$ divides $n$. We have already shown that $d$ divides $e$, we obtain that $d$ divides $e\wedge n$ and that $\frac{e(n-1)}{e\wedge n}$ divides $q$.

The same reasoning shows that $Y\neq 0$ implies that $\frac{e(n-1)}{e\wedge n}$ divides $q$. We obtain $Y=0$, $Z=0$, and $w\in W_{c(e,e,n)^q}$ if and only if $w=i(T)$ with $T\in W'_{c(e,1,n-1)^q}$.
\end{proof}

From now on, we suppose that $\frac{e(n-1)}{e\wedge n}$ does not divide $q$. The isomorphism between $W_{c(e,e,n)^q}$ and $W'_{c(e,1,n-1)^q}$ induced by $i$ induces in turn an isomorphism of posets
\[I_{c(e,1,n-1)}(W'_{c(e,1,n-1)^q})\simeq I_{c(e,e,n)}(W_{c(e,e,n)^q}).\]
On the other hand, we know from Section \ref{sec:G(m,1,n)} that the two posets
\[I_{c(e,1,n-1)}(W)^{c(e,1,n-1)^q}\text{~and~}I_{c(e,1,n-1)}(W'_{c(e,1,n-1)^q})\]
are equal. To show that Theorem \ref{theo:intervalles_identiques_p=1} holds in this case, it only remains to show that the morphism $i$ induces an isomorphism of posets between $I_{c(e,1,n-1)}(W')^{c(e,1,n-1)^q}$ and $I_{c(e,e,n)}(W)^{c(e,e,n)^q}$.
By \cite[Remark after Lemma 1.22]{bescor}, $c(e,e,n)$ induces a free action (of a cyclic group of order $\frac{e(n-1)}{e\wedge n}$) on the set $I_{c(e,e,n)}(W)\setminus i(I_{c(e,1,n-1)}(W'))$. Since $\frac{e(n-1)}{e\wedge n}$ does not divide $q$ by assumption, we get
\[I_{c(e,e,n)}(W)^{c(e,e,n)^q}=i(I_{c(e,1,n-1)}(W'))^{c(e,e,n)^q}=i\left(I_{c(e,1,n-1)}(W')^{c(e,1,n-1)^q}\right),\]
which finishes the proof.

\section{Properties of Springer categories}\label{sec:properties_of_springer_categories}
In essence, Theorem \ref{theo:groupoidtressequivgroupoiddiv} states that Springer categories provide an analogue of the dual braid monoid for centralizers of regular braids in braid groups of well-generated complex reflection groups. This section is dedicated to the study of Springer categories as Garside categories and to the results we deduce for the associated braid groups.

Throughout this section, we fix a finite dimensional complex vector space $V$ with $n:=\dim(V)$, along with an irreducible well-generated complex reflection group $W<\GL(V)$. We denote by $h$ the highest degree of $W$. We also fix a regular number $d$ for $W$ and we set again
\[p:=\frac{d}{d\wedge h},~q:=\frac{h}{d\wedge h}.\]
Moreover, we fix two positive integers $\eta,\mu$ such that $p\eta-q\mu=1$. Lastly, we fix some Coxeter element $c$ in $W$ and we denote by $\Cc$ (resp. by $\Gg$) the Springer category (resp. the Springer groupoid) attached to $(W,c,p,q)$. The Garside map of $\Cc$ will be denoted by $\Delta_p$.

By Theorem \ref{theo:groupoidtressequivgroupoiddiv}, the groupoid $\Gg$ is equivalent to the braid group of the centralizer of some $\zeta_d$-regular element in $W$.

\begin{rem}Recall that the isomorphism $B(W)\simeq G(c)$ identifies the full-twist with $\bbf{c}^h$, where $h$ is the highest degree of $W$. The category of $(1,h)$-periodic elements is then simply the monoid $M(c)$ itself. Thus the content of this section also applies to dual braid monoids of well-generated irreducible complex reflection groups.
\end{rem}

We denote by $I_c^q$ the set of elements of $I_c$ which are invariant under $\phi^q$ (that is, invariant under conjugation by $c^q$).

In this section, we will write the simple elements of $M(c)$ in normal font instead of the bold font used in Section \ref{dbm&dc} for readability purposes. In other words, we completely identify the interval $I_c$ with the set of simples of $M(c)$. Products of elements of $I_c$ should be understood as products in $M(c)$ unless explicitly stated.

\subsection{Elementary properties}\label{sec:3.1}
Recall from Definitions \ref{def:graph_of_periodic_simples} and \ref{def:periodic_elements_category} that $\Cc$ depends on the sets $D_p^q(\bbf{c})$, $D_{2p}^{2q}(\bbf{c})$ and $D_{3p}^{3q}(\bbf{c})$. Thanks to Remark \ref{rem:divided_set_interval_monoid} and Proposition \ref{prop:descriptiondivset}, we can replace these sets with
\begin{align*}
\Ob(\Cc)&:=\Dd_1=\left\{u\in I_c^q~\left|~\begin{cases} p\ell_R(u)=\ell_R(c)=n \\ u(u^{c^{\eta}})\cdots (u^{c^{(p-1)\eta}})=c\end{cases}\right.\right\}\\
\Ss_p^q &:=\Dd_2=\left\{(a,b)\in (I_c^q)^2~|~ ab \in \Ob(\Cc)\text{~and~}p(\ell_R(a)+\ell_R(b))=n\right \}\\
R_p^q&:=\Dd_3=\left\{(x,y,z)\in (I_c^q)^3~|~ xyz\in \Ob(\Cc)\text{~and~}p(\ell_R(x)+\ell_R(y)+\ell_R(z))=n\right \}
\end{align*}

Under these definitions, a simple morphism $(a,b)$ in the graph of simples $\Ss_p^q$ has respective source and target  $ab$ and $ba^{c^\eta}$. An element $(x,y,z)$ of $\Dd_3$ induces the relation $(x,yz)(y,zx^{c^\eta})=(xy,z)$ in $R_p^q$. The following lemma is an obvious consequence of the definition of $\Ob(\Cc)$.

\begin{lem}
Let $u\in \Ob(\Cc)$, the morphism $\Delta_p(u)=(u,1)$ has length $\ell_R(u)=\frac{n}{p}$. This length is independent from from $u$ and every simple morphism in $\Cc$ has length at most $\frac{n}{p}$.
\end{lem}

Let $u\in \Dd_1$ be an object of $\Cc$. The collapse functor $\pi_p:\Cc\to M(c)$ sends $\Delta_p^q(u)$ to some $(p,q)$-periodic element in $M(c)$. We know that $\Delta_p(u)=(u,1)$ and $\Delta_p(\phi_p^n(u))=(u^{c^{\eta n}},1)$. We then have
\begin{align*}
\pi_p(\Delta_p^p(u))&=\pi_p\left((u,1)(u^{c^\eta},1)\cdots (u^{c^{(p-1)\eta}},1)\right)\\
&=uu^{c^\eta}\cdots u^{c^{(p-1)\eta}}=c.
\end{align*}
If $q:=pk+r$ is the Euclidean division of $q$ by $p$, then we get
\[\pi_p(\Delta^q(u))=c^ku^{c^{pk\eta}}\cdots u^{c^{(pk+r-1)\eta}}.\]
We know from Theorem \ref{theo:categoriediviseeetelementsreguliers} that this is a $(p,q)$-regular element in $M(c)$. This gives an explicit formula for roots of the full-twist in $B(W)$, provided that one knows how to compute elements of the sets $\Dd_1$.

\begin{lem}\label{lem:atoms_of_periodic_cat}
The atoms of $\Cc$ are exactly the simples $s=(a,b)\in \Ss_p^q$ such that $a$ admits no proper divisors in $I_c^q$. In particular if $\phi^q$ is trivial, then a simple morphism of $\Cc$ is an atom if and only if it has length 1.
\end{lem}
\begin{proof}The first statement follows directly from Lemma \ref{lem:preceqindivided}. If $\phi^q$ is trivial (that is, if $c^q$ is central in $M(c)$), then Lemma \ref{lem:preceqindivided} gives an isomorphism between $\Ss_p^q(u,-)$ and the interval $\{s\in I_c~|~s\preceq u\}$, where $u$ is the source of $s$. Thus $s$ is an atom if and only if $a$ lies in $R_c$. Since by definition $\ell(s)=\ell_R(a)$, we get the desired result.
\end{proof}

\begin{lem}\label{lem:simples_are_lcm_of_dividing_atoms}
Let $s\in \Cc$ be a simple morphism, and let $A$ be the set of atoms of $\Cc$ which left-divide $s$. The simple morphism $s$ is the right-lcm in $\Cc$ of the elements of $A$.
\end{lem}
\begin{proof}
We start with the case $p=1$. By Corollary \ref{cor:dual_braid_relations_for_centralizer_of_c^q}, it is enough to consider the case of the dual braid monoid $M(c)$. Let $t$ be the right-lcm of the elements of $A$. By definition of $A$ we have $tr=s$ for some simple element $r$. As we mentioned in Section \ref{dbm&dc}, the sets of left- and right-divisors of $s$ coincide and thus $r\preceq s$. If $r$ is nontrivial, then it is left-divisible by some atom $a$. We have $a\preceq s$ and thus $a\preceq t$. Since $t$ is balanced, we obtain that $aa$ is a simple element of $M(c)$, which contradicts Lemma \ref{2.4}. We then have that $r$ is trivial and $s=t$ as we wanted to show.

Assume now that $p>1$. Let $u$ be the source of $s$ in $\Cc$ and write $s=(a,b)\in \Ss_p^q$. By Lemma \ref{lem:preceqindivided}, the collapse functor $\Cc\to M(c)$ induces a poset isomorphism between the divisors of $s$ in $\Cc$ and the divisors of $a$ in $M(c)^{\phi^q}$. By Lemma \ref{lem:atoms_of_periodic_cat}, atoms of $\Cc$ are sent to atoms of $M(c)^{\phi_q}$ by the collapse functor. The element $a\in M(c)^{\phi^q}$ is the right-lcm of the images in $M(c)^{\phi^q}$ of the elements of $A$ after the case $p=1$, and thus $s$ is the right-lcm of the elements of $A$.
\end{proof}

\begin{prop}\label{cor:no_parallel_simples}
If $p>1$, then there is at most one simple morphism in $\Cc$ with given source and target. In particular, for $u\in \Ob(\Cc)$, we have $\Ss_p^q (u,u)=\{1_u\}$ and $\Ss_p^q (u,\phi_p(u))=\{\Delta_p(u)\}$.
\end{prop}

\begin{proof}Let $s:=(a,b)$ be a simple morphism in $\Cc$, and let $u,v$ denote the source and target of $s$, respectively. We have $u=ab$ and $v=ba^{c^\eta}$. We show that $b$ is the left-gcd of $u$ and $v$  in $I_c$.  The result is immediate when $a=1$ and we can assume that $a\neq 1$. We have $ab=ba^b$, so $b$ is an obvious left-divisor of $ab$ and $ba^{c^\eta}$. If $b$ is not the gcd of $ab$ and $ba^{c^\eta}$, then there is a nontrivial common divisor $d$ of $a$ and $a^{c^\eta}$. Since $p>1$, we have $c=aba^{c^\eta}b^{c^\eta}\cdots (ab)^{c^{(p-1)\eta}}$ and $aba^{c^\eta}$ is simple in $M(c)$. Since all elements of $I_c$ are balanced, we obtain that $d^2$ divides $aba^{c^\eta}$, thus $d^2\in I_c$, which contradicts Lemma \ref{2.4}.

Now, let $(a',b')\in \Ss_p^q$ have the same source and target as $(a,b)$. We have $u=ab=a'b'$ and $v=ba^{c^\eta}=b'a'^{c^\eta}$. We have $b\preceq b'$ and $b'\preceq b$, so $b=b'$ and $(a,b)=(a',b')$.
\end{proof}

In the case $p=1$. There is only one object in $\Cc$, which corresponds to $c\in D_1(c)$. The above proposition is then false in this case.

\begin{lem}[\textbf{Lifting words expressing simples}]\label{liftword}~\newline 
Let $s:=(a,b)$ be a simple morphism in $\Cc$, and let $a_1\cdots a_r$ be a word in $I_c^q$ expressing $a$ in $M(c)$. There is a unique path $s_1\cdots s_r$ in $\Ss_p^q$ expressing $s$ in $\Cc$ and such that $\pi_p(s_i)=a_i$ for all $i\in \intv{1,r}$.
\end{lem}
\begin{proof}
Let $u:=ab$ be the source of $s$. We proceed by $\succeq$-induction on $s$. If $s$ is an atom, then $a$ is an atom of $I_c^q$ by Lemma \ref{lem:atoms_of_periodic_cat}. The only word in $I_c^q$ expressing $a$ is then $a$ itself and the result is trivial. Now for the general case, we have $a_1x=a$ with $x=a_2\cdots a_r$. By Lemma \ref{lem:preceqindivided}, $s_1:=(a_1,xb)$ is the only atom in $\Cc$ with source $u$ and such that $\pi_p(s_1)=a_1$. By induction hypothesis, there is a unique path $s_2\cdots s_r$ expressing $(x,ba_1^{c^\eta})$ in $\Cc$ and such that $\pi_p(s_i)=a_i$ for $i\in \intv{2,r}$. The path $s_1s_2\cdots s_r$ is then the unique path expressing $s$ and such that $\pi_p(s_i)=a_i$ for $i\in \intv{1,n}$.
\end{proof}

Let $s:=(a,b)\in \Ss_p^q$ be a simple morphism in $\Cc$. Its source is $ab=b\phi^{\eta}(a^{bc^{-\eta}})$ and its target is $b\phi^{\eta}(a)=a^{c^\eta b^{-1}}b$. By Lemma \ref{lem:preceqindivided} and Lemma \ref{lem:deux_sur_trois}, we deduce the existence of the following simple morphisms in $\Cc$
\[s^\flat:=(a^{bc^{-\eta}},b) \text{~and~} s^\#:=(a^{c^\eta b^{-1}},b).\]

\begin{lem}\label{lem:properties_of_diese_and_flat}
Let $s:=(a,b)$ be a simple morphism in $\Cc$.
\begin{enumerate}[(a)]
\item The target of $s^\flat$ is the source of $s$, and the target of $s$ is the source of $s^\#$.
\item We have $(s^{\flat})^{\#}=(s^{\#})^{\flat}=s$.
\item The paths $s^\flat s$ and $ss^\#$ are both in greedy normal form in $\Cc$.
\item We have $\phi_p(s^\#)=(\phi_p(s))^\#$ and $\phi_p(s^\flat)=(\phi_p(s))^\flat$.
\end{enumerate}
\end{lem}
\begin{proof}$(a)$ The target of $s^\flat$ is $ba^b=ab$ and the source of $s^\#$ is $a^{c^\eta b^{-1}} b=ba^{c^\eta}$. \newline $(b)$ We have $s^{\flat\#}=(a^{bc^{-\eta}},b)^\#=(a^{bc^{-\eta}c^\eta b},b)=(a,b)=s$. The same reasoning proves that $(s^\#)^\flat=s$. \newline $(c)$ Since $s=(a,b)$ is a simple morphism, $ab=u$ is a simple element in $M(c)$. In particular $a$ and $b$ are coprime by Proposition \ref{prop:produit_de_simple_simple_dans_dual}. We get that the path $s^\flat s$ is greedy by Corollary \ref{cor:greedy_normal_form_length_2_in_periodic}. The path $s s^\#=(s^\#)^\flat s^\#$ is also greedy by the same argument. \newline $(d)$ We have
\begin{align*}
(\phi_p(s))^\flat&=(a^{c^\eta},b^{c^\eta})^\flat=(a^{c^\eta b^{c^\eta} c^{-\eta}},b^{c^\eta})=(a^b,b^{c^\eta})=\phi_p(s^\flat)
\end{align*}
and
\[(\phi_p(s))^\#=(\phi_p(s^{\#\flat}))^\#=((\phi_p(s^\#))^\flat)^\#=\phi_p(s^\#).\]
\end{proof}

The transformation $s\mapsto s^\#$ is a bijection of the finite set $\Dd_2$. In particular it has finite order. As the source of $s^\#$ is the target of $s$, there is then a smallest integer $n$ such that $s^{(n\#)}=s^\flat$.

\begin{definition}
Let $s$ be a simple morphism in $\Cc$ with source $u\in \Ob(\Cc)$. The \nit{simple loop} (of the object $u$) associated with $s$ is the morphism
\[\lambda(s):=s s^\# s^{\#\#}\cdots s^{((n-1)\#)}s^{(n\#)}\in  \Cc(u,u),\]
where $n$ is the smallest integer such that $s^{(n\#)}=s^{\flat}$. If $s$ is an atom of $\Cc$, then we say that $\lambda(s)$ is an \nit{atomic loop}.
\end{definition}

\begin{lem}\label{lem:atoloopsinsuppersummitset}
Let $s$ be a simple morphism in $\Cc$. The simple loop $\lambda(s)$ is rigid, in particular it lies in its own super-summit set.
\end{lem}
\begin{proof}
By Lemma \ref{lem:properties_of_diese_and_flat} (c), the greedy normal form of $\lambda(s)$ is given by
\[\lambda(s):=s s^\# s^{\#\#}\cdots s^{\flat\flat}s^{\flat}.\]
Since $s^\flat s$ is greedy, we get that $\lambda(s)$ is rigid. By Lemma \ref{lem:cyclinganddecyclingofrigid},  we have $\cyc(\lambda(s))=\lambda(s^\#)$ and $\dec(\lambda(s))=\lambda(s^\flat)$. In particular we see that cycling and decycling a simple loop does not change its $\inf$ or its $\sup$. By Proposition \ref{prop:sssreachedbycycling}, we deduce that $\lambda(s)$ lies in its own super-summit set.
\end{proof}

\begin{rem}
Let $u$ be an object of $M(c)_p$. We have 
\[\Delta_p(u)^\#=(u,1)^\#=(u^{c^\eta},1)=\phi_p(\Delta_p(u))=\Delta_p(\phi_p(u)).\]
So $\lambda(\Delta_p(u))=\Delta_p^k(u)$ where $k$ is the smallest integer such that $\phi_p^k(u)=u$. By definition of $\Cc$, we have that $k$ divides $q$: there is some $k'$ such that $kk'=q$. The element $(\Delta_p)^q(u)$ is then equal to $\lambda(\Delta_p(u))^{k'}$. By Theorem \ref{theo:categoriediviseeetelementsreguliers}, the collapse functor $\pi_p$ identifies the group $\Gg(u,u)$ with the centralizer in $G(c)$ of $\pi_p(\lambda(\Delta_p(u))^{k'})$.
\end{rem}

\subsection{Presentation by Hurwitz relations} We defined Springer categories via a presentation by generators and relations. For computational application, we would like to have a presentation of Springer categories involving fewer generators and relations than the defining presentation. Such a result is known for dual braid monoids, which are presented by their atoms, endowed with all Hurwitz relations (see Remark \ref{rem:egalitesdesimpledelongueur2}).

In this section, we denote by $\Aa$ the set of atoms of the Springer category $\Cc$.
\begin{lem}\label{lem:comm_square_of_atoms_in_periodic}
Consider a square of atoms in $\Cc$
\[\xymatrix{u \ar[r]^-s\ar[d]_-\sigma & v \ar[d]^-t\\ v'\ar[r]_-\tau &w}\]
The square is commutative if and only if $\pi_p(st)=\pi_p(\sigma\tau)$ in $M(c)^{\phi_q}$.\end{lem}
\begin{proof}
First, assume that $p=1$. In this case, we have $\Cc=M(c)^{\phi_q}$ and the collapse functor $\pi_1$ is the natural morphism $M(c)^{\phi_q}\to M(c)$ (see Section \ref{sec:1.3.1}). Our statement is then a rephrasing of the fact that $\pi_1$ is an embedding $M(c)^{\phi_q}\to M(c)$. We assume that $p>1$ from now on.

The direct implication comes from the functoriality of $\pi_p$. Conversely, let $a,b,\alpha,\beta$ be the respective images of $s,t,\sigma,\tau$ under $\pi_p$, and let $x:=ab$. By Lemma \ref{lem:atoms_of_periodic_cat}, $a,b,\alpha,\beta$ are atoms of the Garside monoid $M(c)^{\phi_q}$ and we have $x=ab=\alpha\beta$. 

If $a=\alpha$, then $b=\beta$ by cancellativity of $M(c)^{\phi_q}$. By Lemma \ref{lem:preceqindivided}, we have $s=\sigma$, and $t,\tau$ share the same source and target. Since we assumed that $p>1$, we obtain $t=\tau$ by Proposition \ref{cor:no_parallel_simples} and the square is commutative. 

If $a\neq \alpha$, then $x$ is a common multiple of the two distinct atoms $a,\alpha$. Let $c,\gamma$ be such that $ac=\alpha\gamma=a\vee \alpha$. By definition of lcm, we have $c\preceq b$ and $\gamma \preceq \beta$. Moreover, since $a\neq \alpha$, the elements $c,\gamma$ are both nontrivial. Since $b,\beta$ are atoms, we deduce that $c=b$ and $\gamma=\beta$. In other words $x$ is the right-lcm of $a$ and $\alpha$ in $M(c)^{\phi_q}$. By Lemma \ref{lem:preceqindivided}, we have $a,\alpha\preceq u$, thus $x \preceq u$. Still by Lemma \ref{lem:preceqindivided}, there is a unique simple morphism $f\in \Ss_p^q(u,-)$ such that $\pi_p(f)=x$. Applying Lemma \ref{liftword} to the words $ab,\alpha\beta$, which both express $x$, yields that $f=st=\sigma\tau$.
\end{proof}

\begin{definition}
We call \nit{Hurwitz relations} on $\Cc$ the relations of the form $st=\sigma\tau$, where $s,t,\tau,\sigma\in \Aa$ are such that $s\neq \sigma$ and $\pi_p(st)=\pi_p(\sigma\tau)$. We denote by $\Hh$ the set of Hurwitz relations in $\Cc$.
\end{definition}

By Lemma \ref{lem:comm_square_of_atoms_in_periodic}, the Hurwitz relations are the relations induced by the (nontrivial) commutative squares of atoms in $\Cc$. In the case where $p=1$ and $q=h$, the set $\Hh$ contains the set of Hurwitz relations in the sense of \cite[Definition 8.7]{beskpi1}. In this case, it is known by \cite[Lemma 8.8]{beskpi1} that the atoms endowed with the Hurwitz relation provide a presentation of the dual braid monoid $M(c)$. We now generalize this result to arbitrary Springer categories:

\begin{theo}\label{prop:hurwitzpresentationforc31} The Springer category is presented by its atoms endowed with the Hurwitz relations, that is, we have $\Cc=\langle \Aa ~|~ \Hh\rangle^+$.
\end{theo}
\begin{proof} 
If $p=1$, the result comes from Corollary \ref{cor:dual_braid_relations_for_centralizer_of_c^q}, since dual braid monoids are presented by their atoms, endowed with Hurwitz relations \cite[Lemma 8.8]{beskpi1}. Suppose now that $p>1$. We know that $\Cc$ is generated by its atoms, and that the defining relations of $\Cc$ imply the Hurwitz relations. It remains to show that the defining relations of $\Cc$ are implied by the Hurwitz relations.

Let $st=u$ be a defining relation of $\Cc$. We consider three paths of atoms in $\Cc$
\[s_1\cdots s_r,~t_1\cdots t_k,~u_1\cdots u_{m},\]
expressing $s,t$ and $u$, respectively. We set $a_i=\pi_p(s_i)$ for $i\in \intv{1,r}$, $b_i:=\pi_p(t_i)$ for $i\in \intv{1,k}$ and $\alpha_i:=\pi_p(u_i)$ for $i\in \intv{1,m}$. In $M(c)^{\phi^q}$, the two words
\[a_1\ldots a_rb_1\cdots b_k\text{~and~}\alpha_1\cdots \alpha_{m}\]
express the same element $\pi_p(u)$. Since $M(c)^{\phi^q}$ is presented by the Hurwitz relations (this is the case $p=1$), there is a sequence of words $\mu_1,\ldots,\mu_n$ on the atoms of $M(c)^{\phi^q}$ such that
\begin{itemize}
\item $\mu_1=a_1\cdots a_rb_1\cdots b_r$
\item $\mu_n=\alpha_1\cdots \alpha_{m}$
\item For $i\in \intv{1,n-1}$, $\mu_i$ and $\mu_{i+1}$ are related by a Hurwitz relation in $M(c)^{\phi^q}$.
\end{itemize}
In particular, each of the $\mu_i$ expresses the element $\pi_p(u)$ in $M(c)^{\phi^q}$. By Lemma \ref{liftword}, each $\mu_i$ admits a unique lift $p_i$ in $\Cc$, and each path $p_i$ expresses $u$ in $\Cc$. By Lemma \ref{lem:comm_square_of_atoms_in_periodic}, the paths $p_i$ and $p_{i+1}$ are related by a Hurwitz relation in $\Cc$ for $i\in \intv{1,n-1}$. In particular the equality $st=u$ is implied by the Hurwitz relations in $\Cc$.
\end{proof}

This new presentation will be useful for computational purposes in Section \ref{sec:presentation_of_b31}.

\subsection{Braided reflections and atomic loops}\label{sec:braid_reflections_and_atomic_loops}
In this section, we give a description of braided reflections in $B(W_g)$ in terms of atomic loops in the category $\Cc$. Namely we prove that the isomorphism of Theorem \ref{theo:groupoidtressequivgroupoiddiv} identifies braided reflections with conjugates of atomic loops. We freely use the notation of Section \ref{sec:2.3}.

Since the homeomorphism between $X_g/W_g$ and $(X/W)^{\mu_d}$ does not depend on a choice, we can identify the two spaces. The choice of a basepoint $x\in \Uu^{\mu_d}$ induces an isomorphism between $B(W_g)\simeq \pi_1(X_g/W_g,x)$ and $\Gg(u,u)$, where $u$ is the connected component of $x$ in $\Uu^{\mu_d}$. By Remark \ref{rem:basepoints_and_braid_reflections}, a change of basepoint preserves the set of braided reflections.

Our main tool is a particular set of neighborhoods of $a$ in $(V/W)^{\mu_d}$, whose image under $\LLL$ consists of a set of small disks around the points of $\LLL(a)$. Using these neighborhoods, we can clarify the definition of braided reflections in the context of the topological tools of Section \ref{sec:2.3} and \ref{sec:2.4}. 

First, we recall that the description of $X/W$ using the map $\LLL$ can be refined into a description of $(X/W)^{\mu_d}$. For $x\in (V/W)^{\mu_d}$, we have $\LLL(\zeta_d.x)=\zeta_d^h\LLL(x)=\LLL(x)$ by Lemma \ref{lem:lll_homog}. In other words,  $\LLL(x)$ lies in the set $(E_n)^{\mu_p}$ of $\mu_p$-invariant multisets. In particular, the number of points of $\LLL(x)$ lying in a sector $P_i$ does not depend on $i$. Now, if $x\in (X/W)^{\mu_d}$ is such that $\LLL(x)$ contains $k$ points in the sector $P_1$, then $\clbl(x)$ lies in $D_{kp}^{kq}(c)$ by \cite[Lemma 11.11]{beskpi1}.

\begin{prop}\label{prop:restriction_image_standard}\cite[Definition 11.12 and Proposition 11.13]{beskpi1} The bijection $(\LLL,\clbl)$ of Theorem \ref{theo:imagestandard} restricts to a bijection between $(X/W)^{\mu_d}$ and the set $(E_n^\circ)^{\mu_p} \boxdot D_{\bullet p}^{\bullet q}(c)$ of compatible pairs where the first term is $\mu_p$-invariant, and the second term belongs to some $D_{kp}^{kq}(c)$.
\end{prop}

Following this proposition, we can now see the cyclic label of some $x\in (X/W)^{\mu_d}$ as an element of $D_{\bullet p}^{\bullet q}(c)$. Following the argument of Lemma \ref{lem:descriptionpointfixesphiq} the cyclic label of $x\in (X/W)^{\mu_d}$ only depends on the terms corresponding to points in the sector $P_1$. We will thus abusively identify $\clbl(x)$ with the terms corresponding to points in the sector $P_1$ from now on.


Recall that $E_n$ is defined as the quotient space $\C^n/\mathfrak{S}_n$. The projection map $\C^n\to E_n$ is open and we can define a neighborhood of some multiset $\kappa=\muls{z_1,\ldots,z_n}$ in $E_n$ as the image in $E_n$ of some neighborhood of $(z_1,\ldots,z_n)$ in $\C^n$. For $\kappa=\muls{z_1,\ldots,z_n}\in E_n$ and $r>0$ we denote by $\D(\kappa,r)$ the image in $E_n$ of the set $\prod_{i=1}^n \D(z_i,r)\subset \C^n$. A disk $\D(z_i,r)$ is called either an central or an outer disk of $\D(\kappa,r)$, depending on wether $z_i=0$ or not.

\begin{definition}[\textbf{Confining neighborhood}] Let $\kappa=\muls{z_1,\ldots,z_n}\in E_n$ be such that $\kappa\cap D\subset \{0\}$ and let $r>0$. We say that $\D(\kappa,r)$ is confining if it satisfies the following conditions:
\begin{enumerate}[(1)]
\item If $z_i\neq z_j$, then $\D(z_i,r)\cap \D(z_j,r)=\varnothing$. In other words two disks composing $\D(\kappa,r)$ are either equal or they intersect trivially.
\item If $z_i\neq 0$, then $\D(z_i,r)\cap D=\varnothing$. In other words outer disks intersect $D$ trivially.
\end{enumerate}
More generally, we call \nit{confining neighborhood} of $\kappa$ any neighborhood contained in some confining $\D(\kappa,r)$. If $a\in V/W$ is such that $\LLL(a)\cap D\subset \{0\}$, then we say that a path connected neighborhood of $a$ in $V/W$ is \nit{confining} if its image under $\LLL$ is contained in a confining neighborhood of $\LLL(a)$.
\end{definition}

The philosophy is that a confining neighborhood $\D(\kappa,r)$ ``isolates'' the points in the central disk as in the example below \begin{center}\begin{tikzpicture}[scale=1]

\draw [dashed] (0,0)--(0,2);
\draw [dashed] (0,0)--(2*0.866,-1);
\draw [dashed] (0,0)--(-2*0.866,-1);

\node[draw,circle,inner sep=0.5pt,fill] (P1) at (0,0) {};
\node[draw,circle,inner sep=0.5pt,fill] (P2) at (0.1987,0.9801) {};
\node[draw,circle,inner sep=0.5pt,fill] (P3) at (0.6*0.9801+1.5*0.1987,-0.6*0.1987+1.5*0.9801) {};
\node[draw,circle,inner sep=0.5pt,fill] (P4) at (0.6*0.9801+0.5*0.1987,-0.6*0.1987+0.5*0.9801) {};
\node[draw,circle,inner sep=0.5pt,fill] (P5) at (1.2*0.9801-0.8*0.1987,-1.2*0.1987-0.8*0.9801) {};
\node[draw,circle,inner sep=0.5pt,fill] (P6) at (-0.866*0.9801-0.5*0.1987,0.866*0.1987-0.5*0.9801) {};

\draw [dashed] (P1) circle (0.1);
\draw [dashed] (P2) circle (0.1);
\draw [dashed] (P3) circle (0.1);
\draw [dashed] (P4) circle (0.1);
\draw [dashed] (P5) circle (0.1);
\draw [dashed] (P6) circle (0.1);

\end{tikzpicture}\end{center}

Let $a\in V/W$ be such that $\LLL(a)\cap D\subset \{0\}$. It is easy to see that the family of confining neighborhoods of $a$ is stable under intersection and that it provides a basis of neighborhoods of $a$ in $V/W$. Moreover, by \cite[Theorem 5.3]{beskpi1}, the map $\LLL$ is a finite map between two analytic varieties $V/W$ and $E_n$, both isomorphic to $\C^n$. In particular, it is an open map. Since $\LLL$ is surjective, any confining neighborhood of $a$ contains the inverse image under $\LLL$ of some confining neighborhood $\D(\LLL(a),r)$.

Let now $a \in (V/W)^{\mu_d}$ be such that $\LLL(a)\cap D\subset \{0\}$, and let $U$ be a confining neighborhood of $a$ in $V/W$. The set $U\cap (V/W)^{\mu_d}$ is a neighborhood of $a$ in $(V/W)^{\mu_d}$. In particular we have that $U\cap \Uu^{\mu_d}$ is nonempty since $\Uu^{\mu_d}$ is dense in $(V/W)^{\mu_d}$.

\begin{lem}\label{lem:confining_neighb_atomic_loops}
Let $a\in (V/W)^{\mu_d}$ be such that $\LLL(a)\cap D\subset \{0\}$, and let $U$ be a confining neighborhood of $a$ in $V/W$. Let also $x\in U\cap \Uu^{\mu_d}$ and let $u$ be the connected component of $x$ in $\Uu^{\mu_d}$. There is some $s=s(x)\in \Ss_p^q(u,-)$ such that the image of $\pi_1(U\cap(X/W)^{\mu_d},x)$ in $\Gg(u,u)$ contains all the simple loops $\lambda(t)$ with $t\preceq s$ in $\Cc$.
\end{lem}

\begin{proof}
We denote by $\beta=\beta(x)$ the product (in clockwise order) of the terms of $\clbl(x)$ corresponding to the outer points in the sector $P_1$. Let $\gamma$ be a path in $U\cap (X/W)^{\mu_d}$ starting from $x$. For all $t\in [0,1]$ such that $\gamma(t)\in \Uu^{\mu_d}$, we have $\beta(\gamma(t))=\beta(x)$ by Proposition \ref{prop:cyclic_hurwitz_moves} and the following discussion.

We consider the path in $E_n$ which consists in sliding the central points in each sector together and then counterclockwise next to the associated half-line:

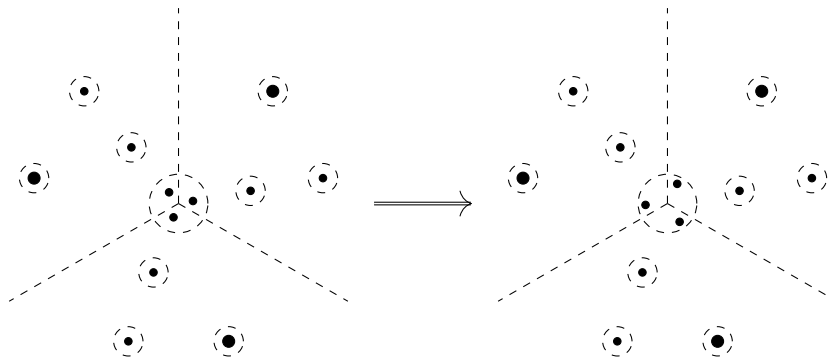
\begin{center}\begin{tikzpicture}[scale=1.3]

\draw [dashed] (-2,0)--(-2,2);
\draw [dashed] (-2,0)--(2*0.866-2,-1);
\draw [dashed] (-2,0)--(-2*0.866-2,-1);

\node[draw,circle,inner sep=1pt,fill] (P2) at (1.477-2,0.260) {};
\node[draw,circle,inner sep=1pt,fill] (P4) at (-0.513-2,-1.410) {};
\node[draw,circle,inner sep=1pt,fill] (P6) at (-0.964-2,1.149) {};

\node[draw,circle,inner sep=1.6pt,fill] (P1) at (0.964-2,1.149) {};
\node[draw,circle,inner sep=1.6pt,fill] (P3) at (0.513-2,-1.410) {};
\node[draw,circle,inner sep=1.6pt,fill] (P5) at (-1.477-2,0.260) {};

\node[draw,circle,inner sep=1pt,fill] (P22) at (0.737-2,0.13) {};
\node[draw,circle,inner sep=1pt,fill] (P44) at (-0.256-2,-0.705) {};
\node[draw,circle,inner sep=1pt,fill] (P66) at (-0.482-2,0.575) {};

\node[draw,circle,inner sep=1pt,fill] (P21) at (0.1477-2,0.0260) {};
\node[draw,circle,inner sep=1pt,fill] (P41) at (-0.0513-2,-0.1410) {};
\node[draw,circle,inner sep=1pt,fill] (P61) at (-0.0964-2,0.1149) {};

\draw [dashed] (0.964-2,1.149) circle (0.15);
\draw [dashed] (P2) circle (0.15);
\draw [dashed] (0.513-2,-1.410) circle (0.15);
\draw [dashed] (P4) circle (0.15);
\draw [dashed] (-1.477-2,0.260) circle (0.15);
\draw [dashed] (P6) circle (0.15);

\draw [dashed] (P22) circle (0.15);
\draw [dashed] (P44) circle (0.15);
\draw [dashed] (P66) circle (0.15);

\draw [dashed] (-2,0) circle (0.3);

\draw [dashed] (3,0)--(3,2);
\draw [dashed] (3,0)--(3+2*0.866,-1);
\draw [dashed] (3,0)--(3+-2*0.866,-1);

\node[draw,circle,inner sep=1pt,fill] (P2) at (3+1.477,0.260) {};
\node[draw,circle,inner sep=1pt,fill] (P4) at (3-0.513,-1.410) {};
\node[draw,circle,inner sep=1pt,fill] (P6) at (3-0.964,1.149) {};

\node[draw,circle,inner sep=1.6pt,fill] (P1) at (3+0.964,1.149) {};
\node[draw,circle,inner sep=1.6pt,fill] (P3) at (3+0.513,-1.410) {};
\node[draw,circle,inner sep=1.6pt,fill] (P5) at (3-1.477,0.260) {};

\node[draw,circle,inner sep=1pt,fill] (P22) at (3+0.737,0.13) {};
\node[draw,circle,inner sep=1pt,fill] (P44) at (3-0.256,-0.705) {};
\node[draw,circle,inner sep=1pt,fill] (P66) at (3-0.482,0.575) {};

\node[draw,circle,inner sep=1pt,fill] (P11) at (3+0.1,0.2014) {};
\node[draw,circle,inner sep=1pt,fill] (P31) at (3-0.05-0.866*0.2014,0.0866-0.1007) {};
\node[draw,circle,inner sep=1pt,fill] (P51) at (3-0.05+0.866*0.2014,-0.0866-0.1007) {};

\draw [dashed] (3+0.964,1.149) circle (0.15);
\draw [dashed] (P2) circle (0.15);
\draw [dashed] (3+0.513,-1.410) circle (0.15);
\draw [dashed] (P4) circle (0.15);
\draw [dashed] (3-1.477,0.260) circle (0.15);
\draw [dashed] (P6) circle (0.15);

\draw [dashed] (P22) circle (0.15);
\draw [dashed] (P44) circle (0.15);
\draw [dashed] (P66) circle (0.15);

\draw [dashed] (3,0) circle (0.3);

\draw [double,->] (0,0) to (1,0);
\end{tikzpicture}\end{center}

Up to homotopy, we can assume that this path is a path in $(E_n)^{\mu_p}\cap \D(\LLL(a),r)$ with $\D(\LLL(a),r)\subset \LLL(U)$. The unique lift $\gamma_0$ of this path in $V/W$ starting from $x$ is a path in $U\cap (X/W)^{\mu_d}$. By construction, the path $\gamma_0$ is homotopically trivial and only affects the cyclic label. Let $\alpha=\alpha(x)$ be the first term of the cyclic label of $\gamma_0(1)$. The cyclic label of $\gamma_0(1)$ is given by $(\alpha,s_1,\ldots,s_k)$ with $s_1\cdots s_k=\beta$. We associate to $x$ the well defined element $s=s(x):=(\alpha,\beta)\in \Ss_p^q$ ($\alpha$ is also the unique element such that the equality $u=\alpha\beta$ holds).

Let $t=(d,e)$ be a left-divisor of $s$ in $\Cc$. There is some $m\in \Ss^{\phi^q}$ such that $dm=\alpha$ and $m\beta=e$. We can desingularize $\gamma_0(1)$ so that we obtain a point $x'\in U\cap \Uu^{\mu_d}$ with cyclic label $(d,m,s_1,\ldots,s_k)$. Consider now the path in $(E_n^\circ)^{\mu_p}$ starting from $\LLL(x')$ and consisting in rotating the central points to the left, so that the first point of the cyclic support goes into the sector $P_p$. The unique lift of this path in $U\cap (X/W)^{\mu_d}$ starting from $x'$ represents the simple morphism $t$ in $\Gg$.

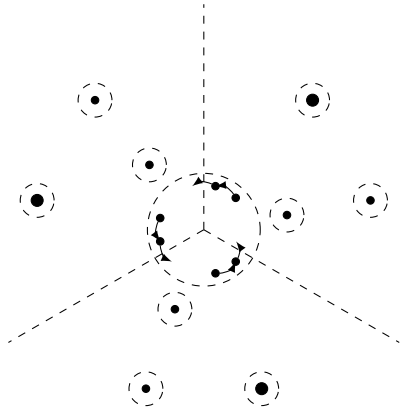
\begin{center}
\begin{tikzpicture}[scale=1.5]

\draw [dashed] (0,0)--(0,2);
\draw [dashed] (0,0)--(2*0.866,-1);
\draw [dashed] (0,0)--(-2*0.866,-1);

\node[draw,circle,inner sep=1pt,fill] (P2) at (1.477,0.260) {};
\node[draw,circle,inner sep=1pt,fill] (P4) at (-0.513,-1.410) {};
\node[draw,circle,inner sep=1pt,fill] (P6) at (-0.964,1.149) {};

\node[draw,circle,inner sep=1.6pt,fill] (P1) at (0.964,1.149) {};
\node[draw,circle,inner sep=1.6pt,fill] (P3) at (0.513,-1.410) {};
\node[draw,circle,inner sep=1.6pt,fill] (P5) at (-1.477,0.260) {};

\node[draw,circle,inner sep=1pt,fill] (P22) at (0.737,0.13) {};
\node[draw,circle,inner sep=1pt,fill] (P44) at (-0.256,-0.705) {};
\node[draw,circle,inner sep=1pt,fill] (P66) at (-0.482,0.575) {};


\node[draw,circle,inner sep=1pt,fill] (P001) at (0.2828,0.2828) {};
\draw[->,>=latex] (P001) to[bend right] (0.1035,0.3864);

\node[draw,circle,inner sep=1pt,fill] (P002) at (0.1035,0.3864) {};
\draw[->,>=latex] (P002) to[bend right] (-0.1035,0.3864);

\node[draw,circle,inner sep=1pt,fill] (P004) at (-0.3864,0.1035) {};
\draw[->,>=latex] (P004) to[bend right] (-0.3864,-0.1035);

\node[draw,circle,inner sep=1pt,fill] (P003) at (-0.3864,-0.1035) {};
\draw[->,>=latex] (P003) to[bend right] (-0.2828,-0.2828);

\node[draw,circle,inner sep=1pt,fill] (P005) at (0.1035,-0.3864) {};
\draw[->,>=latex] (P005) to[bend right] (0.2828,-0.2828);

\node[draw,circle,inner sep=1pt,fill] (P006) at (0.2828,-0.2828) {};
\draw[->,>=latex] (P006) to[bend right] (0.2864,-0.1035);

\draw [dashed] (0.964,1.149) circle (0.15);
\draw [dashed] (P2) circle (0.15);
\draw [dashed] (0.513,-1.410) circle (0.15);
\draw [dashed] (P4) circle (0.15);
\draw [dashed] (-1.477,0.260) circle (0.15);
\draw [dashed] (P6) circle (0.15);

\draw [dashed] (P22) circle (0.15);
\draw [dashed] (P44) circle (0.15);
\draw [dashed] (P66) circle (0.15);

\draw [dashed] (0,0) circle (0.5);


\end{tikzpicture}\end{center}
Let $x_1$ denote the endpoint of this path. The cyclic label of $x_1$ is $(m,s_1,\ldots,s_k,d^{c^{\eta}})$. By applying a new standardization motion $\gamma_1$, we get a point $\gamma_1(1)$ with cyclic label 
\[(d^{c^{\eta}ms_1\cdots s_k},m,s_1,\ldots,s_k)=(d^{c^{\eta}e},m,s_1,\ldots,s_k).\]
We see that the circular tunnel associated with $\gamma(1)$ represents $t^\#$ in $\Cc$. By an immediate induction, we obtain that $\lambda(t)$ lies in the image of $\pi_1(U\cap (X/W)^{\mu_d},x)$ in $\Gg(u,u)$.
\end{proof}

\begin{theo}\label{prop:braidreflectionconjugatetoatomicloop}
Let $x\in \Uu^{\mu_d}$ and let $u$ be the connected component of $x$ in $\mu^{\mu_d}$. The isomorphism $\Gg(u,u)\simeq B(W_g)=\pi_1(X_g/W_g,x)$ maps atomic loops of $u$ to braided reflections. Conversely, any braided reflection in $B(W_g)\simeq \Gg(u,u)$ is conjugate in $\Gg$ to some atomic loop.
\end{theo}
\begin{proof}
First, if $p=1$, then $\Gg=M(c)^{\phi^q}$ is isomorphic to the dual braid monoid associated with $W_g$ by Corollary \ref{cor:dual_braid_relations_for_centralizer_of_c^q}. The result is already known in this case: the braided reflections in $B(W_g)$ are exactly the conjugates of atoms of the dual braid monoid.

From now on we assume that $p>1$. Let $s=(\alpha,\beta)$ be an atom of $\Cc$. We consider the element $x_s$ of $s$ used to define the circular tunnel $b_s$ of Theorem \ref{theo:groupoidtressequivgroupoiddiv}. Let $\gamma$ be the path in $(X/W)^{\mu_d}$ starting from $x_s$ and whose image in $E_n$ consists in sliding the first point of $\LLL(x_s)$ in each sector towards the center
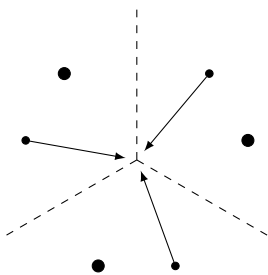
\begin{center}
\begin{tikzpicture}[scale=1]

\draw [dashed] (0,0)--(0,2);
\draw [dashed] (0,0)--(2*0.866,-1);
\draw [dashed] (0,0)--(-2*0.866,-1);
\node[draw,circle,inner sep=1pt,fill] (P1) at (0.964,1.149) {};
\node[draw,circle,inner sep=1.6pt,fill] (P2) at (1.477,0.260) {};
\node[draw,circle,inner sep=1pt,fill] (P3) at (0.513,-1.410) {};
\node[draw,circle,inner sep=1.6pt,fill] (P4) at (-0.513,-1.410) {};
\node[draw,circle,inner sep=1pt,fill] (P5) at (-1.477,0.260) {};
\node[draw,circle,inner sep=1.6pt,fill] (P6) at (-0.964,1.149) {};

\draw[->,>=latex] (P1) to (0.0964,0.1149);
\draw[->,>=latex] (P3) to (0.0513,-0.1410);
\draw[->,>=latex] (P5) to (-0.1477,0.0260);

\end{tikzpicture}\end{center}
The endpoint of $\gamma$ lies on the discriminant hypersurface. Let $U$ be a confining neighborhood of $\gamma(1)$ and let $t_0>0$ be such that $\gamma(t_0)\in U$. The cyclic label of $\gamma(t_0)$ is $(\alpha,\beta)$. By Lemma \ref{lem:confining_neighb_atomic_loops}, the image of $\pi_1(U\cap (X/W)^{\mu_d},\gamma(t_0))$ in $\Gg(u,u)$ contains $\lambda(s)$. 

Next, let $h\in \pi_1(U\cap (X/W)^{\mu_d},\gamma(t_0))$ be represented by a loop $\theta$. We can choose $\theta$ so that, at any given $t\in [0,1]$, at most one point of $\LLL(\theta(t))$ lies in the half-line $i\R^+$ (as points which do not satisfy this form a a subspace of real codimension $2$). This expresses $\theta$ as a concatenation of paths in $U\cap (X/W)^{\mu_d}$ homotopic to circular tunnels (and their inverses).

Let $s'=(d,e)$ be a simple morphism represented by a path inside $U\cap (X/W)^{\mu_d}$ and with source $u$. By Proposition \ref{prop:cyclic_hurwitz_moves} and the following discussion, both the source and target of $s'$ are elements of $M(c)$ divisible by $\beta$. Since $p>1$, the proof of Proposition \ref{cor:no_parallel_simples} gives that $e$ is the gcd in $M(c)$ of the source and target of $s'$. We then have that $\beta$ divides $e$ and $s'$ left-divides $s$. The same reasoning (by induction) gives that the only morphisms represented by a concatenation of paths in $U\cap (X/W)^{\mu_d}$ (starting from $\gamma(t_0)$) homotopic to circular tunnels are of the form $s s^{\#}s^{\#\#}\cdots$ or $(ss^{\#}\cdots)^{-1}$. We obtain that $\pi_1(U\cap (X/W)^{\mu_d},\gamma(t_0))$ is cyclic and generated by a path representing $\lambda(s)$ in $\Gg$, which proves that $\lambda(s)$ is a braided reflection in $\Gg(u,u)$.

Conversely, let $\gamma$ be a path from $x$ to the smooth part of some irreducible divisor of the discriminant in $X_g/W_g\simeq (X/W)^{\mu_d}$. Let also $\rho_\gamma$ denote the braided reflection associated with $\gamma$. Up to conjugacy in $\Gg(u,u)$, we can assume that the endpoint $a$ of $\gamma$ is such that no point of $\LLL(a)$ lies on $D$ and that, for all $t\in[0,1[$, there is some $r\in [t,1[$ such that $\LLL(\gamma(r))\in \Uu^{\mu_d}$.

Let $U$ be a small enough confining neighborhood of $a$ so that the fundamental group of $U\cap (X/W)^{\mu_d}$ is infinite cyclic. Let also $t_0\in [0,1]$ be such that $y:=\gamma(t_0)$ lies in $U\cap \Uu^{\mu_d}$. We denote by $v$ the connected component of $y$ in $\Uu^{\mu_d}$. 

Consider the element $s=s(y)$ introduced in Lemma \ref{lem:confining_neighb_atomic_loops}. We claim that $\lambda(s)$ is an atomic loop. Let $\sigma,\sigma'$ be two atoms dividing $s$ in $\Cc$. Lemma \ref{lem:confining_neighb_atomic_loops} implies that the image of $\pi_1(U\cap (X/W)^{\mu_d},y)$ in $\Gg(v,v)$ contains both $\lambda(\sigma)$ and $\lambda(\sigma')$. Since $\pi_1(U\cap (X/W)^{\mu_d},y)$ is cyclic, there are two integers $k$ and $k'$ such that $\lambda(\sigma)^k=\lambda(\sigma)^{k'}$. Since both $\lambda(\sigma)$ and $\lambda(\sigma')$ are rigid, this implies that $\sigma=\sigma'$. We have shown that there is exactly one atom of $\Cc$ which left-divides $s$. By Lemma \ref{lem:simples_are_lcm_of_dividing_atoms}, this implies that $s$ is an atom and that $\lambda(s)$ is an atomic loop.

Again, as every element of $\pi_1(U\cap (X/W)^{\mu_d},y)$ is represented by a path homotopic to a concatenation of circular tunnels, we obtain that the image of $\pi_1(U\cap(X/W)^{\mu_d},y)$ in $\Gg(v,v)$ is generated by $\lambda(s)$. The path $\gamma_{t_0}:t\mapsto \gamma(t_0t)$ then gives a conjugating element in $\Gg$ between the braided reflection $\rho_{\gamma}$ in $\Gg(u,u)$ and $\lambda(s)\in \Gg(v,v)$.
\end{proof}

\subsection{Conjugacy of atomic loops and centers of finite index subgroups} In this section, we study the conjugacy of atomic loops and their powers in $\Cc$. We show in particular that the centralizer of an atomic loop coincides with those of its nontrivial powers. We deduce a new Garside theoretic proof of \cite[Theorem 1.4]{dmm} in the case of the centralizer of a regular element in a well-generated irreducible complex reflection group. We begin by describing the center of Springer categories and Springer groupoids.

Recall that the \nit{center} $Z(\Dd)$ of a small category $\Dd$ is defined as the monoid of natural endomorphisms of its identity functor. The data of an element $z$ of $Z(\Dd)$ is equivalent to the data, for every object $u$ of $\Dd$, of a morphism $z_u\in \Dd(u,u)$ such that
\[\forall f\in \Dd(u,v), z_uf=fz_v.\]
If $\Dd=M$ is a monoid, we recover the classical definition of the center. Note that if $\Dd$ is a groupoid, then $Z(\Dd)$ is a group.

\begin{lem}
Let $\Dd$ be a connected groupoid, and let $u\in \Ob(\Dd)$. The map
\[\begin{array}{rrcl}r:&Z(\Dd)&\longrightarrow & Z(\Dd(u,u))\\ &z&\longmapsto & z_u\end{array}\]
is an isomorphism of groups.
\end{lem}
\begin{proof} The map $r$ is clearly a morphism of groups. We choose, for every $v\in \Ob(\Dd)$, a morphism $m_v:u\to v$. Let $z_0\in Z(\Dd(u,u))$ be a central element. We obtain an element $z$ of $Z(\Dd)$ by defining $z_v:=m_v^{-1}z_0m_v$ for every $v\in \Ob(\Dd)$. We have in particular that $z_u=m_u^{-1}z_0m_u=z_0$ since $z_0$ is central. Let $f:v\to w$ be a morphism in $\Dd$, we have
\[z_vf=m_v^{-1}z_0m_vf=m_v^{-1}z_0(m_v fm_w^{-1})m_w=m_v^{-1}(m_vfm_w^{-1})z_0m_w=fz_w,\]
and $z$ is indeed in $Z(\Dd)$. The map $z_0\mapsto z$ is the inverse of $r$.
\end{proof}

Now, it is known from \cite[Theorem 2.24]{bmr} and \cite[Theorem 12.3 and Corollary 12.7]{beskpi1} that the center of an irreducible complex braid group is cyclic. The center of the Springer groupoid $\Gg$ is also cyclic since $\Gg$ is connected. Under a combinatorial assumption on the integer $d$, we get that the center of $\Cc$ (and $\Gg$) is actually generated by some power of the Garside map $\Delta$. 

\begin{prop}\label{prop:center_of_categories}If $d$ is the gcd of the degrees of $W$ which are divisible by $d$, then $Z(\Cc)$ (resp. $Z(\Gg)$) is the monoid (resp. group) generated by $\Delta^q_p$.
\end{prop}
\begin{proof}
First, as $\Delta_p$ is a natural transformation from $1_{\Cc}$ to $\phi_p$, and as $\phi_p^q=1_{\Cc}$, we get that $\Delta_p^q\in Z(\Cc)$. Since $\Gg$ is defined as the enveloping groupoid of $\Cc$, we also have $\Delta_p ^q\in Z(\Gg)$.

We claim that $Z(\Gg)$ is generated as a monoid by $Z(\Cc)$ and $\Delta_p^{-q}$. Let $z\in Z(\Gg)$, and let $u\in \Ob(\Gg)$. We set $n_u:=\inf(z_u)$. Since $\Gg$ admits a finite number of objects ($D_p^q(c)$ is finite by definition), there is some $k\in \Z_{\geqslant 0}$ such that $kq+n_u>0$ for all objects $u$. We then have that $\Delta_p^{kq}z\in Z(\Cc)$ and $z=\Delta_p^{-kq}(\Delta_p^{kq}z)$ as claimed.

Now, let $\rho$ be a generator of the center of $B(W_g)\simeq \Gg(u,u)$. As $\Delta_p^q(u)\in Z(\Gg(u,u))$, there is some integer $k$ such that $\rho^k=\Delta_p^q(u)$. That is $\rho$ is a $(pk,q)$-regular element in $\Gg$. By applying $\pi_p$, we get that $\pi_p(\rho)$ is a $(pk,q)$-regular element of $B(W)$. By assumption, we have 
\[\pi_p(\Gg(u,u))=C_{B(W)}(\pi_p(\Delta_p^q))\subset C_{B(W)}(\pi_p(\rho)).\]
Furthermore, we have $C_{B(W)}(\pi_p(\Delta^q))\supset C_{B(W)}(\pi_p(\rho))$ since $\pi_p(\Delta_p^q)$ is a power of $\pi_p(\rho)$. The element $\pi_p(\rho)$ is then a $dk$-regular braid in $B(W)$ with the same centralizer as a $d$-regular braid. Since $d$ is maximal with respect to divisibility, we get that $k=\pm 1$. Thus $\rho=\Delta_p^{\pm q}(u)$ which shows the proposition.
\end{proof}

\begin{rem}
The assumption that $d$ is the gcd of the degrees of $W$ which it divides is important, otherwise $Z(\Cc)$ may be generated by some root of $\Delta_p^q$. For instance in the group $W=G_{37}$, the integers $d_1=5$ and $d_2=10$ are regulars. We have $(p_1,q_1)=(1,6)$ and $(p_2,q_2)=(1,3)$. The associated categories of periodic elements are then monoids, given by $C_{M(c)}(\Delta^3)$ and $C_{M(c)}(\Delta^6)$, respectively. Because $d_1$ and $d_2$ both divide $2$ degrees of $W$, those centralizers are equal, and their centers are both equal to $\langle \Delta^3\rangle^+$.
\end{rem}

\begin{lem}\label{lem:finconjatonongreedy}
Let $u\in \Ob(\Cc)$ and let $\lambda(s)$ be an atomic loop of $u$. Let also $t:u\to v$ be an atom. If the path $s^\flat t$ is not greedy, then we have $\lambda(s)t=t\lambda(s')$ for some atomic loop $\lambda(s')$ of $v$.
\end{lem}
\begin{proof}
First, assume that $p=1$. In this case, $\Cc$ has only one object, and $s=s^{\#}=s^{\flat}$. If $t$ is an atom such that $st$ is not greedy, then $st$ is a simple element of a dual braid monoid. It is then balanced and we have $st=ts'$, with $s'$ an atom. 

Assume now that $p>1$. Let $\lambda(s)=s_n\cdots s_1$ where $s_1=s^\flat$ and $s_{i+1}:=s_i^{\flat}$ for $i\in \intv{1,n-1}$ (in particular we have $s_n=s$). We denote $t:=(\alpha,\beta)$ and $s_i:=(a_i,b)$ for $i\in \intv{1,n}$. By Corollary \ref{cor:greedy_normal_form_length_2_in_periodic}, if $s^\flat t$ is not greedy, then we have $\alpha\preceq b$. Let then $x$ be such that $\alpha x=b$. We have
\[s^\flat t=(a_1,b)(\alpha,\beta)=(a_1\alpha,x)=(\alpha,a_1^{\alpha}x)(a_1^{\alpha},x\alpha^{c^\eta}).\]
We set $t_1:=(\alpha,a_1^{\alpha}x)$ and $\sigma_1=(a_1^\alpha,x\alpha^{c^\eta})$. Again, as $\alpha\preceq b$,  the path $s_2t$ is not greedy, and we have
\[s_2t_1=(a_2,b)(\alpha,a_1^{\alpha}x)=(a_2\alpha,x)=(\alpha,a_2^\alpha x)(a_2^\alpha,x\alpha^{c^\eta}),\]
and we set again $t_2:=(\alpha,a_2^\alpha x)$, $\sigma_2:=(a_2^\alpha,x\alpha^{c^\eta})$. By an immediate induction, we get $\lambda(s)t=t_n\sigma_n\cdots \sigma_1$, with $t_n:=(\alpha,a_n^{\alpha}x)$ and $\sigma_i:=(a_i^\alpha,x\alpha^{c^\eta})$ for $i\in \intv{1,n}$.
The two simple morphisms $t$ and $t_n$ share the same source, and we have $\pi_p(t)=\pi_p(t_n)$. By Lemma \ref{lem:preceqindivided}, we have $t=t_n$. For $i\in \intv{1,n-1}$, the target of $\sigma_{i+1}$ is the source of $\sigma_i$, and the second terms of $\sigma_i$ and $\sigma_{i+1}$ are both equal to $x\alpha^{c^{\eta}}$. Thus  $\sigma_{i}$ is equal to $\sigma_{i+1}^{\#}$, again by Lemma \ref{lem:preceqindivided}. 

We claim that $\sigma_n\cdots\sigma_1=\lambda(\sigma_n)$. Since $\sigma_i=\sigma_{i+1}^\#$ for $i\in \intv{1,n-1}$, we only have to show that $\sigma_1$ is the first $\sigma_i$ with target $v$. The target of $\sigma_i$ is $x(\alpha a_i^\alpha)^{c^\eta})=x(a_i\alpha)^{c^{\eta}}$, and the source of $\sigma_n$ is the target of $\sigma_1$, that is $v=x(a_1\alpha)^{c^\eta}$. By cancellativity of $M(c)$, we get that the target of $\sigma_i$ is $v$ only if $i=1$, for $i\in \intv{1,n}$.
\end{proof}

\begin{theo}\label{theo:conj_of_ato_loops}
Let $\lambda(s)\in \Cc(u,u)$ be an atomic loop of some object $u$, and let $f\in \Gg$. If there is some endomorphism $z\in \Cc$ such that $\lambda(s)^nf=fz$ for some $n\geqslant 1$, then $z=\lambda(s')^n$ for some atomic loop $\lambda(s')$ such that $\lambda(s)f=f\lambda(s')$.
\end{theo}
\begin{proof}
We do the proof in several steps depending on the properties of $f$.\begin{enumerate}[1. ]
\item $f$ is an atom. By assumption, we have $f\preceq \lambda(s)^nf$. If $\lambda(s)^nf$ is in greedy normal form, we get $f\preceq s$. Thus $f$ is trivial or $f=s$. We have $\lambda(s)^f=\lambda(s)$ in the first case, and $\lambda(s)^f=\lambda(s^\#)$ in the second case. If $\lambda(s)^nf$ is not in greedy normal form, then $s^\flat f$ is not in greedy normal form. By Lemma \ref{lem:finconjatonongreedy}, we get that $\lambda(s)^f$ is an atomic loop, hence the desired result.
\item $f$ is a simple morphism. Since $\Cc$ is homogeneous, we can proceed by $\succeq$-induction on $f$. The case where $f$ is an atom has already been dealt with. If $f$ is not an atom, then $\lambda(s)^nf$ cannot be in greedy normal form as $f\preceq s$ is impossible. Thus $s^\flat f$ is not greedy, and there is some decomposition $f=tf'$ with $t$ an atom such that $s^\flat t$ is not greedy. By the first point we have that $\lambda(s)^t$ is of the form $\lambda(s')$. Since $z=\lambda(s)^{f}=\lambda(s)^{tf'}=\lambda(s')^{f'}$, the induction hypothesis gives the desired result.
\item $f$ is an arbitrary morphism in $\Gg$. We have that $\Delta_{p}^q(u)$ is central in $\Gg(u,u)$, thus $(\lambda(s)^{n})^{\Delta_{p}^{kq}(u)f}=(\lambda(s)^{n})^{f}=z$ for all $k\in \Z$. By choosing some $k$ big enough so that $kq+\inf(f)\geqslant 0$, we can replace $f$ with $\Delta_{p}^{kq}(u)f$ and assume that $f\in \Cc$. We claim that both $\lambda(s)^n$ and $z$ lies in $\SSS(\lambda(s)^n)$. First, as atomic loops are rigid, we have
\[\cyc(\lambda(s))=(s^\#\cdots s^{\flat\flat}s^\flat)s=\lambda(s^\#)=\dec(\lambda(s)).\]
Powers of atomic loops are also rigid, and we have $\cyc(\lambda(s)^n)=\dec(\lambda(s)^n)=\lambda(s^\#)^n$.
Because of Proposition \ref{prop:sssreachedbycycling}, we get that $\lambda(s)^n$ lies in $\SSS(\lambda(s)^n)$. Then, we have $\inf(z)\geqslant 0$ as $z$ lies in $\Cc$, so $\inf(z)=0$ as $\inf(z)\leqslant \inf(\lambda(s)^n)$. Lastly, we have $\sup(z)\geqslant \sup(\lambda(s)^n)=nm$ (where $m$ is the length of $\lambda(s)$). Since every simple morphism has length at least $1$, and since $\inf(z)=0$, we get $mn=\ell(z)\geqslant \sup(z)\geqslant nm$ and $\sup(z)=nm$, so $z\in \SSS(\lambda(s)^n)$.
Let $s_1\cdots s_r$ be the greedy normal form of $f\in \Cc$. By Proposition \ref{prop:sssconnectedbysimples}, each $(\lambda(s)^{n})^{s_1\cdots s_i}$ lies in $\SSS(\lambda(s)^n)$, in particular is positive. By an immediate induction using the second case, we get that $z$ is of the form $\lambda(t)^n$, with $\lambda(s)^f=\lambda(t)$.
\end{enumerate}
\end{proof}

This theorem is an analogue of \cite[Proposition 2.2]{dmm} in the context of Springer categories. It gives in particular a complete description of super-summit sets of atomic loops and their powers.

\begin{cor}
Let $u\in \Ob(\Gg)$, and let $\sigma\in B(W_g)\simeq \Gg(u,u)$ be a braided reflection. The super-summit set of $\sigma$ in $\Gg$ consists of all the atomic loops to which $\sigma$ is conjugate in $\Gg$. Furthermore, for $n\geqslant 1$, we have
\[\SSS(\sigma^n)=\{\lambda(s)^n~|~ \lambda(s)\in \SSS(\sigma)\}.\]
\end{cor}
\begin{proof}
We already showed that $\sigma$ is conjugate to some atomic loop $\lambda(s)\in \Gg$, and that $\lambda(s)\in \SSS(\lambda(s))$ for all atomic loops $\lambda(s)$. If $g\in \SSS(\lambda(s))$, then $g$ is a positive conjugate of $\lambda(s)$ and Theorem \ref{theo:conj_of_ato_loops} gives that $g$ is an atomic loop.

We also have $\lambda(s)^n\in \SSS(\sigma^n)$. If $g\in \SSS(\sigma^n)$, then we also have that $g$ is a positive conjugate of $\lambda(s)^n$. We get that $g$ is of the form $\lambda(s')^n$ for some conjugate $\lambda(s')$ of $\lambda(s)$ by Theorem \ref{theo:conj_of_ato_loops}.
\end{proof}

Like in \cite{dmm}, our results on the conjugacy of atomic loops can be used to understand the conjugacy of braided reflections and their powers.

\begin{cor}\label{cor:center_finite_index_subgroups}Let $\sigma\in B(W_g)$ be a braided reflection. We have
\[\forall n\geqslant 1,~ C_{B(W_g)}(\sigma^n)=C_{B(W_g)}(\sigma).\]
In particular, if $U\subset B(W_g)$ is a finite index subgroup, then we have $Z(U)\subset Z(B(W_g))$.
\end{cor}
\begin{proof}
Since an element commuting with $\sigma$ commutes with any of its powers, we only have to show the inclusion $\subset$. We fix $u\in \Ob(\Gg)$ so that $B(W_g)$ is identified with $\Gg(u,u)$. By Theorem \ref{prop:braidreflectionconjugatetoatomicloop}, there is some morphism $f:u\to v$ in $\Gg$ such that $\sigma^f=\lambda(s)$ is an atomic loop. 

Let $x\in \Gg(u,u)$ which commutes with $\sigma^n$. We define $x':=x^f$, which commutes with $\lambda(s)^n$. By Theorem \ref{theo:conj_of_ato_loops}, since $x'^{-1}\lambda(s)^{n}x'=\lambda(s)^n$, we have $\lambda(s)x'=x'\lambda(s')$ for some atomic loop $\lambda(s')$ such that $\lambda(s')^n=\lambda(s)^n$. Since atomic loops are rigid, the equality $\lambda(s)^n=\lambda(s')^n$ implies $\lambda(s)=\lambda(s')$ and $x'$ commutes with $\lambda(s)$. We deduce that $x$ commuted with $\sigma$ as we wanted to show.

Now, let $U\subset B(W_g)$ be a finite index subgroup, and let $z\in Z(U)$. If $\sigma$ is a braided reflection in $B(W_g)$, then some nontrivial power $\sigma^n$ of $\sigma$ lies in $U$ as it has finite index. Since $z$ commutes with $\sigma^n$, it commutes with $\sigma$. We have shown that $z$ commutes with every braided reflection in $B(W_g)$. As braided reflections generate $B(W_g)$, we have $z\in Z(B(W_g))$ as we wanted to show.
\end{proof}

The above corollary allows us to complete the proof of the following general result on complex braid groups. 

\begin{cor}\label{cor:centralizer_braided_reflection}
Let $W$ be a complex reflection group and let $B$ be its associated complex braid group. Let also $\sigma\in B$ be a braided reflection. We have
\[\forall n\geqslant 1,~C_B(\sigma^n)=C_B(\sigma).\]
\end{cor}
\begin{proof}
If $W$ is not irreducible, it decomposes as a product $W=W_1\times\cdots\times W_k$ of irreducible complex reflection groups. Likewise, the braid group $B$ decomposes as $B=B_1\times\cdots \times B_k$, where $B_i=B(W_i)$ for $i\in \intv{1,k}$. The set of braided reflections of $B$ is the union of the sets of braided reflections of the $B_i$. The theorem is then an easy consequence of the case of an irreducible group.

Assume that $W$ is irreducible. The proof of \cite[Corollary 2.5]{dmm} applies to a Garside monoid $M$ satisfying particular assumptions (listed in the beginning of \cite[Section 2]{dmm}). It shows that the centralizer in the enveloping group $G(M)$ of any atom of $M$ coincides with that of its nontrivial powers. If there is an isomorphism $B\simeq G(M)$ which maps braided reflections to (conjugates of) atoms in $G(M)$, then we have the result. As pointed out in the proof of \cite[Corollary 2.7]{dmm}, this is the case when $W$ is well-generated, as well as when $W=G_{12}$ or $G_{22}$.

Assume now that $W$ is badly-generated and different from these two groups. If $W$ belongs to the infinite series and is badly generated, then we have $W=G(m,p,n)$ with $m>1$ and some divisor $p$ of $m$. By \cite[Theorem 3.6 and Proposition 3.8]{bmr}, there is an injective morphism from $B$ to the braid group of the well-generated complex reflection group $G(m,1,n)$. Moreover, this morphism maps braided reflections to powers of braided reflections. The result is then an immediate consequence of the well-generated case. 

Now, by \cite[Theorem 2.2]{beskpi1}, it is sufficient to consider the case where $W$ contains only reflections of order $2$. The only two groups $W$ which contains only reflections of order $2$, are badly generated, and do not belong to the infinite series are $G_{13}$ and $G_{31}$.

For $W=G_{13}$, consider the Artin-Tits monoid of type $I_2(6)$:
\[M:=\langle a,b~|~ababab=bababa\rangle^+.\]
It satisfies the assumptions required for the proof of \cite[Corollary 2.5]{dmm} to hold. By \cite[Section 3.5]{paratresses}, there is an isomorphism $B(G_{13})\to G(M)$, which sends representatives of conjugacy classes of braided reflections to $\sigma_1:=b^{-1}$ and $\sigma_2:=bababa^{-1}$. Since $\sigma_1$ is a nontrivial power of an atom of $M$, we know that the centralizers of $\sigma_1^n$ and $\sigma_1$ coincide for $n\geqslant 1$. For $\sigma_2$, notice that $D=ababab$ is central in $G(M)$, and that $\sigma_2=Da^{-2}$. We then have
\[\forall n\geqslant 1,~ C_{G(M)}(\sigma_2^n)=C_{G(M)}(D^na^{-2})=C_{G(M)}(a^{-2n})=C_{G(M)}(a^{-2})=C_{G(M)}(Da^{-2})=C_{G(M)}(\sigma_2),\]
which finishes the proof in this case.

Lastly, we consider the case of $W=G_{31}$. As in Example \ref{ex:=g31_dans_g37}, we can see $W$ as the centralizer of a $i$-regular element in the well-generated group $G_{37}$ and we can consider the attached Springer groupoid $\Gg$. The result is then a direct application of Corollary \ref{cor:center_finite_index_subgroups}.
\end{proof}

This result is slightly stronger than \cite[Theorem 1.4]{dmm}. It implies in particular that any subgroup of a complex braid group $B$ containing nontrivial powers of braided reflections has its center included in $Z(B)$. Of course, this stronger statement was proven for many cases in \cite{dmm}, and we merely completed the proof for the few remaining cases.

\section{Presentations of the braid group $B(G_{31})$}\label{sec:presentation_of_b31}
The goal of this section is to apply the previous results to the complex braid group $B(G_{31})$. In particular we obtain several positive homogeneous presentations of this group, where the generators are braided reflections.

In this section, we consider the irreducible complex reflection group $W:=G_{37}$. Its highest degree is $h:=30$. If $c$ is a Coxeter element in $W$, then we denote by $I_c$ the interval between $1$ and $c$ in $W$ for the reflection length. We keep the convention of Section \ref{sec:properties_of_springer_categories}: we denote the elements of $I_c$ in normal font and products of elements of $I_c$ should be understood as products in the dual braid monoid $M(c)$.

The integer $d=4$ is regular for $W$, and the centralizer of a $i$-regular element in $G_{37}$ is a group of type $G_{31}$ (see Example \ref{ex:=g31_dans_g37}). Following the notation of previous sections, we set $p=2$ and $q=15$. We denote by $\Bb_{31}$ (resp. by $\Cc_{31}$) the Springer groupoid (resp. the Springer category) attached to $(W,c,p,q)$. The Garside map of $\Cc_{31}$ will be denoted by $\Delta$, and its Garside automorphism will be denoted by $\phi$. By Theorem \ref{theo:groupoidtressequivgroupoiddiv}, the groupoid $\Bb_{31}$ is equivalent to the complex braid group $B(G_{31})$.

The integers $\eta:=8$ and $\mu:=1$ are such that $p\eta-q\mu=1$. As $c^{15}=-\Id$ is central in $W$, we have $I_c^q=I_c$. Following the beginning of Section \ref{sec:3.1}, we have
\begin{align*}
\Ob(\Cc_{31})&:=\Dd_1= \{u\in I_c~|~ u u^{c^8}=c\text{~and~}\ell_R(u)=4\},\\
\Ss_2^{15}&:=\Dd_2=\{(a,b)\in (I_c)^2~|~ ab\in \Dd_1\},\\
R_2^{15}&:=\Dd_3=\{(x,y,z)\in (I_c)^3~|~ xyz\in \Dd_1\}.
\end{align*}
By \cite[Theorem 13.3]{beskpi1}, the category $\Cc_{31}$ has $88$ objects, $2691$ simple morphisms, $16359$ defining relations. 

By Lemma \ref{lem:atoms_of_periodic_cat}, the atoms of $\Cc_{31}$ are exactly the elements of length $1$. Since $4$ is the gcd of the degrees of $G_{37}$ which are divisible by $4$, we can apply Proposition \ref{prop:center_of_categories} and Corollary \ref{cor:center_finite_index_subgroups}:

\begin{theo}\label{theo:centres_b31_p31}
The centers of $\Cc_{31}$ and $\Bb_{31}$ are cyclic and generated by $\Delta^{15}$. If $U$ is a finite index subgroup of $B(G_{31})$, then $Z(U)\subset Z(B(G_{31}))$. In particular, the center of the pure braid group $P(G_{31})$ is cyclic and generated by the full-twist.
\end{theo}
\begin{proof}
The only nontrivial part is that the full-twist is a generator of $Z(P(G_{31}))$. Let us fix an object $u\in \Ob(\Bb_{31})$ so that we have $B(G_{31})\simeq \Bb_{31}(u,u)$. The center of $P(G_{31})$ is cyclic and generated by the smallest power of $\Delta^{15}(u)$ which lies in it. Since the collapse functor $\Cc_{31}\to M(c)$ sends $\Delta^{15}(u_0)$ to some element $g$ with $g^2=c^{15}$, the smallest power of $\Delta^{15}$ lying in $P(G_{31})$ is the full-twist $\Delta^{60}$.
\end{proof}

The remainder of this section is devoted to the study of presentations of $B(G_{31})$. If $u$ is an object of $\Cc_{31}$, then we have an isomorphism $\Bb_{31}(u,u)\simeq B(G_{31})$ which sends atomic loops to braided reflections.

On the one hand, we construct a conjectural presentation of $\Bb_{31}(u,u)$ with atomic loops as generators. On the other hand, the Reidemeister-Schreier method for groupoids (cf. Appendix \ref{app:reidemeisterschreier}) gives a presentation of $\Bb_{31}(u,u)$ which we know holds. We then prove that these two presentations are equivalent.

\subsection{The method}
Let $u$ be an object of $\Bb_{31}$. We start by considering the submonoid $L_u^+$ of $\Cc_{31}(u,u)$ generated by atomic loops of $u$. Our first goal is to construct a (conjectural) group presentation using $L_u^+$. We consider the following algorithm.

\begin{algorithm}\label{alg:shortest_right_multiple_atomic_loops}\caption{Compute shortest right-multiple of atomic loops in $L_u^+$}
\begin{algorithmic}
\STATE \textbf{Input:} Two atomics loops $\lambda(s)$ and $\lambda(t)$ of $u$.
\STATE \textbf{Output:} If $\lambda(s),\lambda(t)$ admit a common right-multiple in $L_u^+$, then the output is a pair of words $\theta(\lambda(s),\lambda(t)),\theta(\lambda(t),\lambda(s))$ in $L_u^+$ such that $\lambda(s)\theta(\lambda(s),\lambda(t))=\lambda(t)\theta(\lambda(t),\lambda(s))$. No output otherwise.
\STATE put $i:=1$
\STATE compute the set $S_i$ of words of length $i$ in $L_u^+$
\WHILE{$\lambda(s)m_1\neq \lambda(s)m_2$ in $\Cc_{31}$ for all $(m_1,m_2)\in S_i\times S_i$}
\STATE put $i:=i+1$
\STATE compute the set $S_i$ of words of length $i$ in $L_u^+$
\ENDWHILE
\STATE put $S:=\{(m_1,m_2)\in S_i\times S_i~|~ \lambda(s)m_1=\lambda(t)m_2\in \Cc_{31}\}$
\RETURN $(m_1,m_2)\in S$ such that $m_1$ is the lowest possible in the lexicographic order, and $m_2$ is the lowest possible in the lexicographic order among the words $m$ such that $(m_1,m)\in S$.
\end{algorithmic}
\end{algorithm}

Algorithm \ref{alg:shortest_right_multiple_atomic_loops}, running on two atomic loops of $u$, terminates if and only if they admit a common right-multiple in $L_u^+$. The fact that it terminates on every pair of atomic loops of $u$ (which we checked by computer using the data of Section \ref{sec:computational_data}) proves that all pairs of atomic loops of $u$ admit a common right-multiple. We now consider the following presentation:
\begin{itemize}
\item The set of generators is a set $X_u:=\{\lambda'(s)\}$ in bijection with atomic loops of $u$ in $L_u^+\subset \Cc_{31}(u,u)$.
\item The set $R_u$ is given by relations of the form
\[\lambda'(s)\theta(\lambda'(s),\lambda'(t))=\lambda'(t)\theta(\lambda'(t),\lambda'(s)),\]
where $\theta(\lambda'(s),\lambda'(t))$ is the image of the word $\theta(\lambda(s),\lambda(t))$ under the natural bijection between the atomic loops of $u$ and the set $X_u$.
\end{itemize}
We define $H_u:=\langle X_u~|~R_u\rangle$ and $H_u^+:=\langle X_u~|~R_u\rangle^+$. By \cite[Lemma II.4.3]{ddgkm}, the presentation of $H_u^+$ is \nit{right-complemented} in the sense of \cite[Definition II.4.1]{ddgkm}. That is for $a\neq b\in X_u$, there is exactly one relation of the form $a\ldots=b\ldots$, namely $a\theta(a,b)=b\theta(b,a)$.

 The goal of this section is to show the following theorem:

\begin{theo}\label{theo:isomorphism_presentation_b31}
Let $u\in \Ob(\Cc_{31})$. The natural map from $X_u$ to the set of atomic loops of $u$ induces a group isomorphism $H_u\simeq \Bb_{31}(u,u)$. In particular, $\langle X_u~|~ R_u\rangle$ gives a positive homogeneous presentations of $B(G_{31})$ with braided reflections as generators.
\end{theo}

We first notice that, by definition of $R_u$ and $X_u$, the natural map from $X_u$ to the set of atomic loops of $u$ induces a group morphism $f_u:H_u\to\Bb_{31}(u,u)$. We only have to show that the said morphism is an isomorphism. We show this by case by case analysis on the objects of $\Cc_{31}$.

\begin{lem}\label{lem:il_suffit_de_considerer_les_phi_orbites}Let $u\in \Ob(\Cc_{31})$. If $f_u$ is an isomorphism, then $f_{\phi(u)}$ is also an isomorphism. In particular we only need to show that Theorem \ref{theo:isomorphism_presentation_b31} holds for a system of representatives of $\phi$-orbits in $\Bb_{31}$.
\end{lem}
\begin{proof}
The automorphism $\phi$ gives an isomorphism between $\Cc_{31}(u,u)$ and $\Cc_{31}(\phi(u),\phi(u))$. Because of Lemma \ref{lem:properties_of_diese_and_flat}, $\phi$ induces a bijection between the sets of atomic loops of $u$ and of $\phi(u)$. We obtain that $\phi$ induces bijections between $X_u$ and $X_{\phi(u)}$, and between $R_u$ and $R_{\phi(u)}$, respectively. Thus $H_u\simeq H_{\phi(u)}$ which shows the claim.
\end{proof}

Note that the monoid $H_u^+$ is homogeneous by definition. In particular we have a solution to the word problem in $H_u^+$, given by Algorithm \ref{alg:word_problem_homogeneous}.

\begin{algorithm}
\caption{Check equality between two words in $H_u^+$}
\label{alg:word_problem_homogeneous}
\begin{algorithmic}
\STATE \textbf{Input:} Two words $m_1,m_2$ in the atoms of $H_u^+$
\STATE \textbf{Output:} \textbf{true} if $m_1$ and $m_2$ represent the same element of $H_u^+$, \textbf{false} otherwise.
\STATE put $S:=\varnothing$.
\STATE put $S_1:=\{m_1\}$.
\WHILE{$S\neq S_1$} 
\STATE put $S:=S_1$.
\STATE add to $S_1$ the set of words obtained from elements of $S_1$ by applying one relation of $R_u$.
\ENDWHILE
\IF{$m_2\in S$}
\RETURN \textbf{true}
\ELSE
\RETURN \textbf{false}
\ENDIF
\end{algorithmic}
\end{algorithm}

Since $H_u^+$ is homogeneous, two words representing the same element have the same length, thus there is a finite number of words that represent the same element of $H_u^+$ and Algorithm \ref{alg:word_problem_homogeneous} always terminates.

In order to prove Theorem \ref{theo:isomorphism_presentation_b31} for the object $u$, we first compute a presentation of $\Bb_{31}(u,u)$ by the Reidemeister-Schreier method for groupoids (see Appendix \ref{app:reidemeisterschreier}). We start with a presentation of the groupoid $\Bb_{31}$, for instance that of Theorem \ref{prop:hurwitzpresentationforc31}. 

As we want atomic loops to appear as generators, we need to choose the Schreier transversal accordingly. We use the following lemmas concerning atomic loops in the category $\Cc_{31}$.

\begin{lem}\label{lem:atomicloopshavelength2}Let $s$ be an atom of $\Cc_{31}$, the atomic loop $\lambda(s)$ has length two in $\Cc_{31}$.
\end{lem}
\begin{proof}
Let $s:=(a,b)$. We denote by $u$ the source of $s$, and by $v$ its target. The morphism $s^\#$ is given by $(a,b)^\#=(a^{c^8b^{-1}},b)$. We know that 
\[a^{c^7}b^{c^7}ab=c\Rightarrow c^7a^{c^7}b^{c^7}a=abc^7a= c^8b^{-1},\] 
and so $a^{c^8b^{-1}}=a^{abc^7a}=a^{(ba^{c^8})c^7}=a^{vc^7}$. The morphism $s^{\#\#}$ is then given by $((a^{vc^7})^{avc^7},b)$. Because of Lemma \ref{lem:reflectionsetconjcommute}, $a$ and $a^{vc^7}$ commute, so $s^{\#\#}=(a^{(vc^7)^2},b)$. Since $(vc^7)^2=c^{15}$ is central, we have $s^{\#\#}=s$ as claimed.
\end{proof}

\begin{lem}
Let $u$ be an object of $\Bb_{31}$. There is a Schreier transversal $T$ rooted in $u$ and containing all atoms with source $u$. In particular the atomic loops of $u$ appear as generators of the presentation of $\Bb_{31}(u,u)$ induced by $T$.
\end{lem}

\begin{proof}First, note that if a Schreier transversal $T$ contains an atom $s$ with source $u$, then Lemma \ref{lem:atomicloopshavelength2} gives that $\gamma(s^\#)=ss^\#=\lambda(s)$ with the notation of Lemma \ref{lem:lemme_de_schreier}.

Now, thanks to Proposition \ref{cor:no_parallel_simples}, all atoms with source $u$ have different targets. We can thus consider a Schreier transversal $T$ rooted in $u$ and containing all atoms with source $u$.
\end{proof}

\begin{rem}
Let $u\in \Ob(\Cc_{31})$. If $s,t$ are two composable atoms with $s\in \Cc_{31}(u,-)$, then by Lemma \ref{lem:atomicloopshavelength2}, we can consider the product $m:=stt^{\#}s^{\#}\in \Cc_{31}(u,u)$. If this product is in greedy normal form, then no atomic loop of $u$ left-divides $m$, and thus $m$ is not an element of $L_u^+$. We can check by computer using the data of Section \ref{sec:computational_data} that this occurs for every object of $\Cc_{31}$, and thus $L_u^+$ is always a strict submonoid of $\Cc_{31}(u,u)$.

Since the atoms of the monoid $H_u^+$ are in bijection with atomic loops of $u$, we also deduce that $H_u^+$ is not isomorphic to the monoid $\Cc_{31}(u,u)$. As a matter of fact, we will see that the monoid $H_u^+$ is never cancellative, and thus cannot be isomorphic to either $L_u^+$ or to $\Cc_{31}(u,u)$. The isomorphism of Theorem \ref{theo:isomorphism_presentation_b31} only occurs at the level of groups.
\end{rem}

Let $T$ be a Schreier transversal rooted in $u$ and containing all atoms with source $u$. Let $\langle S^*~|~R^*\rangle$ denote the presentation of $\Bb_{31}(u,u)$ obtained by the Reidemeister-Schreier method applied to $T$ and to the presentation of Theorem \ref{prop:hurwitzpresentationforc31}. Of course, the presentation $\langle S^*~|~ R^*\rangle$ is quite redundant. We first want to show that every element of $S^*$ can be expressed as a word in the atomic loops. For this we repeatedly apply Tietze transformations, as in Algorithm \ref{alg:bethany}.

\begin{algorithm}\caption{Reduction of generators}\label{alg:bethany}
\begin{algorithmic}
\STATE \textbf{Input:} A group presentation $\langle S~|~R\rangle$ and a subset $S'$ of $S$
\STATE \textbf{Output:} A group presentation $\langle S'~|~R'\rangle$ equivalent to the first by Tietze transformation, or no output.
\WHILE{$S'\neq S$}
\STATE choose $r\in R$ a relator with only one letter $a$ not belonging to $S'\cup S'^{-1}$
\STATE replace in $R$ every occurrence of $a$ by its expression in $S'$ using the relator $r$.
\STATE remove the relator $r$ from $R$
\STATE remove the letter $a$ from $S$.
\ENDWHILE
\RETURN the presentation $\langle S'~|~ R\rangle$
\end{algorithmic}
\end{algorithm}

The fact that this algorithm terminates for each object of $\Bb_{31}$, which is again checked by computer, proves the following result:

\begin{prop}
Let $u$ be an object of $\Bb_{31}$. The atomic loops of $u$ generate the group $\Bb_{31}(u,u)$. In particular the natural morphism $H_u\to \Bb_{31}(u,u)$ is surjective.
\end{prop}

By applying Algorithm \ref{alg:bethany} to the presentation $\langle S^*~|~ R^*\rangle$, we obtain a presentation $\langle X_u~|~ R'_u\rangle$ of the group $\Bb_{31}(u,u)$. In order to prove Theorem \ref{theo:isomorphism_presentation_b31} for the object $u$, it is sufficient to prove that every relator of $R'_u$ is in fact trivial in $H_u$. This will prove that the morphism $H_u\to \Bb_{31}(u,u)$ is injective.

Because the defining presentation of $H_u^+$ is right-complemented, we can consider a \nit{right-reversing} algorithm in the sense of \cite[Definition 4.21]{ddgkm}. The main idea is that the relation $a\theta(a,b)=b\theta(b,a)$ implies that the equality $a^{-1}b=\theta(a,b)\theta(b,a)^{-1}$ holds. We can use these type of relation in order to simplify words in $X_u\cup X_u^{-1}$.

\begin{algorithm}\caption{Right-reversing in $H_u^+$ (\cite[Algorithm II.4.33]{ddgkm})}\label{alg:right-reversing_Hu}
\begin{algorithmic}
\STATE \textbf{Input:} A word $w$ in $X_u\cup X_u^{-1}$
\STATE \textbf{Output:} A fraction $fg^{-1}$ with $f,g\in H_u^+$ representing $w$ in $H_u$, or no output.
\WHILE{there is some subword of the form $a^{-1}b$ in $w$}
\STATE put $j$ the position in $w$ of the first subword of the form $a^{-1}b$ in $w$
\IF{$a=b$}
\STATE remove the subword $a^{-1}b$ from $w$
\ELSE 
\STATE replace $a^{-1}b$ with $\theta(a,b)\theta(b,a)^{-1}$ in $w$ at position $j$
\ENDIF
\ENDWHILE
\RETURN $w$.
\end{algorithmic}
\end{algorithm}

Algorithm \ref{alg:right-reversing_Hu} does not always terminate: it may loop indefinitely. If it terminates, its output is a fraction, and we should check that it is trivial. Of course this may be quite long as algorithm \ref{alg:word_problem_homogeneous} is far from optimal. This process can be sped up by ``partially simplifying'' at each step. The solution to the word problem given by Algorithm \ref{alg:word_problem_homogeneous} allows for the computation of longest common divisors of elements of $H_u^+$. A general word in $X_u\cup X_u^{-1}$ can be written as a product of (short) fractions. We can simplify these fractions at each step of Algorithm \ref{alg:right-reversing_Hu}.

We checked by computer that Algorithm \ref{alg:right-reversing_Hu} terminates for every relator in $R_u'$, that is every relator can be expressed as a right-fraction of elements of $H_u^+$. Finally, we use the following algorithm to prove that every relator, written as a fraction, is trivial in $H_u$.

\begin{algorithm}
\caption{Partial solution to the word problem in $H_u$}\label{alg:partial_solution_word_prblm_homogeneous}
\begin{algorithmic}
\STATE \textbf{Input:} A fraction $fg^{-1}$ with $f,g\in H_u^+$.
\STATE \textbf{Output:} {\tt true} If there is some $n\in H_u^+$ such that $fn=gn$. No output otherwise.
\STATE put $i:=1$
\STATE Compute the set $S_i$ of words of length $i$ in $X_u$
\WHILE{$fn \neq gn$ for all $n\in S_i$}
\STATE put $i:=i+1$
\STATE Compute the set $S_i$ of words of length $i$ in $X_u$
\ENDWHILE
\RETURN true
\end{algorithmic}
\end{algorithm}

This algorithm is useful since the monoid $H_u^+$ may not be cancellative: we can have $fn=gn$ (and thus $fg^{-1}=1$ in $H_u$) without having $f=g$ in $H_u^+$.

The fact that Algorithm \ref{alg:partial_solution_word_prblm_homogeneous} returns {\tt true} for all relators of the presentation of $\Bb_{31}(u,u)$ finally proves that Theorem \ref{theo:isomorphism_presentation_b31} holds for $u$.

\subsection{Computational data}\label{sec:computational_data}
We consider the following elements in $\C^8$

\[\begin{cases}
\alpha_{1}=\frac{1}{2}(1,-1,-1,-1,-1,-1,-1,1),\\ 
\alpha_{2}=(1,1,0,0,0,0,0,0),\\ 
\alpha_{3}=(-1,1,0,0,0,0,0,0),\\ 
\alpha_{4}=(0,-1,1,0,0,0,0,0),\\ 
\alpha_{5}=(0,0,-1,1,0,0,0,0),\\ 
\alpha_{6}=(0,0,0,-1,1,0,0,0),\\ 
\alpha_{7}=(0,0,0,0,-1,1,0,0),\\ 
\alpha_{8}=(0,0,0,0,0,-1,1,0),\\ 
\alpha_{9}=(1,0,1,0,0,0,0,0),
\end{cases}
\begin{cases}
\alpha_{10}=(0,0,-1,0,1,0,0,0),\\ 
\alpha_{11}=(1,0,0,0,0,1,0,0),\\ 
\alpha_{12}=\frac{1}{2}(-1,-1,-1,-1,-1,1,-1,1),\\ 
\alpha_{13}=(0,1,0,0,0,1,0,0),\\ 
\alpha_{14}=\frac{1}{2}(-1,-1,-1,-1,-1,-1,1,1),\\ 
\alpha_{15}=(0,1,0,0,0,0,1,0),\\ 
\alpha_{16}=\frac{1}{2}(1,-1,1,-1,-1,-1,1,1),\\ 
\alpha_{17}=\frac{1}{2}(1,-1,-1,1,-1,-1,1,1),\\ 
\alpha_{18}=\frac{1}{2}(-1,1,-1,-1,1,-1,1,1).
\end{cases}\]
For each $i\in \intv{1,18}$, we consider the reflection $s_i$ of $\C^8$ given by $s_i(x)=x-2\frac{\scal{x}{\alpha_i}}{\scal{\alpha_i}{\alpha_i}}\alpha_i$, where $\scal{.}{.}$ is the standard hermitian product on $\C^8$. The set $s_1,\ldots,s_8$ generates a subgroup $W$ of $\GL_n(\C)$ isomorphic to the complex reflection group $G_{37}$ and which contains all the $s_i$ for $i\in \intv{1,18}$. A Coxeter element of $W$ is given by $c=s_1s_2s_3s_4s_5s_6s_7s_8$.



A system of representatives of $\phi$-orbits of objects of the category $\Cc_{31}$ is given by the following elements of $W$:

\[\begin{cases} u_{1}=s_{10}s_{12}s_{13}s_{18}, \\ u_{2}=s_9s_{10}s_{12}s_{16}, \\ u_{3}=s_3s_{11}s_{12}s_{16}, \\ u_{4}=s_7s_8s_{14}s_{15}, \end{cases} \begin{cases} u_{5}=s_3s_7s_{11}s_{12}, \\ u_{6}=s_3s_4s_7s_{11},\\u_7=s_3s_4s_7s_{17},\\u_{8}=s_3s_7s_{12}s_{17}. \end{cases}.\]

As we want links between the various presentations we obtain for each representative, we give explicit isomorphisms between the different groups $\Bb_{31}(u_i,u_i)$. For this we use the following graph in $\Bb_{31}$:
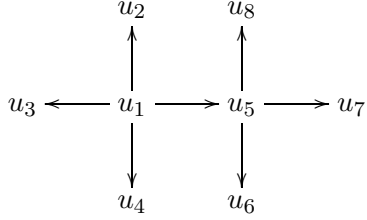
\[\xymatrix{ & u_{2} & u_{8} \\
u_{3} & u_{1} \ar[l] \ar[u] \ar[d]\ar[r] & u_{5} \ar[u] \ar[d] \ar[r] & u_7 \\ & u_{4} & u_{6}} \]

where each arrow is a simple morphism in $\Cc_{31}$ (because of Proposition \ref{cor:no_parallel_simples}, a simple morphism is uniquely determined by its source and target). For each $i,j\in \intv{1,8}$, this graph induces a well defined isomorphism $\varphi_{i,j}:\Bb_{31}(u_i,u_i)\to \Bb_{31}(u_j,u_j)$ which preserves braided reflections.

For $i,j\in \intv{1,8}$, we have by definition $\varphi_{i,j}=\varphi_{i,1}\varphi_{1,j}$, so we only need to describe morphisms of the form $\varphi_{i,1}$ and $\varphi_{1,i}$ for $i\in \intv{1,8}$. 

In the case of the orbit of $u_1$, we give expressions of the atomic loops in the generators $\sigma_1,\ldots,\sigma_8$ of the Artin group associated with $W$. Replacing $\sigma_1,\ldots,\sigma_8$ with $s_1,\ldots,s_8$ gives a set of elements in $W$ which generate a copy of $G_{31}$. We also give a family of vectors in $\C^4$ such that the $2$-reflections associated with the orthogonal hyperplanes of these vectors (in the standard hermitian product) generate a group isomorphic to $G_{31}$.

Furthermore, we know that the full-twist in $\Bb_{31}(u,u)$ is given by $\Delta^{60}(u)$ by Lemma \ref{futw}. By \cite[Theorem 1.2 and Proposition 8.1]{regularbraids}, every root of the full-twist in $B(G_{31})$ is conjugate to a power of either a $20$-th root or a $24$-th root of the full-twist. Furthermore, the full-twist admits $20$-th roots and $24$-th roots. For each presentation we obtain, we give an explicit $20$-th root (resp. $24$-th root) of the full-twist as a word in the generators. If $\rho\in \Bb_{31}(u,u)$ is a $24$-th root of the full-twist $\Delta^{60}(u)$, then $\rho^6$ is a $4$-th root of $\Delta^{60}(u)$. By Remark \ref{rem:pqcoprime} and Theorem \ref{theo:centres_b31_p31}, we get that $\rho^6=\Delta^{15}(u)$ is a generator of $Z(\Bb_{31}(u,u))$.

\subsection{Presentation associated with representatives of the $\phi$-orbits}
\subsubsection{Orbit of $u_1$}\label{sec:4.3.1}
For the first orbit, we recover the presentation of $B(G_{31})$ conjectured in \cite[Table 3]{bmr} and \cite[Conjecture 2.4]{bessismichel}. The object $u_1$ has $5$ atomic loops $s,t,u,v,w$. The relations we obtain are as follows
\[\begin{cases} ts=st,~vt=tv,~wv=vw,\\suw=uws=wsu,\\  svs=vsv,~vuv=uvu,~utu=tut,~twt=wtw. \end{cases}\]
This presentation is usually represented by the following diagram (corresponding to the Broué-Malle-Rouquier diagram for the reflection group $G_{31}$).

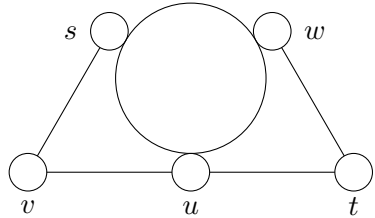
\spa
\begin{center}
\begin{tikzpicture}[scale=0.25]
\draw (0,0) circle (4);

\draw (30:5) circle (1);
\node[right] at (24:6) {$w$};
\node[left] at (156:6) {$s$};
\node[below] at (270:6) {$u$};
\node[below] at (-8.66,-5.96) {$v$};
\node[below] at (8.66,-5.8) {$t$};

\draw (150:5) circle (1);
\draw (270:5) circle (1);
\draw (-8.66,-5) circle (1);
\draw (8.66,-5) circle (1);
\draw (-7.66,-5)--(-1,-5);
\draw (1,-5)--(7.66,-5);
\draw (-8.16,-4.134)--(-4.83,1.634);
\draw (8.16,-4.134)--(4.83,1.634);
\end{tikzpicture}\end{center}

In the Artin group associated with $W$, the atomic loops $s,t,u,v,w$ can be expressed as

\begin{align*}
s&:=(\sigma_1\sigma_4)^{\sigma_2\sigma_3\sigma_1\sigma_4\sigma_5\sigma_4\sigma_2\sigma_3\sigma_4\sigma_5\sigma_6\sigma_7\sigma_8},\\
t&:=(\sigma_4\sigma_2)^{\sigma_2\sigma_3\sigma_4\sigma_5\sigma_6\sigma_7},\\
u&:=(\sigma_4\sigma_2)^{\sigma_2\sigma_3\sigma_1\sigma_4\sigma_5\sigma_4\sigma_2\sigma_3\sigma_4\sigma_5\sigma_6\sigma_5\sigma_7\sigma_6\sigma_8\sigma_7\sigma_6},\\
v&:=(\sigma_1\sigma_3)^{\sigma_3\sigma_4\sigma_5\sigma_4\sigma_6\sigma_7},\\
w&:=(\sigma_2\sigma_3)^{\sigma_1\sigma_4\sigma_2\sigma_3\sigma_1\sigma_4\sigma_5\sigma_4\sigma_2\sigma_3\sigma_1\sigma_4\sigma_5\sigma_6\sigma_5\sigma_4\sigma_2\sigma_3\sigma_1\sigma_4\sigma_5\sigma_6\sigma_7\sigma_6\sigma_5\sigma_4\sigma_2\sigma_3\sigma_4\sigma_5\sigma_6\sigma_7\sigma_8\sigma_7\sigma_6}.
\end{align*}

In $\C^4$, the $2$-reflections associated with the following roots (in the standard hermitian product) generate a subgroup of $\GL_4(\C)$ which is isomorphic to $G_{31}$.

\[\alpha_s:=\frac{1}{2}\matrix{2\\1+i\\-1-i\\0},~\alpha_t:=\matrix{1\\-1\\i\\-1},~\alpha_u:=\matrix{1\\0\\-1-i\\-i},~\alpha_v:=\frac{1}{2}\matrix{2\\1-i\\1-i\\0},~\alpha_w:=\matrix{1\\i\\-1\\-1}.\]

The monoid $H_{u_1}^+$ given by the above presentation is not cancellative: we have $tuwtuw\neq uwtuwt$ and $stuwtuw=suwtuwt$ in $H_{u_1}^+$. Thus $H_{u_1}^+$ cannot be a Garside monoid.

The submonoid $L_{u_1}^+$ of $\Cc_{31}(u_1,u_1)$ generated by $s,t,u,v,w$ is cancellative, but it does not admit right-lcms: we know that $utu=tut$ is the shortest common multiple of $t$ and $u$. If right-lcms exists in $L_{u_1}^+$, then $tut$ must be the right-lcm of $t$ and $u$. We would then have 
\[tut\preceq tuwtuw\Rightarrow t\preceq wtuw\]
in $L_{u_1}^+$ by cancellativity. We can check that $t\preceq wtuw$ does not hold in $L_{u_1}^+$: the element $t^{-1}wtuw$ does lie in $\Cc_{31}(u_1,u_1)$, but not in the submonoid $L_{u_1}^+$.

Maximal roots of the full-twist are given by
\[\begin{array}{r||c|c}\text{Maximal regular number }d & 20 & 24  \\\hline d\text{-th root of }\Delta^{60}(u_1)& stuvws& stuvw\end{array}\]

In particular, we get that $(stuvw)^6=\Delta^{15}(u_1)$ generates $Z(\Bb_{31}(u_1,u_1))$.

\subsubsection{Orbit of $u_{2}$} The object $u_2$ admits $7$ atomic loops $a,b,c,s,t,v,w$, with relations as follows:
\[\begin{cases} sv=vb=bs,~ av=vc=ca,\\
wv=vw,~st=ts,~tv=vt,~tb=bt,\\
was=swa=asw,~wcb=bwc=cbw,\\
wtw=twt,~ata=tat,~aba=bab,~tct=ctc,\\ swav=cabw.
\end{cases}\]
The last line of relations can be omitted, as it is implied by the others. The monoid $H_{u_2}^+$ given by this presentation is not cancellative: we have $bwab\neq abwa$ and $ cbwab= cabwa$ in $H_{u_2}^+$. Thus $H_{u_2}^+$ cannot be a Garside monoid.

The submonoid $L_{u_2}^+$ of $\Cc_{31}(u_2,u_2)$ generated by the atomic loops is cancellative, but it does not admit right-lcms: we know that $aba=bab$ is the shortest common multiple of $a$ and $b$. If right-lcms exists in $L_{u_2}^+$, then $aba$ must be the right-lcm of $a$ and $b$. We would then have
\[aba\preceq abwa\Rightarrow a\preceq wa,\]
which does not hold in $L_{u_2}^+$ (it does not even hold in $\Cc_{31}(u_2,u_2)$).

The relations defining $H_{u_2}^+$ give in particular $b=v^{-1}sv=s^v$ and $c=a^v$ in $H_u$. By deleting these generators, we get that $\Bb_{31}(u_{2},u_{2})$ is generated by $a,s,t,v,w$, with relations as follows:
\[\begin{cases}   wv=vw,~st=ts,~vt=tv, \\swa=was=asw,\\ twt=wtw,~ata=tat,~vav=ava,~svs=vsv.\end{cases}\]
We recover the presentation of Section \ref{sec:4.3.1}, summarized in the diagram
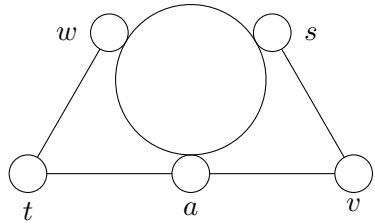
\spa
\begin{center}
\begin{tikzpicture}[scale=0.25]
\draw (0,0) circle (4);

\draw (30:5) circle (1);
\node[right] at (24:6) {$s$};
\node[left] at (156:6) {$w$};
\node[below] at (270:6) {$a$};
\node[below] at (-8.66,-5.96) {$t$};
\node[below] at (8.66,-5.8) {$v$};

\draw (150:5) circle (1);
\draw (270:5) circle (1);
\draw (-8.66,-5) circle (1);
\draw (8.66,-5) circle (1);
\draw (-7.66,-5)--(-1,-5);
\draw (1,-5)--(7.66,-5);
\draw (-8.16,-4.134)--(-4.83,1.634);
\draw (8.16,-4.134)--(4.83,1.634);
\end{tikzpicture}\end{center}
Again, we know that neither the monoid given by this presentation, neither the submonoid of $L_{u_2}^+$ generated by $a,s,t,v,w$ are Garside monoids.

The morphisms $\varphi_{2,1}$ and $\varphi_{2,1}$ are given by
\[\varphi_{2,1}:\begin{cases} s\mapsto s, & a\mapsto u^w,\\ t \mapsto t, & b\mapsto s^v, \\ v\mapsto v, & c\mapsto u^{wv}, \\ w\mapsto w,\end{cases}\text{~~and~~}\varphi_{1,2}:\begin{cases} s\mapsto s,\\ t\mapsto t,\\u\mapsto a^s,\\ v\mapsto b,\\ w\mapsto w. \end{cases}\]

Maximal roots of the full-twist are given by
\[\begin{array}{r||c|c}\text{Maximal regular number }d & 20 & 24  \\\hline d\text{-th root of }\Delta^{60}(u_2)& stwavs& stwav\end{array}\]

In particular, we get that $(stwav)^6=\Delta^{15}(u_2)$ generates $Z(\Bb_{31}(u_2,u_2))$.

\subsubsection{Orbit of $u_3$} The object $u_3$ admits $7$ atomic loops $d,e,f,t,u,v,w$, with relations as follows:
\[\begin{cases} ue=ef=fu,~wt=tf=fw, \\
df=fd,~wv=vw,~tv=vt,~vf=fv,\\
udw=dwu=wud,~dte=ted=edt,\\
utu=tut,~uvu=vuv,~dvd=vdv,~eve=vev,\\ wtud=edwt. 
\end{cases}\]
The last line of relations can be omitted, as it is implied by the others. The monoid $H_{u_3}^+$ given by this presentation is not cancellative: we have $tudt\neq udtu$ and $ wtudt= wudtu$ in $H_{u_3}^+$. Thus $H_{u_3}^+$ cannot be a Garside monoid.

The submonoid $L_{u_3}^+$ of $\Cc_{31}(u_3,u_3)$ generated by the atomic loops is cancellative, but it does not admit right-lcms: we know that $utu=tut$ is the shortest common multiple of $t$ and $u$. If right-lcms exists in $L_{u_3}^+$, then $utu$ must be the right-lcm of $u$ and $t$. We would then have
\[tut\preceq tudt\Rightarrow t\preceq dt,\]
which does not hold in $L_{u_3}^+$ (it does not even hold in $\Cc_{31}(u_3,u_3)$).

The relations defining $H_{u_3}^+$ give in particular $f:=w^t$ and $e=f^u=w^{tu}$ in $H_{u_3}$. By deleting these generators, we get that $\Bb_{31}(u_{3},u_{3})$ is generated by $d,t,u,v,w$, with relations as follows:
\[\begin{cases}wv=vw,~tv=vt,\\
udw=wud=dwu,\\
uvu=vuv,~tut=utu,~dvd=vdv,~wtw=twt,\\
udtu=tudt.
\end{cases}\]
The monoid given by this presentation is not cancellative: we have $tdwt\neq wtdw$ while $tdwtut=wtdwut$. The morphisms $\varphi_{3,1}$ and $\varphi_{1,3}$ are given by
\[\varphi_{3,1}:\begin{cases} t\mapsto t, & d\mapsto s^u,\\ u\mapsto u,&e\mapsto w^{tu}, \\ v\mapsto v,&f\mapsto w^t, \\ w\mapsto w, \end{cases}\text{~~and~~}\varphi_{1,3}:\begin{cases} s\mapsto d^w, \\t \mapsto t, \\ u\mapsto u,\\v\mapsto v, \\ w\mapsto w.\end{cases}\]

Maximal roots of the full-twist are given by
\[\begin{array}{r||c|c}\text{Maximal regular number }d & 20 & 24  \\\hline d\text{-th root of }\Delta^{60}(u_3)& wtudvw& tudvw\end{array}\]

In particular, we get that $(tudvw)^6=\Delta^{15}(u_3)$ generates $Z(\Bb_{31}(u_3,u_3))$.

\subsubsection{Orbit of $u_4$} The object $u_4$ admits $12$ atomic loops $g,h,k,l,m,n,o,p,s,t,u,v$, with relations as follows:
\[\begin{cases} 
gk=hg=kh,~gs=lg=sl,~gn=tg=nt,~gp=vg=pv,~ht=mh=tm,\\
kn=mk=nm,~lo=tl=ot,~so=ns=on,~ut=tp=pu,~uv=nu=vn,\\
gm=mg,~go=og,~hn=nh,~st=ts,~tv=vt,~np=pn,\\
hus=shu=ush,~hvo=ohv=voh,~kpo=okp=pok,\\
lmv=mvl=vlm,~smp=mps=psm,\\
hph=php,~svs=vsv,~mom=omo,\\
gnp=utv,~htvl=lomv,~khnps=lgomp,~knps=somp,\\
khnps=utvsm,~khn=tgm,~khpo=vgok,~lgomp=utvsm,\\
lgo=nst,~lgmp=pvsm,~usht=mhps,~uvsh=ohuv.
\end{cases}\]
The last three lines of relations can be omitted, as they are implied by the others. The monoid $H_{u_4}^+$ given by this presentation is not cancellative: we have $hpsh\neq pshp$ and $mhpsh=mpshp$ in $H_{u_4}^+$. Thus $H_{u_4}^+$ cannot be a Garside monoid.

The submonoid $L_{u_4}^+$ of $\Cc_{31}(u_4,u_4)$ generated by the atomic loops is cancellative, but it does not admit right-lcms: we know that $hph=php$ is the shortest common multiple of $h$ and $p$. If right-lcms exists in $L_{u_4}^+$, then $hph$ must be the right-lcm of $h$ and $p$. We would then have
\[hph\preceq hpsh\Rightarrow h\preceq sh,\]
which does not hold in $L_{u_4}^+$ (it does not even hold in $\Cc_{31}(u_4,u_4)$).The relations defining $H_{u_4}^+$ give in particular 
\[g=n^t=u^{vt},~k=h^g=h^{(u^{vt})},~l=g^s=u^{tvs},~m=h^t,~n=u^v,~o=n^s=u^{vs},~p=u^t\]
in $H_{u_4}$. By deleting these generators, we get that $\Bb_{31}(u_{4},u_{4})$ is generated by $h,s,t,u,v$, with relations as follows:
\[\begin{cases}vt=tv,~st=ts,\\
ush=shu=hus,\\
svs=vsv,~vuv=uvu,~utu=tut,~tht=hth,\\
shvs=vshv.
\end{cases}\]
The monoid given by this presentation does not admit right-lcms. We have $shvs=vshv$ and $vsv=svs$ are common right-multiples of $v$ and $s$, but their longest common divisor is $vs$,which is not a common righ-multiple of $s$ and $v$. This also proves that the submonoid of $\Cc_{31}(u_4,u_4)$ generated by $h,s,t,u,v$ does not admit right-lcms. The morphisms $\varphi_{4,1}$ and $\varphi_{1,4}$ are given by
\[\varphi_{3,1}:\begin{cases} g\mapsto u^{vt}, &o\mapsto u^{vs},\\h\mapsto w^s,&p\mapsto u^t,\\ k\mapsto u^{wtuwvsuvt},&s\mapsto s, \\ l\mapsto u^{tvs},&t\mapsto t,\\ m\mapsto w^{st},&u\mapsto u,\\n\mapsto u^v, &v\mapsto v,\end{cases}\text{~~and~~}\varphi_{1,3}:\begin{cases} s\mapsto s, \\t \mapsto t, \\ u\mapsto u,\\v\mapsto v, \\ w\mapsto h^u.\end{cases}\]
Maximal roots of the full-twist are given by
\[\begin{array}{r||c|c}\text{Maximal regular number }d & 20 & 24  \\\hline d\text{-th root of }\Delta^{60}(u_4)& tuvshv&tuvsh\end{array}\]
In particular, we get that $(tuvsh)^6=\Delta^{15}(u_4)$ generates $Z(\Bb_{31}(u_4,u_4))$. We also have $((tuvsh)^t)^2=(uvsht)^2=\Delta^5(u_4)$ in this case.

\subsubsection{Orbit of $u_5$} The object $u_5$ admits $10$ atomic loops $b,f,g,n,p,s,t,u,v,w$, with relations as follows:
\[\begin{cases}
ut=pu=tp,~uv=nu=vn,~gp=pv=vg,\\
gn=tg=nt,~sv=vb=bs,~wt=fw=tf,\\
st=ts,~wv=vw,~fv=vf,~pn=np,~tv=vt,~tb=bt,\\
uws=suw=wsu,~gfb=fbg=bgf,~spf=fsp=pfs,~wbn=nwb=bnw,\\
ufu=fuf,~gsg=sgs,~suw=wsu,~sns=nsn,~fnf=nfn\\
utv=gpn,~uwsv=bsnw,~gnfb=wtbg,~wsut=pufs,~pfsv=bsgf.
\end{cases}\]
The last line of relations can be omitted, as it is implied by the others. The monoid $H_{u_5}^+$ given by this presentation is not cancellative: we have $ufsu\neq fsuf$ and $pufsu=pfsuf$ in $H_{u_5}^+$. Thus $H_{u_5}^+$ cannot be a Garside monoid.

The submonoid $L_{u_5}^+$ of $\Cc_{31}(u_5,u_5)$ generated by the atomic loops is cancellative, but it does not admit right-lcms: we know that $ufu=fuf$ is the shortest common multiple of $f$ and $u$. If right-lcms exists in $L_{u_5}^+$, then $ufu$ must be the right-lcm of $u$ and $f$. We would then have
\[ufu\preceq ufsu\Rightarrow u\preceq su,\]
which does not hold in $L_{u_5}^+$ (it does not even hold in $\Cc_{31}(u_5,u_5)$).

The relations defining $H_{u_5}^+$ give in particular 
\[b=s^v,~f=w^t,~g=u^{vt},~n=u^v,~p=u^t\]
in $H_{u_5}$. By deleting these generators, we get that $\Bb_{31}(u_{5},u_{5})$ is generated by $s,t,u,v,w$, with the same relations as in $H_{u_1}^+$: we recover once again the obtained in Section \ref{sec:4.3.1}. The morphisms $\varphi_{5,1}$ and $\varphi_{1,5}$ are given by
\[\varphi_{5,1}:\begin{cases} s\mapsto s,&b\mapsto s^v,\\ t\mapsto t, &f\mapsto w^t, \\ u\mapsto u,&g\mapsto u^{vt}, \\ v\mapsto v,&n\mapsto u^v, \\ w\mapsto w, &p\mapsto u^t,\end{cases}\text{~~and~~}\varphi_{1,5}:\begin{cases} s\mapsto s, \\t \mapsto t, \\ u\mapsto u,\\v\mapsto v, \\ w\mapsto w.\end{cases}\]
Maximal roots of the full-twist are given by the same expressions as for the object $u_1$.

\subsubsection{Orbit of $u_6$} The object $u_6$ admits $10$ atomic loops $b,f,g,n,o,p,q,r,s,v$, with \newline relations as follows:
\[\begin{cases}vg=pv=gp,~qn=fq=nf,~qr=gq=rg,~so=ns=on,~vb=sv=bs,\\
vf=fv,~qb=bq,~sr=rs,~np=pn,~og=go,\\
vro=ovr=rov,~qop=pqo=opq,~spf=fsp=pfs,\\
bgf=fbg=gfb,~brn=nbr=rnb,\\
vqv=qvq,~vnv=nvn,~qsq=sqs,~sgs=gsg,~ngn=gng,~pqo=opq,\\
bsgf=pvfb,~bsro=onvr,~fqsp=onpf,~fqbr=rgnb,~pvqo=rovg.
\end{cases}\]
The last line of relations can be omitted, as it is implied by the others. The monoid $H_{u_6}^+$ given by this presentation is not cancellative: we have $nvrn\neq vrnv$ and $onvrn=ovrnv$ in $H_{u_6}^+$. Thus $H_{u_6}^+$ cannot be a Garside monoid.

The submonoid $L_{u_6}^+$ of $\Cc_{31}(u_6,u_6)$ generated by the atomic loops is cancellative, but it does not admit right-lcms: we know that $nvn=vnv$ is the shortest common multiple of $n$ and $v$. If right-lcms exists in $L_{u_6}^+$, then $nvn$ must be the right-lcm of $n$ and $v$. We would then have
\[nvn\preceq nvrn\Rightarrow n\preceq rn,\]
which does not hold in $L_{u_6}^+$ (it does not even hold in $\Cc_{31}(u_6,u_6)$).

The relations defining $H_{u_6}^+$ give in particular 
\[b=s^v,~f=q^n,~o=n^s,~p=v^g,~r=q^g\]
in $H_{u_6}$. By deleting these generators, we get that $\Bb_{31}(u_{6},u_{6})$ is generated by $g,n,q,s,v$, with relations as follows:
\[\begin{cases}nsn=sns,~vgv=gvg,~vsv=svs,~qnq=nqn,~vnv=nvn,\\
qgq=gqg,~ngn=gng,~qsq=sqs,~sgs=gsg,~vqv=qvq,\\
gnvg=vgnv=nvgn,~gqsg=qsgq=sgqs,~nsgn=gnsg=sgns,\\
vqsv=svqs=qsvq,~qnvq=nvqn=vqnv,~vgqnsv=svgqns,\\
gqnsvgs=ngqnsvg.
\end{cases}\]
The monoid given by this presentation does not admit right-lcms. We have that $vgqnsv=svgqns$ and $vsv=svs$ are common right multiples of $s$ and $v$, but their longest common left-divisor is $sv$, which is not a common right-multiple of $s$ and $v$. This also proves that the submonoid of $\Cc_{31}(u_6,u_6)$ generated by $g,n,q,s,v$ does not admit right-lcms. The morphisms $\varphi_{6,1}$ and $\varphi_{1,6}$ are given by
\[\varphi_{6,1}:\begin{cases} b\mapsto s^v, & p\mapsto u^t, \\ f\mapsto w^t, & q\mapsto t^{uwtv}, \\ g\mapsto u^{vt}, & r\mapsto u^{stuwtv}, \\ n\mapsto u^v, & s\mapsto s,\\ o\mapsto u^{vs}, &v\mapsto v,\end{cases}\text{~~and~~}\varphi_{1,6}:\begin{cases} s\mapsto s, \\t \mapsto p^{vn}, \\ u\mapsto v^n,\\v\mapsto v,  \\ w\mapsto p^{vnf}.\end{cases}\]

Maximal roots of the full-twist are given by
\[\begin{array}{r||c|c}\text{Maximal regular number }d & 20 & 24  \\\hline d\text{-th root of }\Delta^{60}(u_6)& svgqns& svgqn\end{array}\]

In particular, we get that $(svgqn)^6=\Delta^{15}(u_6)$ generates $Z(\Bb_{31}(u_6,u_6))$. We also have $svgqns=\Delta^3(u_6)$ in this case.

\subsubsection{Orbit of $u_7$} \label{sec:u7}The object $u_7$ admits $7$ atomic loops $g,k,m,n,o,p,s$, with relations as follows:
\[\begin{cases}on=so=ns,~kn=nm=mk,\\
go=og,~gm=mg,~pn=np,\\
pok=okp=kpo,~psm=smp=mps,\\
omo=mom,~gpg=pgp,~gkg=kgk,~gsg=sgs,~gng=ngn,\\
kpon=somp.
\end{cases}\]
The last line of relations can be omitted, as it is implied by the others. The monoid $H_{u_7}^+$ given by this presentation is not cancellative: we have $ompo\neq mpom$ and $sompo=smpom$ in $H_{u_7}^+$. Thus $H_{u_7}^+$ cannot be a Garside monoid.

The submonoid $L_{u_7}^+$ of $\Cc_{31}(u_7,u_7)$ generated by the atomic loops is cancellative, but it does not admit right-lcms: we know that $omo=mom$ is the shortest common multiple of $m$ and $o$. If right-lcms exists in $L_{u_7}^+$, then $omo$ must be the right-lcm of $m$ and $o$. We would then have
\[omo\preceq ompo\Rightarrow o\preceq po,\]
which does not hold in $L_{u_7}^+$ (it does not even hold in $\Cc_{31}(u_7,u_7)$). The relations defining $H_{u_7}^+$ give in particular $o=n^s$ and $k=n^m$ in $H_{u_7}$. By deleting these generators, we get that $\Bb_{31}(u_{7},u_{7})$ is generated by $g,m,n,p,s$, with relations as follows:
\[\begin{cases}gm=mg,~pn=np,\\
smp=psm=mps,\\
sns=nsn,~pgp=gpg,~gsg=sgs,~mnm=nmn,~gng=ngn\\
sgns=nsgn=gnsg.
\end{cases}\]
The monoid given by this presentation does not admit right-lcms. We have that $sns=nsn$ and $sgns=nsgn$ are common right multiples of $s$ and $n$, but their longest common left-divisor is $ns$, which is not a common right-multiple of $s$ and $n$. This also proves that the submonoid of $\Cc_{31}(u_7,u_7)$ generated by $g,m,n,p,s$ does not admit right-lcms. The morphisms $\varphi_{7,1}$ and $\varphi_{1,7}$ are given by
\[\varphi_{7,1}:\begin{cases} g\mapsto u^{vt}, &p\mapsto u^t,\\ k\mapsto u^{wtuwvsuvt} ,&s\mapsto s,\\ m\mapsto w^{st}, \\ n\mapsto u^v,\\ o\mapsto u^{vs}, \end{cases}\text{~~and~~}\varphi_{1,7}:\begin{cases} s\mapsto s, \\t \mapsto  g^n,\\ u\mapsto p^{gn},\\v\mapsto p^g, \\ w\mapsto p^{gkpn}. \end{cases}\]
Maximal roots of the full-twist are given by
\[\begin{array}{r||c|c}\text{Maximal regular number }d & 20 & 24  \\\hline d\text{-th root of }\Delta^{60}(u_7)& gnmpsm& gnmps\end{array}\]
In particular, we get that $(gnmps)^6=\Delta^{15}(u_7)$ generates $Z(\Bb_{31}(u_7,u_7))$.

\subsubsection{Orbit of $u_8$} The object $u_8$ admits $6$ atomic loops $g,m,n,p,s,t$, with relations as follows:
\[\begin{cases}tg=nt=gn,\\
ts=st,~pn=np,~gm=mg,\\
psm=smp=mps,\\
tpt=ptp,~tmt=mtm,~pgp=gpg,~sns=nsn,~sgs=gsg,~nmn=mnm.
\end{cases}\]
The monoid $H_{u_8}^+$ given by this presentation is not cancellative: we have $psgpsg\neq gpsgps$ and $ mpsgpsg=mgpsgps$ in $H_{u_8}^+$. Thus $H_{u_8}^+$ cannot be a Garside monoid.

The submonoid $L_{u_8}^+$ of $\Cc_{31}(u_8,u_8)$ generated by the atomic loops is cancellative, but it does not admit right-lcms: we know that $pgp=gpg$ is the shortest common multiple of $p$ and $g$. If right-lcms exists in $L_{u_8}^+$, then $omo$ must be the right-lcm of $m$ and $o$. We would then have
\[gpg\preceq gpsgps\Rightarrow g\preceq sgps,\]
which does not hold in $L_{u_8}^+$ (it does even hold in $\Cc_{31}(u_8,u_8)$, but $g^{-1}sgps$ is not generated positively by atomic loops). The relations defining $H_{u_8}^+$ give in particular $t=g^n$ in $H_{u_8}$. By deleting this generator, we get that $\Bb_{31}(u_{8},u_{8})$ is generated by $g,m,n,p,s$, with the same relations as in the end of Section \ref{sec:u7}. The morphisms $\varphi_{8,1}$ and $\varphi_{1,8}$ are given by
\[\varphi_{8,1}:\begin{cases} g\mapsto u^{vt},\\ m\mapsto w^{st}, \\ n\mapsto u^v, \\p\mapsto u^t, \\ s\mapsto s,\\t\mapsto t, \end{cases}\text{~~and~~}\varphi_{1,8}:\begin{cases} s\mapsto s, \\t \mapsto  t,\\ u\mapsto t^p,\\v\mapsto n^{tp},\\ w\mapsto  m^{tmp}.\end{cases}\]
Maximal roots of the full-twist are given by the same expressions as for the object $u_7$. Furthermore, we also have $(gnmps)^2=\Delta^5(u_8)$ in this case.

\begin{rem}
As was pointed out to us by Jean Michel, the presentations with $5$ generators that we give here are related by the Hurwitz action of the usual braid groups on the words giving 24-roots of the full-twist, as in \cite[Section 6]{hurwitz_presentation}. By identifying the atomic loops we consider with their image in $\Bb_{31}(u_1,u_1)$, we obtain for instance that $(s,t,u,v,w)$, $(s,t,w,a,v), (t,u,d,v,w)$ and $(t,u,v,s,h)$ lie in the same Hurwitz orbit. The words $(s,u,v,w,t)$ and $(s,v,g,q,n)$ also lie in the same Hurwitz orbit, as well as $(u,v,w,t,s)$ and $(g,n,m,p,s)$.
\end{rem}

\appendix

\section{The Reidemeister-Schreier method for groupoids}
\subsection{Presentations of categories and groupoids} \label{sec:presentations_categories}

\begin{definition}\cite[Definition II.1.4 and Definition II.1.32]{ddgkm} An \nit{oriented graph} (or precategory) is a pair of sets $(\Oo,\Ss)$ endowed with two maps $s,t:\Ss\to \Oo$. The elements of $\Oo$ are called \nit{objects} and those of $\Ss$ are called elements (or morphisms, or arrows). The maps $s$ and $t$ are called \nit{source} and \nit{target}, respectively.

A \nit{morphism} between two oriented graphs $(\Oo,\Ss)$ and $(\Oo',\Ss')$ is given by two maps $\phi_0:\Oo\to \Oo',$ $\phi_1:\Ss\to \Ss'$ which preserve the source and target:
\[\forall f\in \Ss, s(\phi_1(f))=\phi_0(s(f))\text{~~and~~} t(\phi_1(f))=\phi_0(t(f)).\]
\end{definition}

Note that we make no assumption regarding the existence of loops or multiple arrows with the same source and target. In practice, we often amalgamate an oriented graph $(\Oo,\Ss)$ with its set of arrows $\Ss$. The set of objects will then be denoted by $\Ob(\Ss)$ like for categories. We will sometimes denote by $\Ss(u,-)$ (resp. $\Ss(-,v)$, $\Ss(u,v)$) the set of arrows in $\Ss$ with source $u$ (resp. with target $v$, with source $u$ and target $v$).

\begin{definition}\cite[Definition II.1.28]{ddgkm} Let $\Ss$ be an oriented graph. For $u,v\in \Ob(\Ss)$, a \nit{path} of length $p\geqslant 1$ in $\Ss$ from $u$ to $v$ is a finite sequence $(g_1,\ldots,g_p)$ of elements of $\Ss$ such that
\[s(g_1)=u,~t(g_p)=v\text{~~and~~} t(g_i)=s(g_{i+1}), \forall i\in \intv{1,p-1}.\]
For $u\in \Ob(\Cc)$, we also define an \nit{empty path} from $u$ to itself, denoted by $1_u$, and of length $0$ by definition.
\end{definition}
\begin{definition}\cite[Definition II.1.28 and Proposition 1.33]{ddgkm} Let $\Ss$ be an oriented graph. The \nit{free category} on $\Ss$, denoted by $\Ss^*$, is defined as follows:
\begin{itemize}
\item The objects are the objects of $\Ss$.
\item The morphisms $\Ss^*(u,v)$ are the paths from $u$ to $v$ in $\Ss$.
\item Composition is given by concatenation of paths.
\item The identity of some object $u$ is the empty path $1_u$.
\end{itemize}
\end{definition}

This category is free in the following sense: Let $\Ss$ be an oriented graph, and let $\Cc$ be a category. Any morphism of oriented graphs $\phi:\Ss\to \Cc$ induces a unique functor $\Ss^*\to \Cc$, sending a path $(g_1,\ldots,g_p)$ to the composition $\phi_1(g_1)\cdots\phi_1(g_p)$ in $\Cc$. In practice, a path $(f_1,\ldots,f_p)$ will often be denoted as a formal composition $f_1\cdots f_p$ so that we have
\[\phi(f_1\cdots f_p)=\phi_1(f_1)\cdots\phi_1(f_p).\]

This convenient definition of free category allows for defining relations and presentations of categories. Recall that a \nit{congruence} on a category $\Cc$ is an equivalence relation $\equiv$ on $\Cc$ which is compatible with composition, that is, the conjunction of $f\equiv f'$ and $g\equiv g'$ implies $fg\equiv f'g'$ (if $fg$ and $f'g'$ are defined of course). If $\equiv$ is a congruence on a category $\Cc$, one can form the quotient category $\Cc/\equiv$. It has the same objects as $\Cc$, and its morphisms are $\equiv$-equivalence classes of morphisms in $\Cc$.

\begin{definition}\cite[Definition II.1.36 and Lemma II.1.37]{ddgkm} If $\Cc$ is a category, a \nit{relation} on $\Cc$ is a pair $(g,h)$ of morphisms in $\Cc$ sharing the same source and the same target. If $\Cc$ is a category, and $\Rr$ is a family of relations on $\Cc$, there exists a smallest congruence $\equiv_\Rr^+$ of $\Cc$ which includes $\Rr$.
\end{definition}

The congruence $\equiv_\Rr^+$ is the reflexive-transitive closure of 
\[\{(fgh,fg'h)~|~ (g,g')\in \Rr\text{ or }(g',g)\in \Rr\}.\]

For readability purposes, it is convenient to write a relation $(f,g)$ as an equality $f=g$ instead of a couple of paths. We use this convention from now on.

\begin{definition}\cite[Definition II.1.38]{ddgkm} A \nit{category presentation} is a pair $(\Ss,\Rr)$, where $\Ss$ is an oriented graph, and $\Rr$ is a set of relations on $\Ss^*$. We call $\Ss$ the \nit{generators} and $\Rr$ the \nit{relations}. \newline If $(\Ss,\Rr)$ is a category presentation, the quotient category $\Ss/\equiv_\Rr^+$ is denoted by $\langle \Ss~|~ \Rr\rangle^+$.
\end{definition}
\begin{rem}
If we consider a graph $\Ss$ with one object, we recover the classical notion of monoid presentation.
\end{rem}

Let $(\Ss,\Rr)$ be a categorical presentation, and let $\Cc$ be a category. Any morphisms of oriented graphs $\phi:\Ss\to \Cc$ such that
\[\forall f_1\cdots f_p=g_1\cdots g_q\in \Rr, \phi_1(f_1)\cdots\phi_1(f_p)=\phi_1(g_1)\cdots\phi_1(g_p)\in \Cc\]
induces a unique functor from the presented category $\langle \Ss~|~\Rr\rangle^+$ to $\Cc$. A first application of this property is the following lemma.

\begin{lem}\label{lem:homogeneous_presentation}
Let $\Cc:=\langle S~|~R\rangle^+$ be a presented category. If $\Rr$ consists of relations between paths of the same length in $\Ss^*$, then the category $\Cc$ is homogeneous. We then say that $\langle S~|~ R\rangle^+$ is a \nit{homogeneous presentation}.
\end{lem}
\begin{proof}
First, we note that the category $\Ss^*$ is homogeneous. The length functor is given by the length of the paths. By assumption, this functor induces a well-defined functor from $\Cc$ to $(\Z_{\geqslant 0},+)$, such that elements of $\Ss$ are sent to $1$. This functor is a length functor for $\Cc$ since it is generated by $\Ss$.
\end{proof}

The notion of categorical presentation is useful for defining enveloping groupoids: let $\Cc=\langle \Ss~|~\Rr\rangle^+$ be a presented category. We consider $\bbar{\Ss}$ a formal copy of $\Ss$, with source and target reversed. We consider the set $I(\Ss)$ of relations on $\Ss\sqcup \bbar{\Ss}$ defined by
\[\forall x\in \Ss, x\bbar{x}=1_{s(x)}\text{~and~} \bbar{x}x=1_{t(x)}.\]

\begin{lem}\cite[Definition II.3.3 and Proposition II.3.5]{ddgkm}\label{lem:enveloping_groupoid} The category $\Gg(\Cc):= \langle \Ss\sqcup\bbar{\Ss}~|~ \Rr\cup I(\Ss)\rangle^+$ is a groupoid, called the \nit{enveloping groupoid} of $\Cc$. The inclusion map $\Ss\hookrightarrow \Ss\cup \bbar{\Ss}$ induces a functor $\iota:\Cc\to \Gg(\Cc)$. Every functor $\phi:\Cc\to\Gg$ where $\Gg$ is a groupoid  induces a unique functor $\ttilde{\phi}:\Gg(\Cc)\to \Cc$ such that $\ttilde{\phi}\circ \iota=\phi$.
\end{lem}

This universal property of the enveloping groupoids ensures that it depends only on $\Cc$ and not on its presentation. By convention, if $\Cc=\langle \Ss~|~\Rr\rangle^+$ is a presented category, the presentation of the enveloping groupoid of $\Cc$ will be denoted by $\langle \Ss~|~\Rr\rangle$.

\begin{rem}
Every category $\Cc$ admits a standard presentation, where the generators are all the morphisms in $\Cc$, and the relations are all couples $(fg,h)$ for $f,g,h\in \Cc$ satisfying $fg=h$. Thus the enveloping groupoid can be defined for any category.
\end{rem}

\subsection{Schreier transversal and presentation}\label{app:reidemeisterschreier}
Let $\Gg=\langle \Ss~|~\Rr\rangle$ be a presented connected groupoid. We denote by $\Ff(\Ss):=\langle \Ss~|~\varnothing\rangle$ the free groupoid on the graph $\Ss$. The identity $\Ss\to \Ss$ induces a quotient map $\pphi:\Ff(\Ss)\to \Gg$. Just like for free groups, a morphism in $\Ff(\Ss)$ is represented by a unique \nit{reduced path}, that is a path containing no subpath of the form $ss^{-1}$ or $s^{-1}s$ for $s\in \Ss$.

\begin{definition}\label{def:schreiertransversal}
Let $u$ be an object of $\Ff(\Ss)$. A \nit{Schreier transversal} of $\Ff(\Ss)$ rooted in $u$ is a family of reduced paths $\{m_v\}_{v\in \Ob(\Ff(\Ss))}$ satisfying:
\begin{itemize}
\item For every object $v$ of $\Ss$, the path $m_v$ has source $u$ and target $v$.
\item The family $\{m_v\}$ is stable under prefix. It particular, $m_u=1_u$.
\end{itemize}
\end{definition}

\begin{rem}
Since $\Gg$ is connected, it is also the case of $\Ss$ and $\Ff(\Ss)$. In particular, a Schreier transversal of $\Ff(\Ss)$ rooted in $u$ exists for every object $u$ of $\Ff(\Ss)$.
\end{rem}

Let $u_0$ be an object of $\Ff(\Ss)$, and let $\{m_v\}$ be a Schreier transversal in $\Ff(\Ss)$ rooted in $u_0$. For $s\in \Ss(u,v)$, we define $\gamma(s):=m_u s m_v^{-1}\in \Ff(\Ss)(u_0,u_0)$. Let $S_1$ be the set of all elements $\gamma(s)\neq 1_{u_0}$ for $s\in \Ss$.

\begin{lem}\label{lem:lemme_de_schreier}
The group $\Gg(u_0,u_0)$ is generated by the $\pphi(\gamma(s))$ for $\gamma(s)\in S_1$.
\end{lem}
\begin{proof}
Let $g\in \Gg(u_0,u_0)$. Since $\Gg$ is generated by $\Ss$, we can write $g=s_1^{\epsi_1}\cdots s_k^{\epsi_k}$ with $s_i\in \Ss$ for $i\in \intv{1,k}$ and $\epsi_i\in \{\pm 1\}$ for $i\in \intv{1,k}$. We denote by $u_i$ the target of $s_i^{\epsi_i}$ for $i\in \intv{1,k-1}$. In $\Ff(\Ss)$ we have
\begin{align*}
s_1^{\epsi_1}\cdots s_k^{\epsi_k}&=s_1^{\epsi_1}m_{u_1}^{-1}m_{u_1}\cdots m_{u_{k-1}}^{-1} m_{u_{k-1}}s_k^{\epsi_k}\\
&=m_{u_0}s_1^{\epsi_1}m_{u_1}^{-1}m_{u_1}\cdots m_{u_{k-1}}^{-1} m_{u_{k-1}}s_k^{\epsi_k}m_{u_0}^{-1}\\
&=\gamma(s_1)^{\epsi_1}\cdots \gamma(s_k)^{\epsi_k}.
\end{align*}
Thus, we have $g=\pphi(\gamma(s_1))^{\epsi_1}\cdots \pphi(\gamma(s_k))^{\epsi_k}$ in $\Gg(u_0,u_0)$.
\end{proof}

As $\Gg$ is a groupoid, every relation defining $\Gg$ can be rewritten as a relation of the form
\[s_1^{\epsi_1}\cdots s_k^{\epsi_k}=1_u,\]
where $s_i\in \Ss$ and $\epsi_i\in \{\pm 1\}$ for $i\in \intv{1,k}$, and $u$ is the source of $s_1^{\epsi_1}$. 

\begin{prop}[\textbf{Reidemeister-Schreier method for groupoids}]\label{prop:reidemeister_schreier}
Let $S^*$ be a set of elements $\gamma(s)^*$ in one-to-one correspondence with those of $S_1$. Let also $R^*$ be the set of all relations
\[\gamma^*(s_1)^{\epsi_1}\cdots \gamma^*(s_k)^{\epsi_k}=1,\]
where $s_1^{\epsi_1}\cdots s_k^{\epsi_k}=1_u$ is in $\Rr$. The map $S^*\to \Gg(u_0,u_0)$ sending $\gamma(s)^*$ to $\pphi(\gamma(s))$ induces an isomorphism of groups between $G^*:=\langle S^*~|~ R^*\rangle$ and $\Gg(u_0,u_0)$.
\end{prop}
\begin{proof}
First, we prove that the map $\gamma(s)^*\mapsto \pphi(\gamma(s))$ is compatible with the set of relations $R^*$. Let $\gamma^*(s_1)^{\epsi_1}\cdots \gamma^*(s_k)^{\epsi_k}=1$ be in $R^*$. We have an equality $s_1^{\epsi_1}\cdots s_k^{\epsi_k}=1_u$ in $\Gg$, and we have
\[\pphi(\gamma(s_1))^{\epsi_1}\cdots \pphi(\gamma((s_k))^{\epsi_k}=\pphi(m_u)s_1^{\epsi_1}\cdots s_k^{\epsi_k}\pphi(m_u)^{-1}=1_{u_0}.\]
Let $\pi:G^*\to \Gg(u_0,u_0)$ be the morphism induced by $\gamma(s)^*\mapsto \pphi(\gamma(s))$. We know that $\pi$ is surjective by Lemma \ref{lem:lemme_de_schreier}.

Conversely, the map $\Ss\to S^*$ sending $s$ to $\gamma(s)^*$ induces a functor $\phi:\Ff(\Ss)\to G^*$. Let $v$ be an object of $\Ff(\Ss)$. We show that $\phi$ sends $m_v$ to $1$ by induction on the length of $m_v$ as an $\Ss$-path. First if $m_v=s_1^{\epsi_1}$ has length $1$, we have $\phi(s_1^{\epsi_1})=m_{u_0}m_{v}m_{v}^{-1}=1$. Now if $m_v=s_1^{\epsi_1}\cdots s_k^{\epsi_k}$ is a decomposition of $m_v$ on $\Ss$, we denote by $v'$ the source of $s_k^{\epsi_k}$. By definition of a Schreier transversal we have $m_{v'}=s_1^{\epsi_1}\cdots s_{k-1}^{\epsi_{k-1}}$ and $\phi(m_{v'})=1$ by induction hypothesis. As we also have $\phi(s_k^{\epsi_k})=m_{v'}s_k^{\epsi_k} m_{v}^{-1}=m_vm_v^{-1}=1$, we get $\phi(m_v)=1$.

By definition of the set $R^*$, $\phi$ induces a functor $\bbar{\phi}:\Gg\to G^*$. We call $\iota$ the restriction of this functor to $\Gg(u_0,u_0)$. Let $\gamma(s)^*$ be a generator of $G^*$, with $s\in \Ss(u,v)$. We have
\begin{align*}
\iota(\pi(\gamma(s)^*))&=\iota(\varphi(\gamma(s)))\\
&=\iota(\pphi(m_u) s \pphi(m_v)^{-1})\\
&=\bbar{\phi}(\pphi(m_u))\bbar{\phi}(s)\bbar{\phi}(\pphi(m_v))^{-1}\\
&=\phi(m_u)\gamma(s)^*\phi(m_v)^{-1}=\gamma(s)^*.
\end{align*}
So $\iota\circ\pi$ induces the identity on the generators of $G^*$. We get $\iota \circ \pi=1_{G^*}$ and $\pi$ is injective. 
\end{proof}

\begin{cor}
If $\Gg$ is a finitely presented groupoid, then for every object $u$ of $\Gg$, the group $\Gg(u,u)$ is finitely presented.
\end{cor}

\subsection{The particular case of a subgroup of a presented group} We briefly outline how our result on groupoids can be used to recover the classical result of Reidemeister and Schreier about presentation of subgroups.

\begin{definition}
Let $G$ be a group, with $H$ a subgroup. The \nit{groupoid of cosets} $G_H$ of $G$ and $H$ is the groupoid defined by
\begin{itemize}
\item The objects are the right-cosets of $H$ in $G$.
\item Morphisms between two cosets $Hg$ and $Hg'$ are elements $x$ of $G$ such that $Hgx=Hg'$. 
\item Composition is given by the product in $G$.
\end{itemize}
\end{definition}
The underlying oriented graph of the groupoid of cosets is simply the graph of the action of $G$ on the right-cosets of $H$ in $G$. For each coset $Hg$ and each element $x$ of $G$, we denote by $x[Hg]$ the unique morphism $Hg\to H{gx}$ corresponding to the action of the element $x$. 

By definition of $G_H$, there is a natural functor $\pi:G_H\to G$, sending a morphism $x[Hg]$ to $x$. For two cosets $Hg$ and $Hg'$, the set of morphisms between $Hg$ and $Hg'$ is given by
\[G_H(Hg,Hg')=\{x \in G~|~ Hgx=Hg'\}=g^{-1}Hg'.\]
In particular, we obtain that $G_H(H,H)=H$ and $H$ is an automorphism group in $G_H$.

Since the natural action of $G$ on $H\setminus G$ is transitive, the groupoid $G_H$ is connected. In particular we can apply the Reidemeister-Schreier method to obtain a presentation of $H$ from one of $G_H$.

Now, if $\langle X,R\rangle$ is a presentation of the group $G$. We denote by $F(X)$ the free group over the set $X$. We call $\Xx$ the subgraph of $G_H$ made of those morphisms in $G_H$ whose image under $\pi$ lies in $X$. The restriction of $\pi$ to $\Xx$ induces a functor $\ttilde{\pi}:\Ff(\Xx)\to F(X)$. The inclusion $\Xx\to G_H$ (resp. $X\to G$) induces a functor $\varphi:\Ff(\Xx)\to G_H$ (resp. $\phi:F(X)\to G$). We have the following commutative square
\[\xymatrix{\Ff(\Xx)\ar[r]^-{\ttilde{\pi}} \ar[d]_-{\varphi}& F(X)\ar[d]^-\phi \\ G_H \ar[r]_-{\pi} & G}\] 

Let $m:=x_1^{\epsi_1}\cdots x_k^{\epsi_k}$ be a word in $F(X)$, and let $Hg$ be a right-coset. There is a unique path in $\Ff(\Xx)$ which starts at $Hg$ and whose image under $\ttilde{\pi}$ is $m$. This path is given by
\[x_1[Hg]^{\epsi_1}x_2[Hgx_1^{\epsi_1}]^{\epsi_2}\cdots x_k[Hgx_1^{\epsi_1}\cdots x_{k-1}^{\epsi_{k-1}}]^{\epsi_k}.\]
It will be denoted by $m[Hg]$. Note that $\varphi(m[Hg])=\phi(m)[Hg]$. In particular we see that $\Xx$ generates $G_H$ since $X$ generates $G$.
\begin{lem}The groupoid of cosets $G_H$ admits the presentation $\langle \Xx~|~ \Rr\rangle$, where the set $\Rr$ consists of relations $r[Hg]=1_{Hg}$, where $r=1$ lies in $R$ and $Hg$ lies in $H\setminus G$.
\end{lem}
\begin{proof}Let $\Gg$ be the groupoid presented by $\Xx$ and $\Rr$. Since the graph $\Xx$ generates $G_H$, the natural functor $\Gg\to G_H$ is surjective on morphisms. Consider now two paths in $\Ff(\Xx)$
\[m:=x_1^{\epsi_1}x_2^{\epsi_2}\cdots x_k^{\epsi_k}\text{~~and~~}m':=y_1^{\eta_1}y_2^{\eta_2}\cdots y_m^{\eta_m}.\]
By definition, saying that these two paths induce the same morphism in $G_H$ amounts to saying that they share the same source and that the two words $\ttilde{\pi}(m)$ and $\ttilde{\pi}(m')$ induce the same element in $G$. 

Suppose that $m$ and $m'$ induce the same element in $G_H$. We have $\ttilde{\pi}(m)=\ttilde{\pi}(m')$ in $G$, and there is a finite sequence of words $m_1,\ldots,m_p$ in $F(X)$ such that 
\begin{itemize}
\item $m_1=\ttilde{\pi}(m),m_p=\ttilde{\pi}(m')$ 
\item For $i\in \intv{1,p-1}$, $m_i$ is equivalent to $m_{i+1}$ by the use of one relation in $R$.
\end{itemize}
Let $Hg$ be the common source of $m$ and $m'$ in $\Ff(\Xx)$. We consider the paths $m_i[Hg]$ in $\Ff(\Xx)$. We have $m_1[Hg]=m,m_p[Hg]=m'$. For $i\in \intv{1,p-1}$, the path $m_i[Hg]$ is equivalent to $m_{i+1}[Hg]$ by the use of one relation in $\Rr$ by definition.

Thus $m$ and $m'$ induce the same element in $\Gg$, and the natural functor $\Gg\to G_H$ is then injective.
\end{proof}

Let now $\ttilde{H}=\phi^{-1}(H)\subset F(X)$ be the preimage of $H$ under $\phi$. Recall from \cite[Proposition I.3.8]{reisch} that a Schreier transversal (in the classical sense) for $\ttilde{H}$ in $F(X)$ is a set of words $T$ such that
\begin{itemize}
\item The maps $t\mapsto Ht$ is a bijection between $T$ and $H\setminus G$.
\item The set $T$ is stable under prefix (in particular, the empty word lies in $T$).
\end{itemize}

Let $T$ be a Schreier transversal in the classical sense for $\ttilde{H}$ and $F(X)$. The set of paths $\{t[H]\}_{t\in T}$ is a Schreier transversal rooted in $H$ in $\Ff(\Xx)$ in the sense of Definition \ref{def:schreiertransversal}. Conversely, if $\{m_{Hg}\}_{Hg\in \Ob(\Ff(\Xx))}$ is a Schreier transversal rooted in $H$ in the sense of Definition \ref{def:schreiertransversal}, then the set of words $\{\ttilde{\pi}(m_{Hg})\}_{Hg\in H\setminus G}$ is a Schreier transversal for $\ttilde{H}$ and $F(X)$ in the classical sense.

Let $T$ be a Schreier transversal for $\ttilde{H}$ in $F(X)$, and let $\{m_{Hg}\}_{Hg\in \Ob(\Ff(\Xx))}$ be the associated Schreier transversal of $\Ff(\Xx)$ rooted in $H$. For $g\in G$, let $\bbar{g}$ denote the element of $T$ such that $Hg=H\bbar{g}$.

By definition, we have $\Xx=\{x[Hg]~|~x\in X,Hg\in \Ob(G_H)\}$. The set of generators of $G_H(H,H)=H$ we obtain by our method is
\begin{align*}
&\{\gamma(x[Hg])=m_{Hg}xm_{Hgx}^{-1}~|~x\in X,Hg\in \Ob(G_H)\}\\
=&\{\gamma(x[Ht])=m_{Ht}xm_{Htx}^{-1}~|~x\in X,t \in T\}\\
=&\{tx(\bbar{tx})^{-1}~|~x\in X,t\in T\},
\end{align*}
which is the same set as given in \cite[Proposition 4.1]{reisch}. Following \cite[Proposition 4.1]{reisch}, we denote $\gamma(t,x):=tx(\bbar{tx})^{-1}$. This element is equal to what we denoted by $\gamma(x[Ht])$.

Now for the relators. Let $r=1$ be a relation of $G$. It induces the following family of relation on $G_H$:
\[\{r[Hg]=1_{Hg}~|~Hg\in \Ob(G_H)\}=\{r[Ht]=1_{Ht}~|~ t\in T\}.\]
Each relation $r[Ht]=1_{Ht}$ induces a relation on $G_H(H,H)$, given by $m_{Ht}r[Ht]m_{Ht}^{-1}=1$. If $r=y_1y_2y_3\cdots y_k$ is expressed as a word in $X\cup X^{-1}$, we have
\begin{align*}
m_{Ht}r[Ht]m_{Ht}^{-1}&=\gamma(y_1[Ht])\gamma(y_2[Hty_1])\cdots \gamma(y_k[Hty_1\cdots y_{k-1}])\\
&=\gamma(t,y_1)\gamma(\bbar{ty_1},y_2)\cdots \gamma(\bbar{ty_1\cdots y_{k-1}},y_k)\\
&=\gamma(1,t)\gamma(t,y_1)\gamma(\bbar{ty_1},y_2)\cdots \gamma(\bbar{ty_1\cdots y_{k-1}},y_k)\gamma(\bbar{ty_1\cdots y_k},t^{-1})
\end{align*}
since both $\gamma(1,t)$ and $\gamma(\bbar{ty_1\cdots y_k},t^{-1})$ are trivial. Thus the relations we obtain with our method are the same as those given in \cite[Proposition 4.1]{reisch}, and Proposition \ref{prop:reidemeister_schreier} applied to a groupoid of cosets gives a new proof of \cite[Proposition 4.1]{reisch}.

\printbibliography
\end{document}